\documentclass[11pt,leqno]{amsart}
\usepackage{blindtext}
\usepackage{fancyhdr}
\usepackage{ifpdf}
\usepackage{curves}
\usepackage{tikz}
\usepackage{subcaption}
\usepackage{placeins}
\usepackage{tabularx}
\usepackage{url}

\usepackage{amsfonts}
\usepackage{amsthm}
\usepackage{amsmath}
\usepackage{amscd}
\usepackage{amssymb}
\usepackage{graphicx}
\usepackage{rotating}
\usepackage[small,nohug,heads=LaTeX]{diagrams}
\diagramstyle[labelstyle=\scriptstyle]
\usepackage{ulem}
\usepackage{textcomp}
\usepackage{appendix}
\usepackage{chngcntr}
\usepackage{etoolbox}
\usepackage{lipsum}
\usepackage{tikz}
\usepackage{syntonly}
\usepackage[pagebackref]{hyperref}

\usepackage{leftidx}

\usepackage[utf8]{inputenc}
\usepackage[english, russian]{babel}

\usetikzlibrary{matrix,arrows,decorations.pathmorphing}

\AtBeginEnvironment{subappendices}{%
\chapter*{Appendix}
\addcontentsline{toc}{chapter}{Appendices}
\counterwithin{figure}{section}
\counterwithin{table}{section}
}

\newcommand{\db}[1]{\operatorname{D}^{\operatorname{b}}(#1)}

\newcommand{\mcf}{\mathcal{F}}
\newcommand{\mco}{\mathcal{O}}
\newcommand{\mcc}{\mathcal{C}}

\def\CC {{\mathbb C}}     
\def\FF {{\mathbb F}}     
\def\HH {{\mathbb H}}     
\def\NN {{\mathbb N}}     
\def\PP {{\mathbb P}}     
\def\QQ {{\mathbb Q}}     
\def\RR {{\mathbb R}}     
\def\SS {{\mathbb S}}     
\def\ZZ {{\mathbb Z}}     

\newcommand{\op}[1]{\operatorname{#1}}

\def\ring#1{\ifmmode \mathaccent'027 #1\else \rm\accent'027 #1\fi}

\newcommand{\re}{{\mathrm e}}

\newcommand{\ri}{{\mathrm i}}

\newcommand{\lift}[2]{%
\setlength{\unitlength}{1pt}
\begin{picture}(0,0)(0,0)
\put(0,{#1}){\makebox(0,0)[b]{${#2}$}}
\end{picture}
}
\newcommand{\lowerarrow}[1]{%
\setlength{\unitlength}{0.03\DiagramCellWidth}
\begin{picture}(0,0)(0,0)
\qbezier(-28,-4)(0,-18)(28,-4)
\put(0,-14){\makebox(0,0)[t]{$\scriptstyle {#1}$}}
\put(28.6,-3.7){\vector(2,1){0}}
\end{picture}
}
\newcommand{\upperarrow}[1]{%
\setlength{\unitlength}{0.03\DiagramCellWidth}
\begin{picture}(0,0)(0,0)
\qbezier(-28,11)(0,25)(28,11)
\put(0,21){\makebox(0,0)[b]{$\scriptstyle {#1}$}}
\put(28.6,10.7){\vector(2,-1){0}}
\end{picture}
}

\def\ol  {\overline}
\def\ul  {\underline}

\def\wt  {\widetilde}

\def\cl  {\mathfrak{cl}}

\def\im {\mathfrak{im}}
\def\pr {\mathfrak{pr}}
\def\mc {\mathcal}
\def\mk {\mathfrak}
\def\ms {\mathsf}
\def\Hom {\mathrm{Hom}}

\def\st {\mathrm{Stab}}
\def \Aut {\mathrm{Aut}}
\def \Arg {\mathrm{Arg}}

\def\vol {\mathrm{vol}}

\def\tr {\mathrm{tr}}

\def \bd {\begin{diagram}}
\def \ed {\end{diagram}}
\def\be  {\begin{eqnarray}}
\def\ee  {\end{eqnarray}}
\def\ben {\begin{eqnarray*}}
\def\een {\end{eqnarray*}}

\def\bpr {\begin{proof}[Proof]}
\def\epr {\end{proof}}
\def\bsp {\begin{split}}
\def\esp {\end{split}}
\def\bprr {\begin{proof}[solution]}
\def\bpru {\begin{proof}[hint]}
\def\bpro {\begin{proof}[answer]}
\def\bcd {\begin{CD}}
\def\ecd {\end{CD}}

\def\cl {\rm Cl}

\newcommand{\abs}[1]{\left\vert#1\right\vert}
\newcommand{\scal}[1]{\left\langle#1\right\rangle}
\newcommand{\norm}[1]{\left\Vert#1\right\Vert}

\newtheorem{theorem}{Theorem}[section]
\newtheorem{lemma}[theorem]{Lemma}
\newtheorem{prop}[theorem]{Proposition}
\newtheorem{coro}[theorem]{Corollary}
\newtheorem{remark}[theorem]{Remark}
\newtheorem{df}[theorem]{Definition}

\newtheorem{ex}[theorem]{Example}
\newtheorem{conj}[theorem]{Conjecture}
\newtheorem{quest}[theorem]{Question}

\renewcommand{\contentsname}{Contents}
\renewcommand{\refname}{References}
\begin{document}
\renewcommand{\contentsname}{Contents}
\renewcommand{\refname}{References}
 
 \title[Some new categorical invariants]%
 {Some new categorical invariants}
 
\author{George Dimitrov}
\address[Dimitrov]{Universit\"at Wien\\
 Oskar-Morgenstern-Platz 1, 1090 Wien \\
 \"Osterreich
}
\email{george.dimitrov@univie.ac.at}

\author{Ludmil Katzarkov}
\address[Katzarkov]{Universit\"at Wien\\
Oskar-Morgenstern-Platz 1, 1090 Wien\\
\"Osterreich \\
 National Research University Higher School of Economics, Russian Federation  }
\email{lkatzarkov@gmail.com}

 \dedicatory{Dedicated to the memory of our friends Vasil Tsanov and Vladimir  Voevodskij}
\renewcommand{\abstractname}{Abstract}
\begin{abstract}

In this paper we introduce several new categorical notions and give many examples.

We   prove that   the moduli space of stability conditions on the derived category of representations 
of $K(l)$, the $l$-Kronecker quiver, is biholomorphic to ${\mathbb C} \times \mathcal H$ for $l\geq 3$. 
 This produces   an example  of semi-orthogonal decomposition, SOD, $\mathcal T =\langle \mathcal T_1, \mathcal T_2 \rangle$, where ${\rm Stab}(\mathcal T)$ is not biholomorphic to  ${\rm Stab}(\mathcal T_1)\times {\rm Stab}(\mathcal T_2)$  (whereas  ${\rm Stab}(\mathcal T_1 \oplus \mathcal T_2) \cong {\rm Stab}(\mathcal T_1)\times {\rm Stab}(\mathcal T_2)$ when $\rm rank (K_0(\mathcal T_i))<+\infty$).

   These  calculations suggest   a new notion of a norm, which strictly increases on $\{D^b(K(l)) \}_{l\geq 2}$. To a  triangulated category $\mathcal T$   which has property of a  phase gap we attach a non-negative number  $\left \Vert \mathcal T\right \Vert_{\varepsilon}$.   Natural assumptions on  a SOD $\mathcal T =\langle \mathcal T_1,\mathcal T_2\rangle $ imply  $\left \Vert \langle \mathcal T_1,\mathcal T_2\rangle \right \Vert_{\varepsilon}\geq {\rm max}\{\left \Vert \mathcal T_1 \right \Vert_{\varepsilon},\left \Vert\mathcal T_2\right \Vert_{\varepsilon} \}$. 
		 
		 Using the norm we define a non-trivial topology on the set of equivalence classes of  triangulated categories with a phase gap, in which  the set of discrete derived categories is a discrete subset, whereas the rationality of a smooth surface $S$ ensures  that   $[D^b(point)] \in \cl([D^b(S)] )$. 
		 
		 Viewing $D^b(K(l))$ as a non-commutative curve, we observe that it is reasonable to count non-commutative curves in any  category  which lies in a small neighborhood (w.r. to our topology) of a given non-commutative curve.  Examples show that this idea (non-commutative curve-counting) opens directions to new categorical structures and  connections to number theory and classical geometry. We  give a  general definition, which specializes  to the  non-commutative curve-counting invariants.    
	In an example arising  on the  A side we specialize our  definition to non-commutative Calabi-Yau curve-counting, where the entities we count are a Calabi-Yau modification of  $D^b(K(l))$.
		  
In the last section we speculate  that one might  consider  a holomorphic family of categories, introduced by Kontsevich, as a non-commutative extension with the norm, introduced here, playing a role similar to the classical notion of degree of an extension in Galois theory. 
	 \end{abstract}

\maketitle
\setcounter{tocdepth}{1}
\tableofcontents

\section{Introduction}

Motivated by  M. Douglas's work   in string theory, and especially by the notion of $\Pi$-stability,  T. Bridgeland defined in  \cite{Bridg1} a map:
\begin{gather} \label{Stab map} \bd \left \{ \begin{array}{c} \mbox{triangulated} \\ \mbox{categories} \end{array} \right \} &  \rTo^{\rm Stab} & \left \{ \begin{array}{c} \mbox{complex} \\ \mbox{manifolds} \end{array} \right \}.\ed  \end{gather}
For a triangulated category $\mc T$ the associated  complex manifold $\st(\mc T)$ is referred to as the space of stability conditions (or the stability space or  the  moduli space of stability conditions ) on $\mc T$. 

   Bridgeland's manifolds are expected to provide a rigorous understanding of certain moduli spaces arising in string theory. 
	Homological mirror symmetry predicts  a parallel between dynamical systems and categories, which  is being  established in  \cite{DHKK}, \cite{BS}, \cite{GNM}, \cite{KS}, \cite{Bridg1}, \cite{Bridg2}, \cite{HKK}, \cite{KNPS}. According to this analogy the stability space plays the role of the Teichm\"uller space.  However while the map \eqref{Stab map} being  well defined, it is hard to extract  global information for the stability spaces. In the present paper we determine explicitly the entire stability spaces on a new list of examples.

The map \eqref{Stab map} behaves well with respect to orthogonal decompositions (see Definition \ref{SOD}). This is easy to show and  due to lack of an appropriate  reference in the literature we have given details on this in Section \ref{st cond on Orth Dec}. In particular, there is a  a bijection 
\begin{gather} \label{formula for OD} \st(\mc T_1\oplus \mc T_2\oplus \dots \oplus \mc T_n) \cong \st(\mc T_1)\times  \st(\mc T_2)\times \dots \times \st(\mc T_n). \end{gather}
which is biholomorphism, when the categories are with finite rank Grothendieck groups. Theorem \ref{main} in this paper contains  examples of semi-orthogonal decomposition, SOD, $\mc T=\langle \mc T_1, \mc T_2 \rangle$ where ${\rm rank}(K_0( \mc T))=2 $ however $\st(\mc T)$ is not biholomorphic to $\st(\mc T_1)\times \st(\mc T_2)$. 

 The   behavior of the map \eqref{Stab map} with respect to general SOD has been studied in   \cite{CP}. This  study is  difficult and   so far  a   formula relating $\st(\langle \mc T_1,\mc T_2 \rangle)$ and $\st( \mc T_1)$, $\st(\mc T_2 )$ has not been obtained.   

 In this paper using Bridgeland stability conditions   we define (Definition \ref{main def} ) for any $0<\varepsilon<1$ a function (the domain is explained below and it  does not depend on $\varepsilon$): 
\begin{gather} \label{norm map} \bd \left \{ \begin{array}{c} \mbox{triangulated} \\ \mbox{categories} \\  \mbox{with a phase gap}  \end{array} \right \} &  \rTo^{ \norm{\cdot }_{\varepsilon}} & [0,\pi(1-\varepsilon)] \ed  \end{gather}
and prove (Theorem \ref{prop inequality}) that if $\mc T =\scal{\mc T_1,\mc T_2} $ is a  semi-orthogonal decomposition in which $\mc T$ is  proper,\footnote{by proper we mean that $\sum_{i\in \ZZ} \hom^i(X,Y)<+\infty$ for any to objects $X,Y$ in $\mc T$.} ${\rm rank}(K_0(\mc T))<\infty$,    $\mc T_1$ and  $\mc T_2$ have phase gaps, then $\mc T$ has phase gap as well and 
\begin{gather} \label{result 0} \norm{\langle \mc T_1,\mc T_2\rangle }_{\varepsilon}\geq \max\{\norm{\mc T_1}_{\varepsilon},\norm{\mc T_2}_{\varepsilon} \}. \end{gather}
For the proof of this inequality we employ the method for gluing of stability conditions  in \cite{CP},  crucial role has also \cite[Lemma 4.5]{Bridg2} which ensures certain finiteness property of a stability condition with a phase gap. 

If  $\mc T =\mc T_1 \oplus \mc T_2 $ is an orthogonal decomposition with proper $\mc T $ and ${\rm rank}(K_0(\mc T))<\infty$, then (Corollary \ref{sum of several}):
\begin{gather} \label{result 01} \norm{\mc T_1}_\varepsilon = 0 \ \ \Rightarrow \ \  \norm{\mc T_1 \oplus \mc T_2}_\varepsilon =  \norm{ \mc T_2}_\varepsilon.  \end{gather}

The function \eqref{norm map} depends on $\varepsilon \in (0,1)$, however the three   subsets of its domain   determined by the three conditions  on the first raw in  the following table do not depend on  $\varepsilon$ (Lemma \ref{maximal norms}):
\begin{gather} \label{table with norms}  \begin{array}{|c | c | c | c |} \hline 
 \mbox{Categories with:} & \norm{\cdot}_{\varepsilon}=0 &   0<\norm{\cdot}_{\varepsilon}<\pi (1-\varepsilon) &  \norm{\cdot}_{\varepsilon}=\pi (1-\varepsilon) \\
\hline 
\mbox{examples:} & \begin{array}{c}  \mbox{for any acyclic} \\ \mbox{quiver} \ Q \  D^b(Q)\\  \ \mbox{is here iff}\ Q \\  \ \mbox{is Dynkin or affine,} \\ \mbox{any discrete derived} \\ \mbox{category is here}  \end{array} & \begin{array}{c} D^b(K(l_1))\oplus \dots \oplus D^b(K(l_N)) \\ \mbox{where}\ N\in \ZZ_{\geq 1} \ l_i\geq 3 \ \mbox{for some} \  i \end{array} & \begin{array}{c} D^b(\PP^1 \times \PP^1) \\  D^b(\PP^n)  \ \ n \geq 2 \\ D^b(\FF_m) \ m\geq 0 \\ \mbox{many wild quivers } \\ \mbox{as in Prop. \ref{prop maximal norms} (a)}  \end{array}\\
\hline
\end{array} \nonumber 
\end{gather}

 Further examples can be obtained by using \eqref{result 0} and  \eqref{result 01}. In particular by blowing up the varieties in the last  column  one obtains other elements in this column (see Corollary \ref{blow up}).

In Section \ref{topology section} using \eqref{norm map}   we introduce a non-trivial  topology on the class of triangulated categories with a phase gap up to equivalence. The  function $\pi (1-\varepsilon)-\norm{\cdot}_\varepsilon$ is upper semi-continuous for this topology. The class of discrete derived categories modulo equivalence (we discuss this class in Section \ref{discrete derived categories}) is a discrete  subset w.r. to it. We show also that for any smooth complete rational surface  $S$ holds $[D^b(point)] \in \cl([D^b(S)] )$ (Corollary \ref{point in closure of rational}). These considerations motivate some questions, e.g.:

\begin{quest} Is there a smooth complete  surface  $S$, s. t.  $D^b(S)$ has a full exceptional collection  and   $[D^b(point)] \not \in \cl([D^b(S)] )$   ?  
\end{quest} A weaker question is (positive answer of the previous implies a positive answer of this)
\begin{quest} \label{quistion about homs in surf} Is there a smooth complete  surface  $S$, s. t.  $D^b(S)$  has a full exceptional collection  and  the conditions of Corollary \ref{from exc collections to maximal norm} do not hold, i. e. there exists $N$, s. t.  for any full exceptional collection $(E_0,\dots, E_n)$ in $D^b(S)$  and any $0\leq i<j\leq N$ we have $\hom^{min}(E_i,E_j)\leq N$ (see Section \ref{notations} for the notations) ?
\end{quest}
In view of Corollary \ref{point in closure of rational} and  Proposition \ref{prop maximal norms} any surface which gives a positive answer of some of these two questions would not be rational and would give a counter example to  a folklore Orlov's conjecture stating that  a surface over an algebraically
closed field admits a full exceptional collection only if it is rational.  

In all examples of $D^b(X)$, where $X$ is a smooth projective variety, for which we compute the norm, only for $X=\PP^1$ fails the condition of  Corollary \ref{from exc collections to maximal norm}. This condition fails for $D^b(K(n))$, $n\geq 1$,  and conjecture \ref{special sit} would impy that it  fails for the quivers depicted there. We would like to generalise the question \ref{quistion about homs in surf} as follows:
\begin{quest} \label{quistion about homs} Is there a smooth complete variety  $X$, different from $\PP^1$, and a natural number $N\in \NN$ s. t.   for any  exceptional collection $(E_0,\dots, E_n)$ in $D^b(X)$  and any $0\leq i<j\leq n$ we have $\hom^{min}(E_i,E_j)\leq N$ ?
\end{quest}

	Following Kontsevich-Rosenberg  \cite{KR} we denote sometimes  $D^b(K(l+1))$ by  $N\PP^l$ for $l\geq 0$. By rescaling  $\norm{\cdot}_{\frac{1}{2}}$ (see \eqref{dim map})  we define a function: \begin{gather} \label{dim map intro} \bd \left \{ \begin{array}{c} \mbox{triangulated} \\ \mbox{categories} \\  \mbox{with a phase gap}  \end{array} \right \} &  \rTo^{ \dim_{nc}} & [0,+\infty] 
	, \ed \quad  \mbox{s.t.} \quad  \begin{array}{c}
	 \dim_{nc}(\langle \mc A, \mc B \rangle)\geq \max \{\dim_{nc}(\mc A),\dim_{nc}(\mc B) \} \\ \dim_{nc}(N\PP^l) = l \qquad l \geq 0
	\end{array}   \end{gather}
	Thus the invariant $\dim_{nc}$  takes all natural numbers as values. The last column  in table \eqref{table_intro} contains categories where the invariant $\dim_{nc}$ is $\infty$, and in the other columns $\dim_{nc}$ is finite.  For an acyclic  quiver $Q$ we have $\dim_{nc}\left (D^b(Q) \right ) = 0$ iff $Q$ is Dynkin and  $\dim_{nc}\left (D^b(Q) \right ) = 1$ iff $Q$ is affine ( remark \ref{nc dim of Dynkin}.  We do not know the answer of the following:
	\begin{quest}
		Is there a category $\mc T$ with a phase gap s.t. $\dim_{nc}(\mc T) \not \in \QQ$ ?
	\end{quest}

In our topology  on the domain of \eqref{dim map intro},  up to equivalence,  whenever $\mc T$ is in a small neighborhood of $N\PP^l$ (more precisely whenever $[\mc T] \in B_\delta( N\PP^l)$ for some real $\delta >0$ as defined in \eqref{dimensions inequality})  we have    a SOD of the form $\mc T=\langle N\PP^l, \mc A \rangle $ where $\mc A$ has a phase gap, and hence $\dim_{nc}(\mc T)\geq l$. In particular if $\mc T \in B_{\delta_l}( N\PP^l)$ for arbitrary big $l$, then $\dim_{nc}(\mc T)=+\infty$, and this idea is used 
to obtain the last column of the table above. 
   \vspace{3mm}

 \textit{All examples of categories  considered in this paper  (except $D^b(pt)$ and ${\rm Fuk}(S)$)  satisfy an  incidence  $\mc T \in B_\delta( N\PP^l)$, i.e. they have a SOD $\langle N\PP^l, \mc A \rangle $, for some $l$ and some $\mc A$. Recalling that Gromow-Witten invariants count pseudo-holomorphic curves, we view such embeddings of $N\PP^l$ into $\mc T$ as analogous to a  ``pseudo-holomorphic curve'' in the category $\mc T$, and we ask  a question: can we count such entities in a given $\mc T$, how many are they  ? } 
 \vspace{3mm}

 In Section \ref{NCCC}  we show that the answer is positive.  The idea is: for a linear over a field $k$ category $\mc T$, a subgroup $\Gamma \subset {\rm Aut}(\mc T)$, and a choice of some additional restrictions $P$ to define and study  the set of subcategories of $\mc T$, which are equivalent to another chosen category $\mc A$, which satisfy  $P$, and modulo $\Gamma$.  We denote this set by $C_{\mc A, P}^{\Gamma}(\mc T)$ and define it in Definition \ref{C_l}. 
  We prefer to choose some $\mc A$, which is non-trivial but  well studied. 
  
  The studies in this paper naturally  impose  $N\PP^l$ as our first choice. In particular   we refer to $ N\PP^l$ as a non-commutative curve and  $\dim_{nc}( N\PP^l)=l$ as its ``non-commutative  genus'' (see Remark \ref{motivation for genus} for further motivation).  
   We denote $C_{N\PP^l, P}^\Gamma(\mc T)$ by  $C_{l,P}^{\Gamma}(\mc T)$ and that's the set of non-commutative curves of genus $l\geq 0$ in $\mc T$ satisfying $P$ and modulo the subgroup $\Gamma$. Furthermore, by fixing a   stability condition $\sigma \in \st(\mc T)$ we define the set of $\sigma$-semistable non-commutative curves  of genus $l\geq 1$ in $\mc T$ and denote it by $C_{l,P,\sigma}^{\Gamma}(\mc T)$ (Definition  \ref{Clsigma}). 
Non-commutative curve in $\mc T$ is just an equivalence class of exact fully faithful functors $F:N\PP^l \rightarrow \mc T$ with two functors  $F, F'$  being equivalent if the one is obtained from the other via re-parametrization, i. e.  $F'\cong F\circ \alpha$ for some $\alpha \in {\rm Aut}(N\PP^l)$.
 For the examples, which  we consider,  additional restrictions are not necessary and  $P=\emptyset$.   These examples are two affine quivers (Proposition \ref{the numbers for two quivers}) and $D^b(\PP^2)$ (Proposition \ref{Cl(PP2)}), where we  point out   non-empty and finite sets $C_{l}^{\Gamma}(\mc T)$ and their cardinalities. 
  More precisely, Proposition \ref{Cl(PP2)} concludes that $C_l^{\Aut(D^b(\PP^2))}\left (D^b(\PP^2))\right )$ is finite for all $l$ and non-empty iff $l=3 m -1$ for some Markov number $m$. Furthermore  Corollary \ref{Markov} is that the famous Markov's  conjecture in number theory and a conjecture by Tyurin (\cite[p. 100]{Rudakov1} or \cite[Section 7.2.3 ]{GorKul})  are true iff  for all  Markov numbers $m\neq1, m\neq 2$ we have $\#\left (C_{3 m -1}^{\Aut(D^b(\PP^2))}(D^b(\PP^2))\right )=2$. 
  Via the latter Corollary  in future works (starting with \cite{DK4} where will be written  the  proofs of Propositions \ref{the numbers for two quivers}, \ref{Cl(PP2)},  Corollary \ref{Markov}) we plan to approach  Markov's conjecture using   homological mirror symmetry  and applying  A side techniques  for computing the non-commutative curve-counting invariants introduced here.

 We explain  also a non-trivial  example, where $C_{1,\sigma}^{\{\rm Id\}}(\mc T)$ takes all possible values in $\{0,1,2=C_{1}^{\{\rm Id\}}(\mc T) \}$ as $\sigma$ varies in $\st(\mc T)$ (subsection \ref{dependence on stability condition}). Due to lack of space we will present the full proofs of these examples in a future work \cite{DK4} devoted to  non-commutative curve-counting. 
Here  we   prove (see  part 1.2 of the introduction below)  that $C_l^{\{\rm Id\}}(N\PP^k)=\delta_{l,k}$ for $l,k\geq 0$ and  for $l\geq 1$ we  describe the zones in $\st(N\PP^l)$, where   $C_{l,\sigma}^{\{\rm Id\}}(N\PP^l)$ is zero and one respectively. 

 Section \ref{A-side curve counting} contains an example of finite sets $C_{\mc A, P}^{\Gamma}(\mc T)$ of different origin (the proof is postponed for future work as well). Here again we don't need additional restrictions, i. e. $P=\emptyset$, and $\mc T$ is the so called Fukaya category of an elliptic curve, ${\rm Fuk}(E)$. In this case  the role of $\mc A$ is  played by a  category, denoted by $CY(l)$, which is very well studied by simplectic geometers, beginning with P. Seidel, and seems to be the right substitute for $N\PP^l$.  The question about the cardinality of  $C_{CY(l)}^{\Gamma}\left ({\rm Fuk}(S)\right )$  for higher genus curves   should be realted to counting geodesics on $S$.

Finally (Section \ref{A-side interpretation}), relating our norm to the  notion of holomorphic family of categories introduced by Kontsevich  we suggest a framework in which sequences of holomorphic families of categories are viewed as sequences of extensions of non-commutative manifolds.

 1.1.  We explain know the new  examples where $\st(\mc T_1)\times \st(\mc T_2)$ is not biholomorphic to  $\st(\langle \mc T_1, \mc T_2 \rangle)$ as well as the examples of categories $\mc T$ where we  compute (or estimate) $\norm{\mc T}_{\varepsilon}$.   
 
Let us first give some prehistory.    By definition  each stability condition  $\sigma \in \st(\mc T)$ determines a set of non-zero objects in $\mc T$ (called  \textit{ semi-stable objects})  labeled by real numbers (called \textit{phases of the semistable objects}).    The  semi-stable objects   correspond to  the so called    ``BPS'' branes in string theory.   The set of semi-stable objects  will be denoted by $\sigma^{ss}$, and $\phi_\sigma(X)\in \RR$ denotes the phase of a semi-stable $X$.  For any  $\sigma \in \st(\mc T)$ we denote by $P_\sigma^{\mc T}$    the subset of the unit circle $\{\exp(\ri \pi \phi_\sigma(X)):X\in \sigma^{ss}\} \subset \SS^1$. A categorical analogue of the density of the set of slopes of
closed geodesics on a Riemann surface   was proposed in  \cite{DHKK}. In     \cite[section 3]{DHKK}  the focus falls       on  constructing stability conditions  for which the set $P_\sigma$
is dense  in a non-trivial   arc of the  circle.   The result is  the following  characterization of the map \eqref{Stab map}, when restricted to categories of the form $D^b(Rep_k(Q))$   (from now on  $Q$ denotes an acyclic quiver, $\mc T$ denotes    a triangulated category linear over an algebraically closed field $k$):  
 \begin{gather} \label{table_intro}  \begin{array}{| c | c | c | c | c | c | c |}
  \hline
         \mbox{Dynkin quivers (e.g. \ }  $\bd[1em] \circ &\rTo & \circ\ed$ )            &   P_\sigma  \  \mbox{is always finite}                 \\ \hline
   \mbox{  Affine   quivers  (e.g. \ }  $\bd[1em] \circ & \pile{\rTo \\ \rTo }& \circ  \ed$ ) &  P_\sigma  \  \mbox{is either finite  or} \ \mbox{has exactly two limit points}         \\ \hline
    \mbox{Wild quivers (e.g. \ }  $\bd[1em] \circ & \pile{\rTo \\ \rTo\\ \rTo }& \circ  \ed$ ) &   \ P_\sigma \ \mbox{is dense in an arc }   \mbox{for  a family of stability conditions}\\ \hline
     \end{array}
  \end{gather} 
 \normalsize
In \cite[Proposition 3.29]{Dimitrov}  are  constructed stability conditions $\sigma \in \st(D^b(Q))$  with two limit
points of $P_\sigma$  for any affine quiver Q (by $D^b(Q)$ we mean $D^b(Rep_k(Q))$).

    In   \cite{Woolf} and   \cite{BPP2}  is proved  that the  stability spaces on Dynkin quivers are contractible, but the affine case  is beyond the scope of these papers.  For an integer $l\geq 1$ the $l$-Kronecker quiver  $K(l)$ (the quiver with two vertices and  $l$ parallel  arrows) is in the first, second, and third raw of the table for  $l=1$, $2$, $3$, respectively.       In \cite{Macri} are given arguments that $\st(D^b(K(l)))$ is simply-connected for any $l\geq 1$. We develop further  in \cite{DK1}, \cite{DK3}  the ideas of  Macr\`i from \cite{Macri}, in particular  we give a description of  the  entire  stability space on the acyclic triangular quiver and  prove  that  it is contractible.
 In \cite{Qiu}, \cite{BQS} and earlier by King is shown that  $\st(D^b(K(1)))$ as a complex manifold is $\CC^2$.  Recall that $D^b(K(2))\cong D^b(\PP^1)$ and $D^b(K(l))$ for $l\geq 3$ is referred to as   non-commutative projective space $N\PP^{l-1}$, introduced by Kontsevich and Rosenberg in \cite{KR} and studied further in \cite{Minamato}. In  \cite{Okada} was shown that $\st(D^b(\PP^1))\cong \CC^2$, and hence $\st(D^b(K(2)))\cong\CC^2$ (biholomorphisms). However the question:
\begin{gather} \label{question} \textit{ What is $\st(D^b(K(l)))$ for $l\geq 3$ ? } \end{gather}
was  open after the mentioned  papers.   

 One  result of the present  paper is:

\begin{theorem} \label{main} For each $l\geq 3$ there exists a biholomorphism $\st(D^b(Rep_k(K(l))))\cong \CC\times \mc H$, where $\mc H=\{z\in \CC: \Im(z)>0\}$.  \end{theorem} 
  Thus the  map \eqref{Stab map} has the same value   (up to isomorphism) on all the categories $\{ D^b(K(l)) \}_{l\geq 3}$. Stability conditions on wild quivers  whose set of phases  are dense in an arc were constructed in \cite{DHKK}, however for them the set of phases is still   not  dense in the entire $\SS^1$, i.e. $P_\sigma$ does misses a  non-trivial arc, in which case we say for  short that   $P_\sigma$ has a gap. 
  In particular all the categories in table \eqref{table_intro}  are examples of what we call in this paper a \textit{ triangulated category with phase gap}, this is a triangulated category $\mc T$ for which there exists a full\footnote{we recall what is a full stability condition in Section \ref{full st cond}}  $\sigma \in \st(\mc T)$ whose set of phases $P_\sigma^{\mc T}$ has a gap.  Stability conditions whose set of phases is not dense in $\SS^1$ and their relation to so called algebraic stability conditions have been studied in \cite{Woolf}. In particular the results in \cite{Woolf} imply that when ${\rm rank}K_0(\mc T)<\infty$, then $\mc T$ has a phase gap iff there exists a bounded t-structure in $\mc T$ whose heart is of finite length and has finitely many simple objects (Lemma \ref{criteria for phase gap}).   Whence the domain of the invariant \eqref{norm map}  contains also the   CY3 categories discussed in \cite{BS}.

From the very definition and table \eqref{table_intro} one easily derives that for any acyclic quiver $Q$:  \begin{gather} \label{result 1} \norm{D^b(Q)}_{\varepsilon} = 0  \iff  Q \mbox{ is Dynkin or affine}. \end{gather} 
Thus, we can compose the following table, concerning only the quivers $K(l)$, $l\geq 1$:
\begin{gather}  \label{table_intro1} 
 \begin{array}{| c | c | c |}\hline   Q & \norm{D^b( Q)}_\varepsilon   & \st(D^b( Q)) \\
  \hline
         $\bd[1em] \circ & \pile{\rTo \\ \rTo }& \circ  \ed$ \ \mbox{or} \ $\bd[1em] \circ &\rTo & \circ\ed$  &  \norm{D^b( Q)}_\varepsilon=0    & \CC\times\CC      \\ \hline \begin{array}{c} \\ \\  \end{array}\hspace{-3mm}
     $\bd[1em] \circ & \pile{\rTo \\ \rTo\\ \cdot \\ \rTo }& \circ  \ed$  &   \norm{D^b( Q)}_\varepsilon>0  & \CC\times\mc H \\ \hline
     \end{array}.
  \end{gather} 
 In the present paper we compute $\norm{D^b(K(l))}_\varepsilon$ for any $l$ and any $0<\varepsilon <1$. In particular we derive the following formulas:
  \begin{gather}  \label{result 2}  \norm{D^b(K(l_1))}_\varepsilon < \norm{D^b(K(l_2))}_\varepsilon  \iff  l_1<l_2 \ \mbox{and} \ 3\leq l_2  \\
	 \label{result 3} l\geq 2 \quad \Rightarrow \quad \norm{D^b(K(l))}_{\frac{1}{2}}=\arccos\left (\frac{2}{l}\right ).\end{gather}
Combining \eqref{result 1} and table \eqref{table_intro1}  we deduce that for $l\in \NN_{\geq 1}$
\begin{gather} \label{result 4} \norm{D^b(K(l))}_\varepsilon = 0 \iff \st(D^b(K(l))) \ \mbox{is affine (biholomorphic to $\CC^2$)}. \end{gather}
We expect  that the domains of validity of  \eqref{result 0} and \eqref{result 4} can be extended. Regarding \eqref{result 4} we propose:   \begin{conj}Let $0<\varepsilon <1$ and let   $Q$ be any  acyclic quiver. 

The stability space  $\st(D^b(Q))$   is affine (of the form $\CC^n$) iff $\norm{D^b(Q)}_{\varepsilon}=0$. \end{conj}

We give some criteria ensuring that $\norm{\mc T}_{\varepsilon}=\pi (1-\varepsilon)$, they imply that for many  of the wild quivers $Q$ we have $\norm{D^b(Q)}_{\varepsilon}=\pi (1-\varepsilon)$ (see Proposition \ref{prop maximal norms} (a)) and also  $\norm{D^b(X)}_{\varepsilon} = \pi (1-\varepsilon)$  where $X$ is $\PP^n$, $n\geq 2$, $\PP^1\times\PP^1$, $\mathbb F_a$, $a\geq 0$ or a smooth algebraic variety obtained from these by a sequence of blow ups  in finitely many points (see Proposition  \ref{prop maximal norms} (e), (f)), for $n=1$ we have $\norm{D^b(\PP^1)}_{\varepsilon}=0$.  Actually, the condition $\norm{\mc T}_{\varepsilon}< \pi (1-\varepsilon)$ imposes  restrictions on the full exceptional collections in $\mc T$ (see Corollary \ref{some limitations}).  

The criteria for $\norm{\mc T}_{\varepsilon}=\pi (1-\varepsilon)$  obtained here do not apply to category of the form $\mc T \cong D^b(K(l_1))\oplus D^b(K(l_2))\oplus \dots \oplus  D^b(K(l_N))$ and we do prove that $\norm{\mc T}_\varepsilon <\pi (1-\varepsilon)$ in this case, which is  a generalization of the already discussed wild Kronecker quivers  \eqref{result 2}.    

We expect that  the criterion  in Corollary \ref{from exc collections to maximal norm} does not apply to all wild quivers, and we do know that its corollary, Corollary \ref{prop maximal norms} (a), cannot be applied to  all of them, for example, to the following:
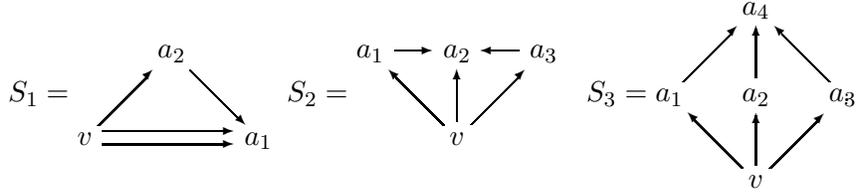
\begin{gather}  \label{special sit} S_1= \begin{diagram}[1.5em]
   &       &  a_2  &       &    \\
   & \ruTo &    & \rdTo  &       \\
v  & \pile{\rTo  \\ \rTo }  &    &       &  a_1
\end{diagram}  \ \ \  S_2= \begin{diagram}[1.5em]
a_1 &\rTo &a_2 &\lTo &a_3\\
&\luTo &\uTo &\ruTo \\
& &v
\end{diagram} \ \  \
 S_3= \begin{diagram}[1.5em]
& &a_4 & \\
&\ruTo &\uTo  &\luTo \\
a_1 &    &a_2 &    &a_3\\
&\luTo &\uTo &\ruTo \\
& &v
\end{diagram}. \nonumber     \nonumber    \ \ \ \end{gather}
 We  conjecture, that:
\begin{conj} \label{conj for S_1,S_2,S_3} For $i=1,2,3$ we have  $0<\norm{D^b(S_i)}_{\varepsilon}<\pi(1-\varepsilon) $. 
\end{conj}

\vspace{3mm}
1.2. It follows a  brief discussion (in this order)  on the proof of Theorem \ref{main},  and  of the computations of  $\norm{\mc T}_\varepsilon$. 

 For any  $\mc T$ Bridgeland defined  actions of  $\widetilde{GL}^+(2,\RR)$ (right) and of ${\rm Aut}(\mc T)$ (left) on  $\st(\mc T)$, which commute. 
The strategy for determining  $\st(D^b(\PP^1))$  in  \cite{Okada}   is to show that the quotient of $\st(D^b(\PP^1))$  for an  action
of $\CC \times \ZZ$ is isomorphic to $\CC^\star$, where  the action of $\CC$ on $\st( D^b(\PP^1))$ comes from an embedding of $\CC$ in $\wt{\rm GL}^+(2,\RR)$ and the action of $\ZZ$ comes from the subgroup of  $ \Aut(D^b(\PP^1))$ generated by the functor $(\cdot)\otimes \mc O(1)$. 
   On the one hand  in \cite{Okada} Okada  relies on the commutative geometric nature of  $D^b(K(2))$ ($\cong (D^b(\PP^1)$) and on the other hand, implicitly, he relies on  the affine  nature of the root system of $K(2)$, which are  obstacles to  answer the question \eqref{question}.  In this paper we use the ideas  of Okada in \cite{Okada} and  we   go through the mentioned   obstacles by observing how to apply simple  geometry of the action of modular subgroups on $\mc H$ and by employing the interplay between exceptional collections and stability conditions  developed in  \cite{Macri}, \cite{DK1},  
	, \cite{DK3}, \cite{Dimitrov}.
	
	In Section \ref{recollection} we  recall  facts about the action of the modular group ${\rm SL}(2,\ZZ)$ on $\mc H$, which we need. 
	
	In Section \ref{BSC}  we recall what are the actions of $\CC$ and $\Aut(\mc T)$ on $\st(\mc T)$ for any $\mc T$. It is  known that the action of $\CC$ is free and holomorphic (\cite{Okada}, \cite{BS}).  When $K_0(\mc T)$ has finite rank,  we show that the action of $\CC$ is proper on each connected component of $\st(\mc T)$, and in particular, when $\st(\mc T)$ is connected, then $\st(\mc T)\rightarrow \st(\mc T)/\CC$ is a principal holomorphic  $\CC$-bundle: Proposition \ref{principal bundle}.    
	
Section \ref{there are no Ext-nontrivial...} is devoted to the exceptional objects in $\st(\mc T_l)$, where $\mc T_l=D^b(K(l))$ for $l\geq 2$. We utilize here the method of helices \cite{BP}. Up to shifts,  there is only one helix in $\mc T_l$, which follows from \cite{WCB1}. We denote by $\{s_i\}_{i\in \ZZ}$ the helix, for which $s_1$ is the object in $Rep_k(K(l))$, which is both simple and projective.  Lemma \ref{actionLEMMA}  is the place, where  we invoke the action of ${\rm SL}(2,\ZZ)$,	here we  give a formula relating the fractions   $\{\frac{Z(s_{i+1})}{Z(s_{i})}\}_{i\in \ZZ}$  of a central charge\footnote{i.e. a group homomorphism $\bd Z:K_0(\mc T)&\rTo[1em]&\CC\ed$ such  that $Z(X)\neq 0$ for any $X\in Ob(\mc T)$.} $\bd[1em] Z:K_0(\mc T)&\rTo &\CC\ed$, it is a simple but  important  observation for  the present paper. For any two exceptional objects $E_1,E_2$ in $\mc T_l$  Lemma \ref{nonvanishings}   describes exactly those $p\in \ZZ$ for which $\hom^p(E_1,E_2)$ does not vanish. In \cite[Lemma 4.1]{Macri} is given a statement, but no proof.  The statement  of Lemma \ref{nonvanishings}  is a slight modification of  \cite[Lemma 4.1]{Macri}  and  we give a proof here.  By the equivalence $D^b(\PP^1)\cong D^b(K(2))$ the sequence $\{\mc O(i)\}_{i\in \ZZ}$ corresponds, up to translation, to the helix $\{ s_i\}_{i\in \ZZ}$.  In the non-commutative case $l\geq 3$ we still have the helix $\{s_i\}_{i\in \ZZ}$ and it plays the role of $\{\mc O(i)\}_{i\in \ZZ}$.  Corollary  \ref{the auto-eq F}   is a corollary of \cite[Theorem 0.1]{Miyachi}) and it ensures existence of a functor  $A_l\in \Aut(\mc T_l)$ for any $l\geq 2$, which is analogous to the functor $ (\cdot)\otimes \mc O(1)$ for $l=2$, more precisely it satisfies $A_l(s_i)=s_{i+1}$ for all $i\in \ZZ$. In particular from the subgroup  $\langle A_l \rangle  \subset \Aut(\mc T_l)$ we  get an action of $\CC\times \langle A_l \rangle \cong \CC\times \ZZ$ on $\st(\mc T_l)$ for all $l\geq 2$, which for $l=2$ coincides with the action on $\st(D^b(\PP^1))$ used by Okada for studying $\st(D^b(\PP^1))$. 

  In sections \ref{the union} and \ref{final section},   with the help of  ideas and results of \cite{Macri}, \cite{DK1},  
	we go on adapting   arguments of Okada about the $\CC\times \ZZ$-action on $\st(D^b(\PP^1))$   to  the non-commutative case, i.e. to   the  $\CC\times \langle A_l \rangle $-action on $\st(\mc T_l)$ for $l\geq 3$. We will explain briefly  how we utilize the   ${\rm SL}(2,\ZZ)$-action on $\mc H$.  \cite[p. 497,498]{Okada}  contain arguments about  choosing a representative in the $\CC\times \ZZ$-orbit of  any  stability condition $\sigma$ for which  $\mc O$, $\mc O(-1)$ are semi-stable and  $0<\phi_\sigma(\mc O)<\phi_\sigma(\mc O(-1)[1])\leq 1$.  These arguments rely on the fact that  for any central charge  $Z$ the vectors $\{Z(\mc O(i))\}_{i\in \ZZ}$ lie on a line in $\CC \cong \RR^2$ (see  \cite[figures 3,4,5]{Okada}), more precisely $Z(\mc O(i+1))-Z(\mc O(i))$ is the same vector for all $i\in \ZZ$, or in other words    $Z(s_{i+1})-Z(s_i)$ does not change as $i$ varies in $ \ZZ$, when $l=2$. For $l\geq 3$ this property of $\{s_i\}_{i\in \ZZ}$ fails and the  arguments and pictures on  \cite[p. 497,498]{Okada}  cannot be applied anymore. We  avoid this obstacle   by        translating this  problem to the problem of  finding a fundamental domain in $\mc H$ of a subgroup of the form $\langle \alpha_l \rangle \subset {\rm SL}(2,\ZZ)$.\footnote{In our paper  fundamental domain  is  as defined on \cite[p. 20]{Miyake}, in particular it is a  closed subset of $\mc H$.}   This translation   is encoded in  formula  \eqref{action1} in   Lemma \ref{biholomorphisms}, whose derivation  spreads throughout Subsections \ref{recollection}, ..., \ref{final section}.   The matrix $\mc \alpha_l$ appears first   in  Lemma  \ref{actionLEMMA}.  For $l=2$ this matrix is a parabolic element and for $l\geq 3$ it is a hyperbolic element in ${\rm SL}(2,\ZZ)$,\footnote{Which  are the  hyperbolic and the parabolic elements in ${\rm SL}(2,\ZZ)$ is recalled  in Subsection \ref{hpe}} which determines the difference of the type of the fundamental domains of $\langle \alpha_l \rangle$ (see Figure \eqref{1FD}) in $\mc H$, which in turn determines  the difference between the pictures on Figure  \eqref{FD}. The colored parts   in Figures \eqref{2FDhyperbolic} and \eqref{2FDparabolic} with only one of the two boundary curves  included    are in 1-1 correspondence with the   set  $\st(\mc T_l)/\CC\times \langle A_l\rangle$ for $l\geq 3$ and $l=2$, respectively. Further properties of the sets given in Figure \eqref{FD} are derived  in Corollary \ref{fundamental domain}. In the rest of Section \ref{final section} is shown how these properties and the presence of the  non-trivial real segment $(-\Delta_l,\Delta_l)$ seen on Figure \eqref{2FDhyperbolic}     imply Theorem \ref{main}. 
	
 In Sections \ref{the set of phases}, \ref{norms on kroneckers} we  compute  $\norm{{\mc T}_l}_\varepsilon$.  	By definition $\norm{{\mc T}_l}_\varepsilon$ is the supremum of\footnote{For a Lebesgue measurable subset $X\subset \mathbb S^1$ we denote by $\vol(X)$ its Lebesgue measure  with $\vol(\SS^1)=2 \pi$.} $\vol\left ( \ol{P_\sigma^l} \right )/2$ as $\sigma$ varies in the  subset $\st_\varepsilon(\mc T_l) \subset \st(\mc T_l)$ of those stability conditions $\sigma$ for which  $P_\sigma^l$ misses at least one  closed $\varepsilon$-arc (see Definitions \ref{definitin of st_varepsilon}).
	
	Sections \ref{there are no Ext-nontrivial...} and \ref{the union} are a prerequisite for the  explicit determining of the set of phases  $P_\sigma^l$ for each $\sigma \in \st(\mc T_l)$ and each $l\geq 2$, which is done  in Section \ref{the set of phases} (Proposition \ref{lemma for P_sigma}). It turns out that for $l\geq 3$   a stability condition    has $\vol\left (\ol{P_\sigma^l}\right )\neq 0$ and   satisfies  $\sigma \in \st_{\varepsilon}(\mc T_l)$ iff there  exists $j\in \ZZ$ such that $s_j, s_{j+1} \in \sigma^{ss}$ and $\varepsilon < \phi_\sigma( s_{j+1})- \phi_\sigma( s_{j})<1$,  the set $P_\sigma^l$ for such a $\sigma$  is    the set of fractions $\{n/m: (n,m) \in \Delta_+(K(l))\}$ appropriately embedded in the circle via a function  depending on the stability condition. 
	In  Lemma \ref{lemma for nmi} we  shed    light on  the structure of  the set  $\{n/m: (n,m) \in \Delta_+(K(l))\}$ (see  formulas \eqref{nmi}, \eqref{roots for K(l)}) and use it in the proof of Proposition \ref{lemma for P_sigma}. 
	
	We start  Section \ref{norms on kroneckers} by deriving  a formula expressing the non-vanishing numbers   $\vol\left ( \ol{P_\sigma^l} \right )/2$ as a smooth function  depending on  $\frac{\abs{Z(s_{j+1})}}{
	\abs{Z(s_j)}}$ and $\phi_{\sigma}(s_{j+1})-\phi_\sigma(s_j)$ for any $j\in \ZZ$ (see Proposition \ref{formula for the volume}), which is a straightforward application of the results in Section \ref{the set of phases}.   After computing partial  derivatives  of this function  we find that the supremum of   $\vol\left ( \ol{P_\sigma^l} \right )/2$ as $\sigma$ varies in $\st_\varepsilon(\mc T_l)$   is equal to $\vol\left ( \ol{P_\sigma^l} \right )/2$ where  $\sigma$ has   $s_j, s_{j+1} \in \sigma^{ss}$, $\frac{\abs{Z(s_{j+1})}}{
	\abs{Z(s_j)}}=1$ and $\phi_{\sigma}(s_{j+1})-\phi_\sigma(s_j)=\varepsilon$.   The precise formula for  $\norm{\mc T_l}_{\varepsilon}$ is in Proposition \ref{norm of Kroneckers}  and it produces \eqref{result 2}, \eqref{result 3}. In particular it follows that 
	 \begin{gather}\label{limit of norms} \lim_{l\rightarrow +\infty} \norm{\mc T_l}_{\varepsilon}=\pi (1-\varepsilon).\end{gather}
 Section  \ref{criteria} contains  examples of $\mc T$ with $\norm{\mc T}_{\varepsilon}=\pi (1-\varepsilon)$ (Proposition \ref{prop maximal norms}).  This  section is based on \eqref{limit of norms} and the observation (Proposition \ref{inequality for any exc pair}) that for any exceptional pair  $(E_1,E_2)$ in a proper $\mc T$   holds $ \norm{\scal{E_1,E_2}}_{\varepsilon}\geq \norm{\mc T_l}_{\varepsilon}$ where  $l=\hom^{min}(E_1,E_2)$. From the arguments leading to these examples it follows that  the condition $\norm{\mc T}_{\varepsilon}< \pi (1-\varepsilon)$ imposes  restrictions on $\hom^{min}(E_i,E_j)$ in a full exceptional collection $(E_0,\dots,E_n)$  (see Corollary \ref{some limitations}).

Section \ref{section sum of Kroneckers}   is devoted to the proof  (for any $N\in \ZZ_{\geq 1}$ and any $0<\varepsilon<1$)  of  
\begin{gather} \label{direct sum in intro} \norm{D^b(K(l_1))\oplus \dots \oplus  D^b(K(l_N))}_\varepsilon <\pi (1-\varepsilon).\end{gather}
Using  the results for the sets $P_\sigma^{l}$ from  Sections \ref{the set of phases}, \ref{norms on kroneckers} we show here that,  whenever $P_\sigma^l$ is contained in $C \cup - C$  for an open arc  $C\subset\SS^1$ with length less than $\pi$, then for some closed arc $p_\sigma^l\subset C \cap \ol{P_\sigma^l}$ the set  $\ol{P_\sigma^l}\setminus (p_\sigma^l\cup-p_\sigma^l)$ is at most  countable, and furthermore, provided that the  length of $C$ is fixed,  we show that when some of the end points of  $p_\sigma^l$ is very close to some of the end points of $C$, then $p_\sigma^l$ itself has very small length (Corollary \ref{the closed arcs}).   
 Due to the  fact, proven in   Section \ref{st cond on Orth Dec}, that for any orthogonal decomposition $\mc T=\mc T_1\oplus\dots\oplus \mc T_n $ and any $\sigma \in \st(\mc T)$ holds $P_\sigma^{\mc T}=\bigcup_{i=1}^n P_{\sigma_i}^{\mc T_i}$, where  $(\sigma_1, \dots , \sigma_n)$ is the value of the map \eqref{formula for OD} at $\sigma$ (see Proposition \ref{lemma for orthogonal composition 1} and Corollary \ref{gaps in products}), the proof of \eqref{direct sum in intro} reduces to proving that the measure of union of arcs $\cup_{i=1}^n p_\sigma^{l_i}\subset C$  of the type explained above, cannot become arbitrary close to the length of  $C$. Having proved this for one arc (in Section \ref{norms on kroneckers}) we perform induction and the tool  for the induction step is the  already discussed Corollary \ref{the closed arcs}.
 
 In  Section \ref{discrete derived categories} we discuss the class of discrete derived categories and show that  $\norm{\mc T}_\varepsilon=0$ for any such category. These categories  were  introduced  by Vossieck \cite{Vossieck}, they were classified in \cite{BGS} and thoroughly studied in \cite{BPP1}, whereas the topology of the stability spaces on them were studied in \cite{BPP2}, \cite{Woolf}, in particular it was shown that these spaces are all contractible.   This class contains the categories $\{D^b(Q): Q \ \mbox{is Dynkin}\}$, and the discrete derived categories not contained in this list are of the form $D^b(\Lambda(r,n,m))$ for $n\geq r \geq 1$ and $m\geq 0$, where $\Lambda(r,n,m)$ is the path algebra of the quiver with relations shown on \cite[Section 4.3, Figure 1]{Woolf}. Actually we show that  if $\mc T$ is a  category with phase gap, s.t. every heart of a bounded t-structure has finitely many indecomposable objects up to isomorphism (the discrete derived categories have this property) , then $\norm{\mc T}_\varepsilon = 0$.  
 
 1.3. Having explained what is the helix $\{s_j\}_{j\in \ZZ}$ in $D^b(K(l+1))\cong N\PP^l$ for $l\geq 1$ we can explain what we mean by a \textit{$\sigma$-semistable non-commutative curve} (see Definition \ref{Clsigma} for precise statement). Let $\st(\mc T)\neq \emptyset$ and $\sigma \in \st(\mc T)$. Recall that a non-commutative curve of genus $l$ in $\mc T$ is equivalence class of fully faithful exact functors from $N\PP^l$ to $\mc T$ (equivalence is re-parametrization in the domain), we will say that the curve is $\sigma$-semistable if for infinitely many $j\in \ZZ$ the object $F(s_j)\in \mc T$ is $\sigma$-semistable object (it does not matter which functor $F$ we take as a representative). We denote the set of $\sigma$-semistable non-commutative curves of genus $l$, satisfying the additional restrictions $P$ and modulo subgroup $\Gamma\subset {\rm Aut}(N\PP^l)$, by $C_{l,P,\sigma}^{\Gamma}(\mc T)$. The basic example is  $C_{l,\sigma}^{\{\rm Id\}}(N\PP^l)$, $l\geq 1$. First note that from Remark \ref{functor from NP tp NP} it follows that  $C_j^{\{\rm Id\}}(N\PP^j)=\delta_{i,j}$ for $i,j\geq 0$. In Lemmas \ref{sections}, \ref{sections1}  is shown that for any $\sigma \in \st( N\PP^l)$ we have one of the following possibilities

\begin{itemize}
	\item only two elements in the helix, of the form   $s_j, s_{j+1}$, are semi-stable and $\phi_\sigma(s_{j+1})>\phi_\sigma(s_j)+1$, in particular  $C_{l,\sigma}^{\{\rm Id\}}(N\PP^l)=0$ 

	\item all elements $\{s_j\}_{j\in \ZZ}$
 are semistable  and $\phi_\sigma(s_{j+1})=\phi_\sigma(s_j)+1$ for some $j\in \ZZ$, hence $C_{l,\sigma}^{\{\rm Id\}}(N\PP^l)=1$

\item all elements $\{s_j\}_{j\in \ZZ}$
 are semistable  and $\phi_\sigma(s_j)<\phi_\sigma(s_{j+1})<\phi_\sigma(s_j)+1$ for all $j\in \ZZ$, hence $C_{l,\sigma}^{\{\rm Id\}}(N\PP^l)=1$
\end{itemize}
For details on this example see Proposition \ref{kroneckers wc}. 
Another example where we discuss the numbers $C_{l,\sigma}^{\{\rm Id\}}(\mc T)$ can be seen in the end of subsection \ref{dependence on stability condition}. Here   $C_{l,\sigma}^{\{\rm Id\}}(\mc T)$ takes all possible values in $\{0,1,2=C_{1}^{\{\rm Id\}}(\mc T) \}$ as $\sigma$ varies in $\st(\mc T)$. 
 
  In the last section  we pose some questions and conjectures  about possible interplay between our norm, the  notion of holomorphic family of categories,  introduced by Kontsevich, and  the classical notion of uni-rationality.

\textit{{\bf Acknowledgements:}}
The authors wish to express their deep gratitude to  Maxim Kontsevich   for his  interest.
 The authors  are very thankful to Denis Auroux for his help on the last Section \ref{A-side}. 

  Parts of this  work  was carried out    during the stay 01.07.15-30.06.16 of the first author  at the Max-Planck-Institute f\"ur Mathematik Bonn and  his stay 01.07.16-30.06.17   at International Centre for Theoretical Physics Trieste and  G.  Dimitrov gratefully acknowledges the support and the excellent conditions  at these institutes. 
The first author  was partially supported by ERC Gemis grant and by a FWF Project P 27784.
 The second author   was  supported by Simons research grant, NSF DMS 150908, ERC Gemis, 
DMS-1265230, 
 DMS-1201475 
OISE-1242272 PASI. 
Simons collaborative Grant - HMS.  HSE-grant, HMS and automorphic  forms. The second  author is partially supported by Laboratory of Mirror Symmetry NRU HSE, RF Government grant, ag. № 14.641.31.0001

 \section{Notations} \label{notations}  In this paper the letters ${\mathcal T}$ and $\mc A$ denote   a triangulated category and an abelian category, respectively, linear over an  algebraically closed  field  $k$. The shift functor  in ${\mathcal T}$ is designated by $[1]$.   We write $\Hom^i(X,Y)$ for  $\Hom(X,Y[i])$ and  $\hom^i(X,Y)$ for  $\dim_k(\Hom(X,Y[i]))$, where $X,Y\in \mc T$.  For $X,Y\in\mc  A$,  writing $\Hom^i(X,Y)$, we consider $X,Y$ as elements in   $\mc T=D^b(\mc A)$, i.e.  $\Hom^i(X,Y)={\rm Ext}^i(X,Y)$.

 A triangulated category $\mc T$ is called  \textit{ proper} if $\sum_{i\in \ZZ} \hom^i(X,Y)<+\infty$ for any two objects $X,Y$ in $\mc T$.
For  $X,Y\in \mc T$  in a proper $\mc T$, we denote:
\begin{gather} \label{hom min}\hom^{min}(X,Y) = \left  \{ \begin{array}{c c} \hom^i(X,Y) & \mbox{if} \ i= \min \{j: \hom^j(X,Y)\neq 0\}>-\infty \\ 0 & \mbox{otherwise}.\end{array} \right  .  \end{gather}

  We write $\langle  S \rangle  \subset \mc T$ for  the triangulated subcategory of $\mc T$ 
   generated by $S$, when $S \subset Ob(\mc T)$.

  An \textit{exceptional object}  is an object $E\in \mc T$ satisfying $\Hom^i(E,E)=0$ for $i\neq 0$ and  $\Hom(E,E)=k $. We denote by ${\mc A}_{exc}$, resp. $D^b(\mc A)_{exc}$,  the set of all
    exceptional objects of  $\mc A$, resp. of  $D^b(\mc A)$.

An \textit{exceptional collection} is a sequence $\mc E = (E_0,E_1,\dots,E_n)\subset \mc T_{exc}$ satisfying $\hom^*(E_i,E_j)=0$ for $i>j$.    If  in addition we have $\langle \mc E \rangle = \mc T$, then $\mc E$ will be called a full exceptional collection.  For a vector $\textbf{p}=(p_0,p_1,\dots,p_n)\in \ZZ^{n+1}$ we denote $\mc E[\textbf{p}]=(E_0[p_0], E_1[p_1],\dots, E_n[p_n])$. Obviously  $\mc E[\textbf{p}]$ is also an exceptional collection. The exceptional collections of the form   $\{\mc E[\textbf{p}]: \textbf{p} \in \ZZ^{n+1} \}$ will be said to be shifts of $\mc E$. 

If an exceptional collection  $\mc E = (E_0,E_1,\dots,E_n)\subset \mc T_{exc}$ satisfies  $\hom^k(E_i,E_j)=0$ for any $i,j$ and for $k\neq 0$, then it is said to be \textit{strong exceptional collection}.   

For two exceptional collections $\mc E_1$, $\mc E_2$ of equal length we  write $\mc E_1 \sim \mc E_2$ if $\mc E_2 \cong \mc E_1[\textbf{p}]$ for some $\textbf{p} \in \ZZ^{n+1}$.

 \textit{An abelian category $\mc A$ is said to be hereditary, if ${\rm Ext}^i(X,Y)=0$ for any  $X,Y \in \mc A$ and $i\geq 2$,  it is said to be of finite length, if it is Artinian and Noterian.
}

By $Q$ we denote an acyclic quiver and   by  $D^b(Rep_k(Q))$, or just $D^b(Q)$, -  the derived category of the category of representations of $Q$.
	
For an integer $l\geq 1$ the $l$-Kronecker quiver  (the quiver with two vertices and  $l$ parallel  arrows)  will be denoted by  $K(l)$.

 The real and the imaginary parts of a complex number $z\in \CC$ will be   denoted by $\Re(z)$ and $\Im(z)$, respectively,  and by $\mc H$  we denote the upper half plane, i. e. $\mc H=\{ z\in \CC : \Im(z)>0  \}$. For any  complex number $z \in \mc H$ we denote by $\Arg(z)$  the unique $\phi \in (0,\pi)$ satisfying $z=\abs{z} \exp(\ri \phi)$.  

For a complex number $z=(a+\ri b)$, $a,b\in \RR$ we denote $\Im(z)=b$, $\Re(z)=a$.

The letter $\HH$ will denote the upper half plane  with the negative real axis included, i. e.  $\HH=\{ r \exp(\ri \pi t): r>0 \ \ \mbox{and} \ \ 0< t \leq 1 \}$.

  For an element $\alpha$ in a group $G$ we denote by $\langle \alpha \rangle \subset G$ the subgroup  $\langle \alpha \rangle = \{\alpha^i\}_{i\in\ZZ}$.  
	
	If $A\subset B$ are subsets in a top. space $X$, we denote by  ${\rm Bd}_{ B}(A)$ the boundary of $A$ w.r. to $B$, and by $L_B(A)$ the set of limit points of $A$ w.r. to  $B$.

	\section{ On Bridgeland stability conditions} \label{BSC}  We use freely the axioms and notations  on stability conditions introduced by Bridgeland in \cite{Bridg1} and some additional notations  used in \cite[Subsection 3.2]{DK1}.  In particular, the underlying set of the manifold $\st(\mc T)$ is the set of locally finite stability conditions on $\mc T$ and for $\sigma = (Z, {\mathcal P}) \in \st(\mathcal T)$ we  denote by
$\sigma^{ss}$ the  set of  $\sigma$-semistable
objects, i. e. \be
\label{sigma^{ss}} \sigma^{ss}=\cup_{t \in \RR} {\mathcal P}(t)\setminus \{0\}.
\ee 
Also for a hearth $\mc A$ of bounded t-structure  in $\mc T$ we denote by $\HH^{\mc A}\subset \st(\mc T)$ the subset of the stability conditions $(Z, \mc P)\in \st(\mc T)$ for which $\mc P(0,1]=\mc A$ (see  \cite[Definition 2.28]{Dimitrov}).

Recall  that one of Bridgeland's axioms   \cite{Bridg1}  is: for any nonzero $X \in Ob({\mathcal T})$ there exists  a
diagram of distinguished  triangles called \textit{Harder-Narasimhan filtration}:
\be
\begin{diagram}[size=1em] \label{HN filtration}
0 & \rTo      &     &       &   E_1 & \rTo      &     &        & E_2 & \rTo & \dots & \rTo & E_{n-1} & \rTo      &     &       & E_n=X \\
  & \luDashto &     & \ldTo &       &\luDashto  &     &  \ldTo &     &      &       &      &         & \luDashto &     & \ldTo &       \\
  &           & A_1 &       &       &           & A_2 &        &     &      &       &      &         &           & A_n &
\end{diagram}
\ee where $\{ A_i  \in {\mathcal P}(t_i) \}_{i=1}^n $, $t_1 > t_2
> \dots > t_n $ and $A_i$ is non-zero object for any $i=1,\dots,n$
(the non-vanishing  condition makes the factors $\{ A_i  \in
{\mathcal P}(t_i) \}_{i=1}^n$ unique up to isomorphism). Following 
\cite{Bridg1} we denote  $\phi_\sigma^-(X):=t_n$,
$\phi_\sigma^+(X):=t_1$, and the phase of a semistable object $A
\in {\mathcal P}(t)\setminus \{0\}$
 is denoted by $\phi_{\sigma}(A):=t$. The positive integer:
\begin{gather} \label{mass} m_{\sigma}(X)=\sum_{i=1}^n \abs{Z(A_i)}\end{gather} 
is called the mass of $X$ w.r. to $\sigma$(\cite[p.332]{Bridg1}).  We will use also the following axioms \cite{Bridg1}:
\begin{gather}
\label{phase formula}  X\in \sigma^{ss} \qquad
\Rightarrow \qquad  Z(X)=m_{\sigma}(X)\ \exp(\ri \pi \phi_\sigma(X)), \ m_{\sigma}(X)=\abs{Z(X)} >0 \\
\label{vanishing formula} X, Y\in \sigma^{ss} \ \ \phi_\sigma(X)>
\phi_\sigma(Y) \quad  \Rightarrow \quad  \Hom(X,Y)=0.
\end{gather}

Finally we note that: 
\begin{gather} \label{direct sums} X\cong X_1\oplus X_2 \Rightarrow \begin{array}{c}  m_{\sigma}(X)=m_\sigma(X_1)+m_{\sigma}(X_2) \\   \phi_\sigma^-(X)={\rm min}\{ \phi_\sigma^-(X_1),  \phi_\sigma^-(X_2)\} \\ \phi_\sigma^+(X)={\rm max}\{ \phi_\sigma^+(X_1),  \phi_\sigma^+(X_2)\}  \end{array}, \end{gather}
which follows easily from \eqref{mass} and the arguments  for the proof of \cite[Lemma 2.25]{Dimitrov}.

\subsection{Actions on \texorpdfstring{$\st(\mc T)$}{\space}} \label{actions}
\subsubsection{The universal covering group of $GL^+(2,\RR)$.}   The universal covering group  $\wt{GL}^+(2,\RR)$ of  $GL^+(2,\RR)$ can be constructed as follows (we  point  the steps without proving them). First step is to show that the following set  with the specified bellow operations and metric is a topological group: 
\begin{gather} \wt{GL}^+(2,\RR)= \left \{(G,\psi ) : \begin{array}{c} G \in GL^+(2,\RR),  \quad \psi \in C^{\infty}(\RR) \\  \forall t\in \RR \  \psi'(t) > 0, \ \psi(t+1)=\psi(t)+1,   \frac{G(\exp(\ri \pi t))}{\abs{G(\exp(\ri \pi t))}}=\exp(\ri \pi \psi(t) ) \end{array} \right \} \\
\mbox{unit element:} \ \ ({\rm Id}_{\CC}, {\rm Id}_{\RR}) \\
\mbox{multiplication:} \ \ \left ( (G_1, \psi_1), (G_2, \psi_2) \right ) \mapsto (G_1 \circ G_2, \psi_1 \circ  \psi_2) \\ 
\mbox{inverse element:} \ \  (G , \psi)  \mapsto (G^{-1}, \psi^{-1} )\\
\mbox{metric:} \label{metric on univ cover} \ \ d \left ( (G_1, \psi_1), (G_2, \psi_2) \right ) = \sup_{t\in \RR}\{\abs{G_1(\exp(\ri \pi t))-G_2(\exp(\ri \pi t))}, \abs{\psi_1(t)-\psi_2(t)} \}.
 \end{gather}
Second step is to show that the following  is a covering map:
\begin{gather} \bd \wt{GL}^+(2,\RR) & \rTo^{\pi } &  GL^+(2,\RR) \ed  \ \ \  \ (G , \psi)  \mapsto G. \end{gather}  
The subset  $U_\varepsilon=\{G\in GL^+(2,\RR) ; \sup_{t\in \RR}\{\abs{G(\exp(\ri \pi t))-\exp(\ri \pi t)} \}<\sin(\pi \varepsilon) \}$\footnote{ neighborhood of  ${\rm Id}_\CC\in  GL^+(2,\RR)$}  is evenly covered by a family of  open subsets  $\{(G,\psi); G \in U_\varepsilon \ \  \sup_{t\in \RR} \abs{\psi(t)-t-2 k} <\varepsilon  \}$ indexed by   $k\in \ZZ$ for small enough $\varepsilon$. 
 In particular one obtains a structure or a Lie Group on $\wt{GL}^+(2,\RR)$ such that $\pi $ is a morphism of Lie groups.
 
Finally, one can show that  $ \wt{GL}^+(2,\RR) $ is simply connected  by recalling that $\pi_1( GL^+(2,\RR))\cong \ZZ$ is generated by  $\SS^1= SO(2)\subset GL^+(2,\RR) $  and then by finding  the lifts of this path in  $\wt{GL}^+(2,\RR)$. 

\begin{remark}\label{an element in tilde} For any  $0<\varepsilon <1$,  $0<\varepsilon' <1$ there exists unique  $g_{\varepsilon,\varepsilon'}=(G,\psi)\in \wt{\rm GL}^+(2,\RR)$ such that
 $G^{-1}(1)=1$ and  $G^{-1}(\exp(\ri \pi \varepsilon ))=\exp(\ri \pi \varepsilon')$ and 
$\psi(0)=0 , \psi(1)=1,  \psi(\varepsilon')=\varepsilon  $, in particular :
\begin{gather} \label{an element in tilde remark}  \psi ([0,\varepsilon'])=[0,\varepsilon], \ \ \psi ([\varepsilon', 1])=[\varepsilon, 1]\end{gather}
 Furthermore, $(g_{\varepsilon, \varepsilon'})^{-1}=g_{\varepsilon', \varepsilon}$.
 \end{remark}

The right action of $\wt{GL}^+(2,\RR)$ on $\st(\mc T)$ is defined by (recall \cite{Bridg1} ): 
\begin{gather} \label{GL tilde} \st(\mc T) \times \wt{GL}^+(2,\RR) \rightarrow  \st(\mc T) \qquad  \left ( ( Z, \mc P), (G, \psi ) \right ) \mapsto (Z, \mc P) \cdot (G, \psi ) =  ( G^{-1} \circ Z, \mc P \circ \psi ). \end{gather}
 
Using the formula \eqref{metric on univ cover} determining the topology on $\wt{GL}^+(2,\RR)$   and the basis of the topology in $\st(\mc T)$ explained on \cite[p. 335]{Bridg1} one can show that the function   in \eqref{GL tilde} is continuous. 

 We recall also (see \cite[Theorem 1.2]{Bridg1}) that the projection $\bd \st(\mc T)&\rTo^{proj}&\Hom(K_0(\mc T),\CC)\ed$,   $proj(Z, \mc P)=Z$ restricts to a local biholomorphism between each connected component of $\st(\mc T)$ and a corresponding vector subspace of  $\Hom(K_0(\mc T),\CC)$ with a well defined linear topology (when $\rm rank(K_0(\mc T))<+\infty$ this is the ordinary linear topology). Note also  that the results in \cite{Bridg1} imply  that  $\st(\mc T)$ is locally path connected (follows from the results in \cite[Section 6]{Bridg1} and \cite[Theorem 7.1]{Bridg1}), therefore  the components and the  path components  of   $\st(\mc T)$ coincide and they are  open subsets in $\st(\mc T)$. 

Finally, assume for simplicity that  ${\rm rank}(K_0(\mc T))<+\infty$.  Due to continuity of \eqref{GL tilde} it follows that for each connected component $\Sigma$ of  $\st(\mc T)$ the action \eqref{GL tilde} restricts to a continuous action  $\Sigma \times \wt{GL}^+(2,\RR) \rightarrow  \Sigma $ and  it is easy to show that there is a commutative diagram: 
\be \label{diagram for actions} \bd[size=1.5em] \Sigma \times \wt{GL}^+(2,\RR) & \rTo               & \Sigma \\
                    \dTo^{proj_{\vert}\times \pi}               &                           &      \dTo^{proj_{\vert}}                      \\
											V(\Sigma) \times   GL^+(2,\RR)        &   \rTo & 	V(\Sigma)   \ed \ee
											where $V(\Sigma)\subset \Hom(K_0(\mc T),\CC)$ is the corresponding to $\Sigma$ vector subspace, such that the vertical arrows are local diffeomorphisms (the right arrow is local biholomorphism), and the lower horizontal arrow is  an action of the form $(A,G) \mapsto A\circ G^{-1}$ on $V(\Sigma)$. Now it follows that the upper horizontal arrow is smooth, and therefore \eqref{GL tilde} is smooth as well.   

\subsubsection{The action of \texorpdfstring{$\CC$}{\space}} \label{the action of CC} There is a Lie group homomorphism $\CC \rightarrow \wt{GL}^+(2,\RR)$ given by $\lambda \mapsto\left  ({\rm e}^{-\lambda}, {\rm Id}_\RR - \Im(\lambda)/\pi\right )$. And composing the action \eqref{GL tilde}  with this  homomorphism results in the action  \eqref{star} below. This action  is free  \cite[Definition 2.3, Proposition 4.1]{Okada}.  It is easy to show that for any $X\in \mc T$, $\sigma \in \st(\mc T)$, $z\in \CC$ hold  the properties in \eqref{star properties 1},  \eqref{star properties 2} below, and  the HN filtrations of $X$ w.r. to $\sigma$ and  to $z\star \sigma$ are the same: 
\begin{gather}\label{star}\bd \CC\times \st(\mc T) &\rTo^{\star}& \st(\mc T)\ed \quad  z\star\left (Z,\{\mc P(t)\}_{t\in\RR}\right )= \left (\re^z Z, \{ \mc P(t-\Im(z)/\pi)\}_{t\in \RR} \right ) \\ \label{star properties 1}
(z\star\sigma)^{ss}=\sigma^{ss} \qquad  \phi_{z\star \sigma}(X)= \phi_{ \sigma}(X)+\Im(z)/\pi \qquad X\in \sigma^{ss} \\ 
\label{star properties 2}  \phi^{\pm}_{z\star \sigma}(X)=\phi^{\pm}_{ \sigma}(X)+\Im(z)/\pi; \qquad m_{z\star \sigma}(X) =\re^{\Re(z)}m_{\sigma}(X).\end{gather} 
 
\begin{prop} \label{principal bundle} Let   $K_0(\mc T)$ have finite rank and let $\Sigma \subset \st(\mc T)$ be a  connected component. Then the action \eqref{star} restricted to $\Sigma$ is proper. In particular, $\Sigma/\CC$ with the quotient topology carries a structure of a complex manifold, s. t. the projection $\pr  : \Sigma \rightarrow  \Sigma/\CC$ is a holomorphic $\CC$-principal bundle of complex dimension $\dim_\CC(\st(\mc T))-1$.
\end{prop}
\bpr The action is holomorphic (now we get a diagram \eqref{diagram for actions} with $\CC$ instead of $\wt{GL}^+(2,\RR)$, $\CC^{\star}$ instead of $GL^+(2,\RR)$, and now both the vertical arrows  are local biholmorphisms, whereas  the lower horizontal  arrow is holomorphic) and free  (see e.g. \cite[Proposition 4.1]{Okada}). If we show that the function 
\begin{gather} \gamma : \CC \times \Sigma \rightarrow  \Sigma\times  \Sigma \qquad (\lambda, \sigma) \mapsto (\lambda \star \sigma, \sigma)\end{gather} 
is proper, then the proposition follows from \cite[Proposition 1.2]{Huck}. Let  $K_i\subset \Sigma$, $i=1,2$  be two compact subsets.  Since $\st(\mc T)\times \st(\mc T)$ is locally compact, it is enough to show that  $\gamma^{-1}(K_1\times K_2)\subset\CC\times\Sigma$ is compact.  

\cite[Proposition 8.1.]{Bridg1} says that assigning to any two $\sigma_1,\sigma_2 \in \st(\mc T)$ the following:
\begin{gather} \label{metric}  d(\sigma_1,\sigma_2)=\sup_{0\neq X\in\mc T}\left  \{\abs{\phi_{\sigma_1}^-(X) -\phi_{\sigma_2}^-(X)},\abs{\phi_{\sigma_1}^+(X) -\phi_{\sigma_2}^+(X)},\abs{\log \frac{m_{\sigma_2}(X)}{m_{\sigma_1}(X)}} \right \}\in [0,+\infty]\end{gather} 
defines a generalized metric whose topology coincides with the topology of $\st(\mc T)$. Furthermore, since  $\Sigma$ is a connected component, it follows that $d(\sigma_1,\sigma_2)<+\infty$ for any two  $\sigma_1,\sigma_2 \in \Sigma$, hence $d$ is a usual metric on $\Sigma$. In particular, the function $d:\Sigma \times \Sigma \rightarrow \RR $ is continuous, hence there exists $M\in \RR_{>0}$, s.t.
\begin{gather} \label{interm} \forall x,y\in K_1\times K_2  \qquad d(x,y)\leq M. \end{gather}

 From \eqref{star properties 2} and \eqref{metric}  we see that 
\begin{gather}\label{interm1}  \forall \sigma \in \st(\mc T) \ \forall \lambda \in \CC \quad   \max\left \{ \abs{\Re(\lambda)}, \abs{\Im(\lambda)/\pi}\right \}\leq d\left (\sigma, \lambda \star \sigma\right ). \end{gather}
Asume that $(\lambda,\sigma)\in \gamma^{-1}(K_1\times K_2)$, then $\gamma(\lambda,\sigma)=(\lambda \star \sigma, \sigma) \in K_1\times K_2$, therefore  $ \sigma\in  K_2$ and  by \eqref{interm}, \eqref{interm1} it follows that $ \max\left \{ \abs{\Re(\lambda)}, \abs{\Im(\lambda)/\pi}\right \}\leq d(\lambda \star \sigma, \sigma) \leq M$, thus we obtain:
\begin{gather} \gamma^{-1}(K_1\times K_2)\subset \{\lambda \in \CC; \abs{\Re(\lambda)}\leq M,  \abs{\Im(\lambda)}\leq \pi M\} \times K_2, \end{gather} 
therefore $\gamma^{-1}(K_1\times K_2)$ is a closed  subset of a compact subset, hence   $\gamma^{-1}(K_1\times K_2)$  is compact and the proposition is proved.
\epr

\begin{coro} \label{trvial bundle 1} Let $\mc T$ has a contractible  stability space, and let  $\dim_\CC(\st(\mc T))={\rm rank}(K_0(\mc T))=2$. Then  $\st(\mc T)$ is biholomorphic to one of the two: $\CC \times \CC$ or  $\CC \times \mc H$ and the quotient map   $\pr: \st(\mc T) \mapsto \st(\mc T)/\CC$ is a  trivial $\CC$-principal bundle. 
\end{coro}
\bpr Proposition \ref{principal bundle} imply that $\pr: \st(\mc T) \mapsto \st(\mc T)/\CC$ is a $\CC$-principal bundle and now we are given that  the total space is  contractible and the  base is a one dimensional connected complex manifold.  Now from the long sequence of the homotopy groups associated to a fibration  one can deduce that  $\st(\mc T)/\CC$ is contractible one dimensional complex manifold (see e.g. \cite[p. 82]{Tromba} for more details on the arguments).  By uniformisation theorem $\st(\mc T)/\CC$  is biholomorphic either to $\CC$ or to $\mc H$. Now the corollary follows from the theorem that every fiber bundle over a non-compact Riemann
surface is trivial, provided the structure group G is connected.(\cite{Grauert}, \cite{Rohrl}). \epr

\subsubsection{The action of \texorpdfstring{$\Aut(\mc T)$}{\space}} There is a left action  of the group of exact auto-equivalences $\Aut (\mc T)$ on $\st(\mc T)$, which commutes with the action \eqref{star}\cite[Lemma 8.2]{Bridg1}. This action is determined as follows:
\begin{gather} \label{left action} \Aut(\mc T) \times \st(\mc T) \ni \left  (\Phi, (Z,\{ \mc P(t)\}_{t\in \RR}) \right ) \mapsto \left (Z\circ [\Phi]^{-1}, \left
\{  \overline{\Phi(\mc P(t))} \right \}_{t\in \RR}\right )\in \st(\mc T), \end{gather}
where $[\Phi]:K_0(\mc T) \rightarrow K_0(\mc T)$ is the induced isomorphism (we will often omit specifying the square brackets) and $\overline{\Phi(\mc P(t))} $ is the full isomorphism closed subcategory containing $\Phi(\mc P(t))$. 

When ${\rm rank}(K_0(\mc T))<+\infty$, it is easy to show that ${\rm Aut}(\mc T)$ acts via   \eqref{left action}  biholomorphically  on $\st(\mc T)$. 

 For any $\Phi \in \Aut(\mc T)$, $\sigma\in \st(\mc T)$ let us denote $\sigma=(\mc P_{\sigma}, Z_\sigma)$, $\Phi\cdot \sigma=(\mc P_{\Phi \cdot \sigma}, Z_{\Phi\cdot \sigma})$, then  we have: 
\begin{gather} \label{aut properties 1}
(\Phi\cdot \sigma)^{ss}=\overline{\Phi(\sigma^{ss})} \qquad  \phi_{\Phi \cdot \sigma}(X)= \phi_{ \sigma}(\Phi^{-1} X) \qquad X\in (\Phi\cdot \sigma)^{ss} \\ 
\label{aut properties 2} Z_{\Phi \cdot \sigma}(X)=Z_\sigma(\Phi^{-1}(X)) \qquad X\in \mc T.\end{gather}

\section{Triangulated categories with phase gaps and their norms   }

\subsection{Full stability conditions} \label{full st cond}

We start this section by recalling what is meant when saying that a stability condition is full.

\textit{Full stability condition}  on $K3$ surface is defined in   \cite[Definition 4.2]{Bridg1}. Analogous definition can be given for any triangulated category $\mc T$ and locally finite stability condition  whose central charge factors  through a given group homomorphism  $ch: K_0(\mc T) \rightarrow  \ZZ^n $.  

 When $K_0(\mc T)$ has finite rank, we choose always the trivial homomorphism $K_0(\mc T)\rightarrow K_0(\mc T)$. 
 Now the projection $\bd \st(\mc T)&\rTo^{proj}&\Hom(K_0(\mc T),\CC)\ed$,   $proj(Z, \mc P)=Z$ restricts to a local biholomorphism between each connected component of $\st(\mc T)$ and a corresponding vector subspace of  $\Hom(K_0(\mc T),\CC)$ (see \cite[Theorem 1.2]{Bridg1}). A stability condition $\sigma \in  \st(\mc T)$ in this case is a  full stability condition, if the vector subspace of  $\Hom(K_0(\mc T),\CC)$ corresponding to the connected component $\Sigma$ containing $\sigma$ is the entire
 $\Hom(K_0(\mc T),\CC)$, which is equivalent to the equality $\dim_\CC(\Sigma) = {\rm rank} (K_0(\mc T))$. 

As we will see later all  stability conditions on $K(l)$ are full, for all $l\geq 1$ (see table \ref{table_intro1}).  It is reasonable to hope that, whenever $\st(\mc T) \neq \emptyset$,  there are always full stability conditions on $\mc T$ and, to the best of our knowledge, there are no counterexamples of this statement so far.   

\subsection{The $\varepsilon$-norm of a triangulated category} 

Recall that  for $\sigma = (Z, \mc P) \in \st(\mc T)$ we denote (see \cite[Section 3]{DHKK}):
\begin{equation} \label{the set of phases def} P_\sigma^{\mc T}= \{\exp(\ri \pi \phi_\sigma(X)) : X \in \sigma^{ss}\}=\{\exp(\ri \pi t) : t\in \RR \ \mbox{and} \ \mc P(t)\neq \{0\} \},\end{equation} 

Here we will use also the notation:
\begin{equation} \label{the set of real phases def} \wt{P}_\sigma^{\mc T}=\{t\in \RR :  \ \mc P(t)\neq \{0\} \} \ \ \Rightarrow \ \   P_\sigma^{\mc T}=\exp\left (\ri \pi \wt{P}_\sigma^{\mc T}\right ).\end{equation} 
The sets $ P_\sigma^{\mc T}$ and $\wt{P}_\sigma^{\mc T}$ satisfy $P_\sigma^{\mc T}=- P_\sigma^{\mc T}$,   $\wt{P}_\sigma^{\mc T}+1=\wt{P}_\sigma^{\mc T}$. In particular  the closures  $\ol{ P_\sigma^{\mc T}}$, $\ol{\wt{P}_\sigma^{\mc T}}$ satisfy:
\begin{gather} \label{translate the volume}  \vol\left (\ol{P_\sigma^{\mc T}}\right ) = 2 \pi\mu \left (\ol{\wt{P}_\sigma^{\mc T} } \cap [0,1]\right )= 2 \pi \int_{[0,1]\cap \ol{\wt{P}_\sigma}} {\rm d}\mu,  \end{gather}
 where $\mu$ is the  Lebesgue measure  in $\RR$ and $\vol$ is the corresponding  measure in $\SS^1$ with $\vol(\SS^1)=2 \pi$.  Due to  \eqref{GL tilde},  \eqref{star properties 1},  for any $z\in \CC$, any $g=(G,\psi)\in \wt{GL}^+(2,\RR)$,  and any $\sigma \in  \st(\mc T)$ we have: \begin{equation} \label{rotating P_sigma} P_{(z\star \sigma)}^{\mc T}= \exp(\ri \Im(z)) P_\sigma^{\mc T} \qquad  \wt{P}_{( \sigma\cdot g)}^{\mc T}= \psi^{-1} (\wt{P}_\sigma^{\mc T} ).\end{equation}

\begin{df} \label{arc} Let $0<\varepsilon <1$. Any subset of $\mathbb S^1$ of  the form $\exp(\ri \pi [a,a+\varepsilon])$, where $a\in \RR$ will be referred to as   a \textit{\textbf{closed $\varepsilon$-arc}} in $\mathbb S^1$.
\end{df}
\begin{remark} \label{phases and equivalences} The action of  ${\rm Aut}(\mc T)$  on $\st(\mc T)$ was recalled in the end of the previous section.  Following this definition one defines straightforwardly a biholomorphism $[F]:\st(\mc T_1)\rightarrow \st(\mc T_2)$ for any equivalence $F$  between triangulated categories $\mc T_1$ and  $\mc T_2$ satisfying $P_{[F](\sigma)}^{\mc T_2}=P_{\sigma}^{\mc T_1}$ for each $\sigma \in \st(\mc T_1)$. 
\end{remark}

In Definition \ref{main def}  we will use the following subset of the set of stability conditions: 
\begin{df} \label{definitin of st_varepsilon} For any $0<\varepsilon<1$ and any triangulated category $\mc T$ we denote:
\begin{gather} \st_{\varepsilon}(\mc T)=\{\sigma \in \st(\mc T): \sigma \ \ \mbox{is full and} \ \  {\mathbb S}^1\setminus P_\sigma^{\mc T} \ \ \mbox{contains a closed $\varepsilon$ arc}  \}\nonumber \\ \st_{[a,a+\varepsilon]}(\mc T) = \{\sigma \in \st(\mc T): \sigma \ \ \mbox{is full and} \ \   \wt{P}_\sigma^{\mc T}\cap [a,a+\varepsilon] = \emptyset  \}. \nonumber \end{gather}
\end{df}
 It is obvious that (recall also \eqref{rotating P_sigma}):
\begin{gather} \label{reduce to 0 varepsilon} \st_{\varepsilon}(\mc T) = \cup_{a\in \RR} \st_{[a,a+\varepsilon]}(\mc T)=\CC \star  \st_{[0,\varepsilon]}(\mc T)
  \end{gather}

The next  simple observation is: 
\begin{lemma} \label{for any varepsilon varepsilon prime} Let $g_{\varepsilon, \varepsilon'}\in \wt{GL}^+(2,\RR)$  be as in Remark \ref{an element in tilde}.  For any  $0<\varepsilon <1$,  $0<\varepsilon' <1$ holds: 
\begin{gather} \st_{[0,\varepsilon]}(\mc T)\cdot g_{\varepsilon, \varepsilon'} = \st_{[0,\varepsilon']}(\mc T). \end{gather}

\end{lemma}
\bpr   Using \eqref{rotating P_sigma}, \eqref{an element in tilde remark}, and the fact that $\psi$ is diffeomorphism  we compute $$\wt{P}_{(\sigma\cdot g_{\varepsilon,\varepsilon'})}^{\mc T}\cap [0,\varepsilon']= \psi^{-1}(\wt{P}_\sigma^{\mc T})\cap \psi^{-1}([0,\varepsilon])=\psi^{-1}\left (\wt{P}_\sigma^{\mc T}\cap [0,\varepsilon]\right ).$$ 
 Now the lemma follows from the very Definition \ref{definitin of st_varepsilon} and the property $g_{\varepsilon,\varepsilon'}^{-1}=g_{\varepsilon',\varepsilon}$.
 \epr

\begin{coro} \label{some equivalences}  Let $\mc T$ be any triangulated category. The following  are equivalent: 

{\rm (a)} $\st_{\varepsilon}(\mc T) \neq \emptyset$ for some $\varepsilon \in (0,1)$

{\rm (b)}  $\st_{\varepsilon}(\mc T) \neq \emptyset$ for each $\varepsilon \in (0,1)$

{\rm (c)} $P^{\mc T}_{\sigma} $ is not dense in $\mathbb S^1$ for some full $\sigma \in \st(\mc T)$. 

\end{coro}
\bpr (a) $\Rightarrow$ (b). Follows from \eqref{reduce to 0 varepsilon}  and Lemma \ref{for any varepsilon varepsilon prime}. 

(b) $\Rightarrow$ (c). It is obvious from the definitions that for any $0<\varepsilon<1$ and any $\sigma \in \st_{\varepsilon}(\mc T) $ the set $P_\sigma^{\mc T}$ is not dense in $\mathbb S^1$. 

(c) $\Rightarrow$ (a). If $P_{\sigma}^{\mc T}$ is not dense, then ${\mathbb S}^1\setminus P_{\sigma}^{\mc T}$ contains an open arc, but then it contains a closed arc  as well and then $\sigma \in \st_{\varepsilon}(\mc T) $ for some $\varepsilon \in (0,1)$.
\epr

\begin{df} \label{cat with phase gap} A triangulated category $\mc T$ will be called a \ul{category with phase gap} if  $P^{\mc T}_{\sigma} $ is not dense in $\mathbb S^1$  for some full $\sigma \in \st(\mc T)$ (by Corollary \ref{some equivalences} then $\st_{\varepsilon}(\mc T)$ is not empty for any $0<\varepsilon<1$). 
\end{df} 
\begin{lemma} \label{criteria for phase gap} If $K_0(\mc T)$ has finite rank, then $\mc T$ has a phase gap iff there exists a bounded t-structure in $\mc T$ whose heart is of finite length and has finitely many simple objects.    
\end{lemma}
\bpr

Let $\mc A$ be such a heart and let $s_1, s_2, \dots, s_n$ be  the simple objects in $\mc A$. Under the given assumptions $K_0(\mc T)\cong K_0(\mc A)\cong \ZZ^n$. \cite[Proposition 2.4, Proposition 5.3]{Bridg1} imply   that  for any sequence of vectors $z_1, z_2, \dots , z_n$ in $\HH$ there exists unique stability condition $\sigma = (Z, \mc P) $ with $\mc P(0,1]=\mc A$ and $Z(s_i)=z_i$, $i=1,\dots, n$. For this $\sigma$ we have $Z(\mc P(0,1]\setminus\{0\})=\{\sum_{i=1}^n a_i z_i: (a_1,a_2,\dots, a_n)\in \NN^n\setminus\{0\}\}$  and therefore $Z(\sigma^{ss})\subset \pm \{\sum_{i=1}^n a_i z_i: (a_1,a_2,\dots, a_n)\in \NN^n\setminus\{0\}\}$, now from  \cite[Lemma 1.1]{CP} it follows that $\sigma$ is locally finite and therefore, using the  notation explained after \eqref{sigma^{ss}}, we can write $\sigma \in \HH^{\mc A}$.  Varying the vector $(z_1,z_2,\dots,z_n)\in \HH^n$ we obtain a biholomorphism between $\HH^n$ and the  subset $\HH^{\mc A}\subset\st(\mc T)$. In particular the  stability conditions in $\HH^{\mc A}$ are full.   Since for $\sigma \in  \HH^{\mc A}$ holds   $Z(\sigma^{ss})\subset \pm \{\sum_{i=1}^n a_i z_i: (a_1,a_2,\dots, a_n)\in \NN^n\setminus\{0\}\}$,   it follows that $Z(\sigma^{ss})\subset \pm \{x \exp(\ri \pi a)+y \exp(\ri \pi (a+1-\varepsilon): x,y\in (0,+\infty) \}$ for some  $a\in \RR$ and some $0<\varepsilon<1$, therefore by  \eqref{phase formula} it follows that $\sigma \in \st_\varepsilon(\mc T)$, hence  $\mc T$ has a phase gap. 

Conversely,  let $\sigma'=(Z, \mc P)\in \st_\varepsilon(\mc T)$ be  a full stability condition, then due to \eqref{reduce to 0 varepsilon} for some $\lambda \in \CC $ we have $\sigma = \lambda \star \sigma'$ satisfying  $\mc  P_{\sigma}(t)=\{0\}$ for $t\in [0,\varepsilon]$ and  $\mc P_{\sigma}(0,1]= \mc P_{\sigma}(\varepsilon/2,1]$. From \cite[Lemma 4.5]{Bridg2}  it follows that $\mc P_{\sigma}(0,1]= \mc P_{\sigma}(\varepsilon/2,1]$ is a finite length quasi-abelian  category (here the property of $\sigma$  being full  is used), and since $\sigma$ is a stability condition,  $\mc P_{\sigma}(0,1]$  is a heart of a bounded t-structure, therefore  it is a finite length abelian  category whose   simple objects are a basis of $K_0(\mc T)$, in particular  the  simple objects are finitely many.
\epr 

\begin{remark} \label{algebraic } The elements  $\sigma \in \st (\mc T)$ for which $\mc P(0,1]$ is of finite length and with finitely many simple objects are called \ul{algebraic stability conditions} and have been discussed extensively in \cite{Woolf}. 
\end{remark} 
 \begin{remark}  Due to {\rm (i), (ii) } in the beginning of  \cite[Subsection 7.1]{BS} the   CY3 categories discussed in \cite{BS} have phase gaps. 
\end{remark} 

\begin{remark} \label{from exc coll to phase gap} Let  $\mc T$ be proper and with a full exceptional collection.  \cite[Remark 3.20]{DK1}  and Corollary \ref{some equivalences} imply that  $\st_\varepsilon(\mc T) \neq \emptyset$ for any $0<\varepsilon<1$, i. e. $\mc T$ is a category with a phase gap. \end{remark}   
  The main definition of this section is:
\begin{df} \label{main def} Let $\mc T$ be a triangulated category with phase gap. Let $0<\varepsilon<1$. We define:
\begin{equation} \label{main def eq} \norm{\mc T}_{\varepsilon} =  \sup \left  \{ \frac{1}{2} \vol\left (\ol{P_\sigma}\right ): \sigma \in \st_{\varepsilon}(\mc T)\right  \}. \end{equation}
\end{df} 
\begin{remark} \label{remark for dense categories} For a category $\mc T$ which carries a full stability condition, but has no phase gap (i. e. $P_{\sigma}^{\mc T}$ is dense in $\mathbb S^1$ for all full stability conditions $\sigma$)   it seems reasonable to define $\norm{\mc T}_{\varepsilon}=  \pi(1-\varepsilon)$, but we will restrict our attention to categories with phase gaps in the rest.
\end{remark}
In   remarks \eqref{norms with integral}, \eqref{remark leq pi(1-varepsilon)}  $\varepsilon$ and $\mc T$ are  as in Definition \ref{main def}.
\begin{remark} \label{norms with integral} Using \eqref{translate the volume}, \eqref{rotating P_sigma}, \eqref{reduce to 0 varepsilon}  one shows that ( $\mu$ is the Lebesgue measure of $\RR$):
\begin{equation} \label{main def eq int} \norm{\mc T}_{\varepsilon} =  \sup \left  \{\pi 
  \mu \left ( [\varepsilon ,1]\cap \ol{\wt{P}_\sigma} \right ) :\sigma \in \st_{[0,\varepsilon]}(\mc T) \right \}. \end{equation}
\end{remark}

\begin{remark} \label{remark leq pi(1-varepsilon)}  We have always $ 0\leq \norm{\mc T}_{\varepsilon}\leq  \pi(1-\varepsilon) $. \end{remark}

\begin{remark} \label{equal norms for isomorphic} Using Remark \ref{phases and equivalences} we see that if $\mc T_1$, $\mc T_2$ are equivalent triangulated categories  with finite rank Grothendieck groups, then for any $0<\varepsilon<1$
holds  $\norm{{\mc T}_{1}}_\varepsilon=\norm{{\mc T}_{2}}_\varepsilon$.  \end{remark}
\begin{lemma} \label{maximal norms}  Let $\varepsilon$, $\varepsilon'$ be any two numbers in $(0,1)$. 

{\rm (a)} There exist  $0<m<M$ such that $m \norm{\mc T}_{\varepsilon}\leq \norm{\mc T}_{\varepsilon'}\leq M \norm{\mc T}_{\varepsilon}$ for any category with a phase gap $\mc T$. In particular,  for any category with a phase gap $\mc T$ we have: $\norm{\mc T}_{\varepsilon}=0 \iff \norm{\mc T}_{\varepsilon'}=0$
  .

{\rm (b)} For  any category with a phase gap $\mc T$ we have  $\norm{\mc T}_{\varepsilon}=\pi (1-\varepsilon)\iff \norm{\mc T}_{\varepsilon'}=\pi (1-\varepsilon')$.

\end{lemma}
\bpr We will use the element $g_{\varepsilon,\varepsilon'}=(G,\psi)\in \wt{GL}^+(2,\RR)$ from Remark \ref{an element in tilde}. In particular the function $\psi \in C^{\infty}(\RR)$ restricts to a diffeomorphsim  $\psi_{\vert}:[\varepsilon',1]\rightarrow [\varepsilon,1]$. Let us denote the inverse function by $\kappa$, then we choose   $m, M\in \RR$ as follows:
\begin{gather} \label{maximal norms formula 1} \psi_{\vert}^{-1}=\kappa:[\varepsilon,1]\rightarrow [\varepsilon',1] \ \ \forall t \in [\varepsilon, 1] \ \ \  0<m\leq\kappa'(t)\leq M.\end{gather} 

 With the help of  \cite[formula (15) on page 156]{Rudin},  we  see that for any Lebesgue measurable  subset $A\subset [\varepsilon,1]$ holds (for a subset $E\subset [\varepsilon',1]$ or $E\subset [\varepsilon,1]$ we denote by $\chi_E$ the function equal to $1$ at the points of $E$ and $0$  elsewhere):
\begin{gather}  \mu \left (\kappa\left ( A \right )\right ) = \int_{\varepsilon'}^1 \chi_{\kappa\left ( A \right )}(t) {\rm d} t = \int_{\varepsilon}^1 \chi_{\kappa\left ( A\right )}( \kappa ( t)) \kappa'(t) {\rm d} t  = \int_{\varepsilon}^1 \chi_{ A }( t) \kappa'(t) {\rm d} t \nonumber
  \end{gather} which by  \eqref{maximal norms formula 1} implies:
\begin{gather}\label{measures with kappa} m \mu \left ( A \right ) \leq \mu \left (\kappa\left ( A \right )\right ) \leq M \mu \left ( A \right ).
  \end{gather}
 Using Remark \ref{norms with integral}, Lemma \ref{for any varepsilon varepsilon prime}, and the second equality in  \eqref{rotating P_sigma}  we get: 
\begin{gather}  \label{maximal norms formula 21}  \norm{\mc T}_{\varepsilon}/\pi =  \sup \left  \{  \mu \left ( [\varepsilon ,1]\cap \ol{\wt{P}_{\sigma}}  \right ):\sigma \in \st_{[0,\varepsilon]}(\mc T) \right \}\\
 \label{maximal norms formula 2} \norm{\mc T}_{\varepsilon'}/\pi =  \sup \left  \{  \mu \left (\kappa \left ( [\varepsilon ,1]\cap \ol{\wt{P}_{\sigma}} \right ) \right ):\sigma \in \st_{[0,\varepsilon]}(\mc T) \right \}. \end{gather} 
Now (a) follows from \eqref{measures with kappa}, \eqref{maximal norms formula 21},  \eqref{maximal norms formula 2}.

 (b) Let $\norm{\mc T}_{\varepsilon}=\pi (1-\varepsilon)$ and $\delta>0$.  We will prove that \eqref{maximal norms formula 2} equals $(1-\varepsilon')$ by finding $\sigma \in \st_{[0,\varepsilon]}(\mc T)$ such that $\mu \left (\kappa \left ( [\varepsilon ,1]\cap \ol{\wt{P}_{\sigma}} \right )\right )>1-\varepsilon'-\delta$. Since 
$ 1-\varepsilon' = \mu ([\varepsilon',1]) = \mu \left (\kappa\left ( [\varepsilon ,1]\cap \ol{\wt{P}_{\sigma}} \right )\right ) + \mu \left (\kappa\left ( [\varepsilon ,1]\setminus \ol{\wt{P}_{\sigma}} \right )\right ) $, we need to find $\sigma \in \st_{[0,\varepsilon]}(\mc T)$ such that: \begin{gather} \label{maximal norms part 3} \mu \left (\kappa\left ( [\varepsilon ,1]\setminus \ol{\wt{P}_{\sigma}} \right )\right ) < \delta. \end{gather} Since $\norm{\mc T}_{\varepsilon}=\pi (1-\varepsilon)$, \eqref{maximal norms formula 21} ensures that there is $\sigma \in \st_{[0,\varepsilon]}(\mc T)$ such that $\mu\left ([\varepsilon,1]\cap \ol{\wt{P}_\sigma}\right )> 1-\varepsilon - \frac{\delta}{M}$, which due to the equality  $\mu\left ([\varepsilon,1]\cap \ol{\wt{P}_\sigma}\right )+\mu\left ([\varepsilon,1]\setminus \ol{\wt{P}_\sigma}\right ) = 1-\varepsilon$ is the same as 
$  \mu\left ([\varepsilon,1]\setminus \ol{\wt{P}_\sigma}\right ) < \frac{\delta}{M} . $
  We combine \eqref{measures with kappa} and the latter inequality to deduce  the desired \eqref{maximal norms part 3}:
	$  \mu \left (\kappa\left ( [\varepsilon ,1]\setminus \ol{\wt{P}_{\sigma}} \right )\right )\leq M  \mu\left ([\varepsilon,1]\setminus \ol{\wt{P}_\sigma}\right )<\delta.\nonumber $
\epr

\cite[Corollary 3.28]{DHKK} (see \cite[Corollary 3.25]{Dimitrov} for general algebraically closed field $k$) amounts to the following criteria for non-vanishing of $\norm{\mc T}_{\varepsilon}$ 
\begin{prop} \label{from kronecker to nonvan norm}
 Let $(E_0,E_1,\dots,E_n)$ be a full exceptional collection in a $k$-linear proper triangulated category 
$\mc D$. If for some $i$ the pair    $(E_i,E_{i+1})$  satisfies  $\hom^1(E_i,E_{i+1})\geq 3$ and  $\hom^{\leq 0}(E_i,E_{i+1})=0$, then 
$\norm{\mc D}_{\varepsilon} >0$. 
\end{prop}
\begin{coro} \label{norm on wild quivers} Let  $\varepsilon\in (0,1)$. Then:

{\rm (a) } If $Q$ is an acyclic quiver, which is neither Dynkin nor affine, then $\norm{D^b(Q)}_{\varepsilon}>0$. 

{\rm (b)} $\norm{D^b(coh(X))}_{\varepsilon}>0$, where $X$ is  a smooth projective variety over $\CC$, such that $D^b(coh(X))$ is generated by a strong exceptional collection of
three elements
\end{coro}
\bpr  (a) Follows from the previous proposition, \cite[Proposition 3.34]{DHKK}, and the fact that each exceptional collection in $D^b(Q)$ can be extended to a full exceptional collection (see \cite{WCB1}).

(b) It follows from  proposition \ref{from kronecker to nonvan norm} and \cite[3.5.1]{DHKK}.
\epr

In Section \ref{criteria} we will refine Proposition \ref{from kronecker to nonvan norm}, which will help us   to prove that 
 $\norm{D^b(coh(X)}_{\varepsilon}=\pi(1-\varepsilon)$ if $X$ is $\PP^1 \times \PP^1$, $\PP^n$ with $n\geq 2$ or some of these  blown  up in finite number of points.

\begin{prop} \label{norm of Dynkin} Let  $\varepsilon \in (0,1)$.  For acyclic  quiver $Q$ we have $\norm{D^b(Q)}_{\varepsilon} = 0$ iff $Q$ is affine or Dynkin.  In particular  $\norm{D^b(coh(\PP^1))}_{\varepsilon}=0$.  
\end{prop}
\bpr If $Q$ is affine or Dynkin, then  from the first and the second raws of table \eqref{table_intro}  we see that  $\vol\left (\ol{P_\sigma}\right )=0$   for any $\sigma \in \st(D^b(Q))$, therefore   $\norm{D^b(Q)}_{\varepsilon} = 0$, and in Corollary \ref{norm on wild quivers}  we showed that  $\norm{D^b(Q)}_{\varepsilon} > 0$ for the rest quivers. 
\epr

\section{Stability conditions on  orthogonal decompositions} \label{st cond on Orth Dec}

First we recall the  definition of a semi-orthogonal, resp. orthogonal, decomposition  of a triangulated category:

\begin{df}\label{SOD}  If $\mc T$ is a  triangulated category,  $\mc T_1$,  $\mc T_2$, $\dots$,  $\mc T_n$  are triangulated subcategories in it  satisfying the equalities   $\mc T = \langle \mc T_1, \mc T_2,\dots \mc T_n \rangle$ and  $\Hom(\mc T_j, \mc T_i)=0$ for $j>i$, then we say that   $\mc T = \langle \mc T_1, \mc T_2,\dots \mc T_n \rangle $ is a \ul{semi-orthogonal decomposition}.  If in addition holds  $\Hom(\mc T_i,\mc T_j)=0$ for $i<j$, then we say that $\mc T = \langle \mc T_1, \mc T_2,\dots \mc T_n \rangle  $ is an \ul{orthogonal decomposition}, in which case we will write sometimes $\mc T=\mc T_1 \oplus \mc T_2\oplus \dots 
\oplus \mc T_n$. 
Obviously, if  $\mc T = \langle \mc T_1, \mc T_2,\dots \mc T_n \rangle  $ is an orthogonal decomposition, then  $\mc T = \langle \mc T_{s(1)}, \mc T_{s(2)},\dots \mc T_{s(n)} \rangle  $ is an orthogonal decomposition for any permutation $s:\{ 1,\dots,n \} \rightarrow \{ 1,\dots,n \}$.\end{df}

\begin{prop} \label{lemma for orthogonal composition 1} Let  $\mc T = \langle \mc T_1, \mc T_2, \dots, \mc T_n \rangle$ be any orthogonal decomposition.  Let\\ $\bd K_0(\mc T_i) & \rTo^{in_i} & K_0(\mc T) & \rTo^{pr_j} & K_0(\mc T_j) \ed $, $1\leq i,j \leq n$ be the natural biproduct diagram.  Then:

{\rm (a)} The following map is a  bijection:
\begin{gather} \label{map for products1}  \st(\mc T) \rightarrow \st(\mc T_1)\times \st(\mc T_2) \times \cdots \times \st(\mc T_n)\\  (Z,\{ \mc P(t) \}_{t \in \RR}) \mapsto \left ( (Z\circ pr_1,\{ \mc P(t)\cap \mc T_1 \}_{t \in \RR}), \dots,  (Z\circ pr_n,\{ \mc P(t)\cap \mc T_n \}_{t \in \RR}) \right ).\end{gather}

{\rm (b) } For any  $(Z,\{ \mc P(t) \}_{t \in \RR})\in  \st(\mc T)$ and any $t\in \RR$ the subcategory $\mc P(t)$ is non-trivial iff for some $j$ $\mc P(t)\cap \mc T_j$ is non-trivial.

{\rm (c) } If ${\rm rank}(K_0(\mc T_i))<+\infty$ for all $i=1,2,\dots,n$, then the map defined above is biholomorphism.

{\rm (d)} For each $\sigma \in \st(\mc T)$ holds $P_\sigma^{\mc T}=\cup_{i=1}^n P_{\sigma_i}^{\mc T_i}  $, where $(\sigma_1,\dots,\sigma_n)$ is the value of \eqref{map for products1} at $\sigma$.
\end{prop}
\bpr We will give all details for the proof of (a), (b), (c)  in the case $n=2$. The general case follows easily by induction. (d) follows from the very definition \eqref{the set of phases def} and (a), (b). 

 It is well known that for each $X\in \mc T$ there exists unique up to isomorphism triangle $E_2 \rightarrow X \rightarrow E_1 \rightarrow E_2[1]$ with $E_i\in \mc T_i$, $i=1,2$. By $\Hom(\mc T_1, \mc T_2)=0$ it follows that each of these triangles is actually part of a direct product diagram and $pr_i([X])=[E_i]$ for $i=1,2$. 

Now let $X\in \mc T_1$ and  $U\rightarrow X \rightarrow B \rightarrow U[1]$ be a triangle in $\mc T$. Using $\Hom(\mc T_2, \mc T_1)=\Hom(\mc T_1, \mc T_2)=0$ and decomposing $U$ into direct summands $U_1 \oplus U_2$ with $U_i\in \mc T_i$ one easily  concludes that the triangle  $U\rightarrow X \rightarrow B \rightarrow U[1]$  is isomorphic to a  triangle of the form $U_1 \oplus U_2 \rightarrow X \rightarrow B'\oplus  U_2[1] \rightarrow U_1[1] \oplus U_2 [1]$. If we apply these arguments to the last triangle in \eqref{HN filtration} and using that $\hom(E_{n-1}, A_n[i])=0$ for $i\leq 0$, we immediately obtain $E_{n-1}, A_n \in \mc T_1$ and then by induction it follows that the entire HN filtration of $X$ lies in $\mc T_1$, in particular $A_i \in \mc P(t_i) \cap \mc T_1$ for $i=1,2,\dots,n$, furthermore we have $Z_1 ([X]) =Z (pr_1([X])) $ for each $X \in \mc P(t) \cap \mc T_1$ and now it is obvious that $ (Z\circ pr_1,\{ \mc P(t)\cap \mc T_1 \}_{t \in \RR})=(Z_1, \mc P_1)$ is a stability condition on $\mc T_1$.  The same arguments apply to the  case  $X \in \mc T_2$ and show that   $ (Z\circ pr_2,\{ \mc P(t)\cap \mc T_2 \}_{t \in \RR})=(Z_2, \mc P_2)$  is a stability condition on  $\mc T_2$. We will  show that $\sigma_i$ are locally finite for $i=1,2$. Indeed, since $\sigma$ is locally finite  stability condition on $\mc T$, then there exists $\frac{1}{2} > \varepsilon >0$ such that  $\mc P(t-\varepsilon,t+\varepsilon)$ is quasi-abelian category of finite length for each $t\in \RR$. One easily shows that  $\mc P_i(t-\varepsilon,t+\varepsilon)=\mc T_i \cap \mc P(t-\varepsilon,t+\varepsilon)$ for each $t$. From \cite[Lemma 4.3]{Bridg1}   we know that a sequence $A\rightarrow B \rightarrow C$ in  $\mc P_i(t-\varepsilon,t+\varepsilon)$ is a strict short exact sequence iff it is part of a triangle  $A\rightarrow B \rightarrow C \rightarrow A[1]$  in $\mc T_i$. Since for $A,B,C$ in $\mc T_i$    $A\rightarrow B \rightarrow C \rightarrow A[1]$ is triangle in $\mc T_i$ iff it is  a triangle in $\mc T$, we deduce that for  $A, B,  C \in \mc P_i(t-\varepsilon,t+\varepsilon)$  $A\rightarrow B \rightarrow C $  is a  strict exact  sequence in $\mc P_i(t-\varepsilon,t+\varepsilon)$ iff it is a strict exact sequence in   $\mc P(t-\varepsilon,t+\varepsilon)$, and now from the fact that  $\mc P(t-\varepsilon,t+\varepsilon)$ is of finite length it follows that  $\mc P_i(t-\varepsilon,t+\varepsilon)$ is of finite length and $\sigma_i \in \st(\mc T_i)$ for $i=1,2$. Thus the map is well defined. Since for any interval $I\subset \RR$ the subcategory $\mc P(I)$ is thick (see e.g. \cite[Lemma 2.20.]{Dimitrov}), it follows that $\mc P(t)=\mc P_1(t)\oplus \mc P_2(t)$ for each $t \in \RR$ and hence follows the injectivity of the map. Furthermore, using the terminology of \cite[Definition before Proposition 2.2]{CP} we see that $\sigma$ is glued from $\sigma_1$ and $\sigma_2$. From the given arguments it follows also that for $X\in \mc T_i$ the HN filtrations w.r. to $\sigma$ and w.r. to $\sigma_i$ coincide, in particular:
\begin{gather}\label{help for the metrics}  X\in \mc T_i \Rightarrow \phi_{\sigma_i}^{\pm}(X)=\phi_{\sigma}^{\pm}(X) \qquad m_{\sigma_i}(X)=m_{\sigma}(X)
\end{gather}
on the other hand any $X\in \mc T$ can be represented uniquely (up to isomorphism) as a biproduct $X\cong X_1\oplus X_2 $ with $X_i\in \mc T_i$ for $i=1,2$ and \eqref{direct sums} imply
\begin{gather} \label{help for metrics} X\in \mc T \Rightarrow X\cong X_1\oplus X_2, \ \ X_i \in \mc T_i \Rightarrow \begin{array}{c} m_\sigma(X) = m_{\sigma_1}(X_1) + m_{\sigma_2}(X_2)  \\  \phi_\sigma^-(X)={\rm min}\{ \phi_{\sigma_1}^-(X_1),  \phi_{\sigma_2}^-(X_2)\} \\  \phi_\sigma^+(X)={\rm max}\{ \phi_{\sigma_1}^+(X_1),  \phi_{\sigma_2}^+(X_2)\} \end{array} \end{gather}

Conversely, if $(\sigma_1,\sigma_2) \in \st(\mc T_1)\times \st(\mc T_2)$, then \cite[Proposition 3.5]{CP}  ensures existence of a locally finite  stability condition $\sigma \in \st(\mc T)$ glued from $\sigma_1$, $\sigma_2$ and using  \cite[(3) in Proposition 2.2]{CP}) one easily shows  that our map sends the glued $\sigma$ to the pair $(\sigma_1,\sigma_2)$, hence the surjectivity of the map follows. 

Now we will show that if   ${\rm rank}(K_0(\mc T_i))<+\infty$ for  $i=1,2$, then the map defined above is biholomorphism. 
First we show that it is continuous. We denote by  $d$, $d_1$, $d_2$ the generalized  metrics on $\st(\mc T)$, $\st(\mc T_1)$, $\st(\mc T_2)$ (as defined in \eqref{metric}). For any $\sigma, \sigma'\in \st(\mc T)$ let  $(\sigma_1,\sigma_2)$ and $(\sigma'_1,\sigma'_2)$ be the pairs assigned via the bijection.  To show that the  map is homeomorphism we will show that :
\begin{gather} \label{metrics1} \max\{d_1(\sigma_1,\sigma'_1), d_2(\sigma_2,\sigma'_2)\}  \leq d(\sigma,\sigma') \\  \label{metrics2}  d(\sigma,\sigma')\leq d_1(\sigma_1,\sigma'_1)+ d_2(\sigma_2,\sigma'_2)\end{gather}
The first \eqref{metrics1} follows easily from \eqref{help for the metrics}. The second requires a bit more computations, which we will present partly. Take any $X\in \mc T$ and decompose it $X \cong X_1 \oplus X_2$, $X_i\in \mc T_i$, then from \eqref{help for metrics} we see that 
\begin{gather} \abs{\log\frac{m_\sigma(X)}{m_{\sigma'}(X)}}= \abs{\log\frac{m_{\sigma_1}(X_1)+m_{\sigma_2}(X_2)}{m_{\sigma'_1}(X_1)+m_{\sigma'_2}(X_2)}}\leq \abs{\log\frac{m_{\sigma_1}(X_1)}{m_{\sigma'_1}(X_1)}}+\abs{\log\frac{m_{\sigma_2}(X_2)}{m_{\sigma'_2}(X_2)}}\\ \leq  d_1(\sigma_1,\sigma'_1)+ d_2(\sigma_2,\sigma'_2),\end{gather} where we used, besides the definition of the generalized metrics \eqref{metric}, the following lemma:
\begin{lemma} For any positive  real numbers $x_1,x_2,y_1,y_2$ holds the inequality:
 $$\abs{\log\frac{x_1+x_2}{y_1+y_2}}\leq \abs{\log\frac{x_1}{y_1}}+\abs{\log\frac{x_2}{y_2}}.$$
\end{lemma}
\bpr We can assume that $\frac{x_1+x_2}{y_1+y_2}\geq 1$ (otherwise take $\frac{y_1+y_2}{x_1+x_2}$). Now we consider three cases:

If $\frac{x_1}{y_1}\geq 1$ and  $\frac{x_2}{y_2}\geq 1$, then the desired inequality becomes 
$ \log\frac{x_1+x_2}{y_1+y_2}\leq \log\frac{x_1}{y_1}+\log\frac{x_2}{y_2}$
 which after exponentiating is equivalent to 
 \begin{gather} \frac{x_1+x_2}{y_1+y_2}\leq \frac{x_1 x_2}{y_1 y_2} \Longleftrightarrow  (x_1+x_2)y_1 y_2 \leq x_1 x_2 (y_1+y_2) \Longleftrightarrow0 \leq   x_1 y_1(x_2 -y_2)+x_2 y_2(x_1 -y_1)\nonumber \end{gather}
the latter inequality follows from $x_1\geq y_1$,  $x_2\geq y_2$.

 If $\frac{x_1}{y_1}\leq 1$ and  $\frac{x_2}{y_2}\geq 1$, then the desired inequality becomes 
$ \log\frac{x_1+x_2}{y_1+y_2}\leq \log\frac{y_1}{x_1}+\log\frac{x_2}{y_2} $
 which after exponentiating is equivalent to
 \begin{gather} \frac{x_1+x_2}{y_1+y_2}\leq \frac{y_1 x_2}{x_1 y_2} \Longleftrightarrow  (x_1+x_2)x_1 y_2 \leq y_1 x_2 (y_1+y_2) \Longleftrightarrow 0 \leq   y_1^2 x_2-x_1^2 y_2+ x_2 y_2(y_1 -x_1)\nonumber \end{gather}
the latter inequality follows from $y_1\geq x_1$,  $x_2\geq y_2$.

 If $\frac{x_1}{y_1}\leq 1$ and  $\frac{x_2}{y_2}\leq 1$, then the desired inequality becomes 
$ \log\frac{x_1+x_2}{y_1+y_2}\leq \log\frac{y_1}{x_1}+\log\frac{y_2}{y_2} $
 which after exponentiating is equivalent to
 \begin{gather} \frac{x_1+x_2}{y_1+y_2}\leq \frac{y_1 y_2}{x_1 x_2} \Longleftrightarrow  (x_1+x_2)x_1 x_2 \leq y_1 y_2 (y_1+y_2) \Longleftrightarrow 0 \leq   y_1^2 y_2-x_1^2 x_2+y_2^2 y_1-x_2^2 x_1 \nonumber \end{gather}
the latter inequality follows from $y_1\geq x_1$,  $y_2\geq x_2$.
\epr
Now in order to prove \eqref{metrics2}  it is enough to show that      $ \abs{\phi_\sigma^\pm(X)-\phi_{\sigma'}^\pm(X)}\leq d_1(\sigma_1,\sigma'_1)+ d_2(\sigma_2,\sigma'_2) $ which in turn via \eqref{help for metrics} is the same as 
\begin{gather} \abs{{\rm max}\{ \phi_{\sigma_1}^+(X_1),  \phi_{\sigma_2}^+(X_2)\} -  {\rm max}\{ \phi_{\sigma'_1}^+(X_1),  \phi_{\sigma'_2}^+(X_2)\} } \leq  d_1(\sigma_1,\sigma'_1)+ d_2(\sigma_2,\sigma'_2) \\
 \abs{{\rm min}\{ \phi_{\sigma_1}^-(X_1),  \phi_{\sigma_2}^-(X_2)\} -  {\rm min}\{ \phi_{\sigma'_1}^-(X_1),  \phi_{\sigma'_2}^-(X_2)\} } \leq  d_1(\sigma_1,\sigma'_1)+ d_2(\sigma_2,\sigma'_2),\end{gather}
  which in turn follow from the following: 
\begin{lemma} For any real numbers $x_1,x_2,y_1,y_2$ we have:
\begin{gather} \abs{{\rm max}\{x_1,x_2\}-{\rm max}\{y_1,y_2\}} \leq \abs{x_1-y_1}+\abs{x_2-y_2} \nonumber \\
 \abs{{\rm min}\{x_1,x_2\}-{\rm min}\{y_1,y_2\}} \leq \abs{x_1-y_1}+\abs{x_2-y_2} \nonumber \end{gather}
\end{lemma}
\bpr If ${\rm max}\{x_1,x_2\}=x_i$ and  ${\rm max}\{y_1,y_2\}=y_i$ for the same $i$, then the inequalities  follow immediately. So let ${\rm max}\{x_1,x_2\}=x_i$  ${\rm max}\{y_1,y_2\}=y_j$, $i\neq j$, e.g. let $i=1$, $j=2$. Then $x_1\geq x_2$, $y_1\leq y_2$, and the lemma follows from:
\begin{gather} \abs{{\rm max}\{x_1,x_2\}-{\rm max}\{y_1,y_2\}}=  \abs{x_1-y_2}=\\ 
= \left \{\begin{array}{c c}  x_1-y_2= x_1-y_1+y_1-y_2\leq  x_1-y_1=\abs{ x_1-y_1} & \mbox{if} \ x_1\geq y_2\nonumber \\ y_2-x_1= y_2-x_2+x_2-x_1 \leq  y_2-x_2=\abs{ x_2-y_2} & \mbox{if} \ x_1\leq y_2\end{array} \right. \nonumber \end{gather}
\begin{gather} \abs{{\rm min}\{x_1,x_2\}-{\rm min}\{y_1,y_2\}}=  \abs{x_2-y_1}=\\ 
= \left \{\begin{array}{c c}  x_2-y_1=x_2-x_1+x_1-y_1\leq  x_1-y_1=\abs{ x_1-y_1} & \mbox{if} \ x_2\geq y_1 \nonumber \\ y_1-x_2= y_1-y_2+y_2-x_2  \leq  y_2-x_2=\abs{ x_2-y_2} & \mbox{if} \  x_2\leq y_1\end{array} \right. .\nonumber \end{gather}
\epr
Thus, we have \eqref{metrics1}, \eqref{metrics2} and they imply that \eqref{map for products1}  is homeomorphism for $n=2$.

Let $\bd \st(\mc T)&\rTo^{proj}&\Hom(K_0(\mc T),\CC)\ed$, $\bd \st(\mc T_i)&\rTo^{proj_i}&\Hom(K_0(\mc T_i),\CC)\ed$, $i=1,2$ be the projections  $proj(Z, \mc P)=Z$. Then the following diagram (the first row is  the map \eqref{map for products1} and the second row is the assignment $Z\mapsto(Z\circ pr_1,Z\circ pr_2)$)  is commutative: \bd[size=1.5em] \st(\mc T) & \rTo^{\varphi}                & \st(\mc T_1)\times \st(\mc T_2)\\
                    \dTo^{proj}               &                           &    \dTo^{proj_1\times proj_2}                         \\
												\Hom(K_0(\mc T),\CC)          &   \rTo^{\varphi'} & 	\Hom(K_0(\mc T_1),\CC)  \times 	\Hom(K_0(\mc T_2),\CC).   \ed   
												If we take any connected component $\Sigma \subset \st(\mc T)$, then (since $\varphi$ is homeomorphism) $\varphi(\Sigma)=\Sigma_1 \times \Sigma_2$ is a connected component of $\st(\mc T_1)\times \st(\mc T_2)$, resp. $\Sigma_i$ are connected components of $\st(\mc T_i)$, and furthermore $m=\dim_\CC(\Sigma)=\dim_\CC(\Sigma_1\times \Sigma_2)$.  From the Bridgeland's main theorem we know that $proj$  restricts to local biholomorphisms between $\Sigma$  and an $m$-dimensional  vector subspace $V\subset \Hom(K_0(\mc T),\CC) $  and $proj_1\times proj_2$ restricts to local biholomorphisms between $\Sigma_1\times \Sigma_2$  and an $m$-dimensional  vector subspace   $V_1\times V_2 \subset \Hom(K_0(\mc T_1),\CC) \times  \Hom(K_0(\mc T_2),\CC) $. It follows (using that $\varphi'$ is a linear isomorphism and that each open subset in a vector subset contains a basis of the space) that $\varphi'(V)=V_1\times V_2$. Thus, the diagram above restricts to a diagram with vertical arrows which are local biholomorphisms, the bottom arrow is biholomorphism,   and the top arrow is a homeomorphism, it follows with standard arguments  that the the top arrow  must be biholomorphic. It follows that $\varphi$ is biholomorphism  and  we proved the proposition. 
 \epr

From this proposition and Definition \ref{definitin of st_varepsilon}  it  follows:
\begin{coro} \label{gaps in products}  Let $\mc T =  \mc T_1 \oplus \mc T_2 \oplus \dots \oplus \mc T_n $ be an orthogonal decomposition (Definition \ref{SOD}) and let ${\rm rank}(K_0(\mc T_i))<+\infty$ for $i=1,\dots,n$.    Let $ \st(\mc T) \rightarrow \st(\mc T_1)\times \cdots \times   \st(\mc T_n),  \sigma \mapsto \left ( \sigma_1, \sigma_2,\dots,\sigma_n \right )$ be the biholomorphism  from Proposition  \ref{lemma for orthogonal composition 1} {\rm (a)}. 

For any $0<\varepsilon<1$ the following are equivalent:  {\rm (a)} $\sigma \in \st_{\varepsilon}(\mc T)$;  {\rm (b)}  $\{ \sigma_i \in \st_{\varepsilon}(\mc T_i) \}_{i=1}^n$ and there exists  a closed $\varepsilon$-arc $\gamma$ such that $P_{\sigma_i}^{\mc T_i}\cap \gamma= \emptyset$ for each $1\leq i \leq n$. 

In particular $\mc T$ has a phase gap iff   $\mc T_i$ has a phase gap for each $1\leq i \leq n$.
\end{coro}
Since the closure of $A \cup B$ equals the union of  closures of $A$ and $B$ and $\vol(A)\leq \vol(A\cup B)\leq \vol(A)+\vol(B)$, from Corollary \ref{gaps in products} it follows:
\begin{coro} \label{sum of several} Let  $\mc T ={\mc T}_{1}\oplus{\mc T}_{2}\oplus\dots\oplus{\mc T}_{n}$  be an orthogonal decomposition with finite rank Grothendieck groups of the factors,  and  let $0<\varepsilon<1$.

If  $\mc T$ has a phase gap and  $\norm{\mc T_j}_{\varepsilon}=0$ for some $j$, then $\norm{\mc T}_{\varepsilon} = \norm{ \langle {\mc T}_{1},{\mc T}_{2},{\mc T}_{j-1}, {\mc T}_{j+1},\dots,{\mc T}_{n}\rangle}_{\varepsilon}$.
\end{coro}

\section{ The inequality \texorpdfstring{$\norm{\langle \mc T_1,\mc T_2\rangle }_{\varepsilon}\geq \max\{\norm{\mc T_1}_{\varepsilon},\norm{\mc T_2}_{\varepsilon} \}$ }{\space} }

Here we  show  conditions which ensure $\norm{\langle \mc T_1,\mc T_2\rangle }_{\varepsilon}\geq \max\{\norm{\mc T_1}_{\varepsilon},\norm{\mc T_2}_{\varepsilon} \}$ for any $\varepsilon \in (0,1)$, where $ \langle \mc T_1, \mc T_2 \rangle$ is  a semi-orthogonal decomposition (see Definition \ref{SOD}) of some $\mc T$.

\begin{theorem} \label{prop inequality} Let  $\mc T$ be proper  and let $K_0(\mc T)$ has finite rank. Assume $0<\varepsilon <1$. 
Let $\mc T = \langle \mc T_1, \mc T_2 \rangle$ be  a semi-orthogonal decomposition.  If  $\mc T_1$, $\mc T_2$ are categories with phase gaps,  then  $\mc T$ is a category with phase gap and for any $0<\varepsilon <1$ holds : 
\begin{gather}\label{main inequality} \norm{\langle \mc T_1, \mc T_2 \rangle}_{\varepsilon} \geq  \max\left \{ \norm{\mc T_1}_{\varepsilon
}, \norm{\mc T_2}_{\varepsilon
} \right \}.\end{gather}
\end{theorem}
\bpr
Take any $0<\mu$.  Let $\sigma_i=(Z_i, \mc P_i)\in \st_\varepsilon(\mc T_i)$ be full stability conditions, s. t. $\frac{\vol\left (\ol{P^{\mc T_i}_{\sigma_i}}\right )}{2} >\norm{\mc T_i}_{\varepsilon
}-\mu$ for $i=1,2$. Due to \eqref{reduce to 0 varepsilon}  we can assume that $\exp(\ri \pi [0,\varepsilon]) \subset \SS^1\setminus P^{\mc T_i}_{\sigma_i}$. By the same arguments as in the last paragraph of the proof of Lemma \ref{criteria for phase gap} it follows that  $\mc P_{\sigma_i}(0,1]$ are finite length abelian categories, therefore   the simple objects in them are a basis of $K_0(\mc T_i)$ for $i=1,2$, and these abelian categories are the extension closures of their simple objects.  In particular  the sets of simple objects are finite  and it follows that for some $j\in \ZZ$ $\Hom^{\leq 1} (\mc  P_{\sigma_1}(0,1],\mc  P_{\sigma_2}(0,1][j])=\Hom^{\leq 1} (\mc  P_{\sigma_1}(0,1],\mc  P_{\sigma_2}(j,j+1])=0$. Recalling \eqref{star} we deduce that $\Hom^{\leq 1} (\mc  P_{\sigma_1}(0,1],\mc  P_{(-\ri j \pi)\star\sigma_2}(0,0+1])=0$.  By replacing  $\sigma_2$ with   $(-\ri j \pi)\star\sigma_2$ we  obtain    stability conditions  $\sigma_i\in \st_\varepsilon(\mc T_i)$ for $i=1,2$ satisfying the following conditions:

 \begin{gather} \label{first step in main ineq} \frac{\vol\left (\ol{P^{\mc T_i}_{\sigma_i}}\right )}{2} >\norm{\mc T_i}_{\varepsilon
}-\mu \quad  \mbox{for}  \quad i=1,2, \\  
\Hom^{\leq 1} (\mc  P_{\sigma_1}(0,1],\mc  P_{\sigma_2}(0,1])=0,   \\ 
\label{are of finite length} \mc  P_{\sigma_2}(0,1]  \quad \mbox{and}   \quad \mc  P_{\sigma_2}(0,1] \quad  \mbox{are of finite length and with f.m. simples}, \\
\label{help formula 4} \mc  P_{\sigma_i}(t)=\{0\}\quad \mbox{for}   \quad t\in [j,j+\varepsilon] \quad \mbox{for}   \quad i=1,2,  j\in \ZZ.
\end{gather}

In the listed  properties of $\sigma_i \in \st(\mc T_i)$ with  the given semi-orthogonal decomposition $\mc T=\langle \mc T_1, \mc T_2 \rangle$   are contained the conditions of  \cite[Proposition 3.5 (b)]{CP}. This proposition  ensures a  glued (see \cite[Definition ]{CP}) locally finite stability condition $\sigma =(Z,\mc P)\in \st(\mc T)$.  The glued stability condition satisfies the following (we use \cite[Proposition 2.2 (3)]{CP} and write   $\mc P_i$ instead of $\mc P_{\sigma_i}$)
\begin{gather}\label{P01 is ext closure} \mc P(0,1] \quad \mbox{is extension closure of} \quad \mc P_1(0,1], \mc P_2(0,1] \\
  \label{P_i subset}  \forall i \in \{1,2\} \ \ \forall t \in \RR   \qquad  \mc P_i(t) \subset   \mc P(t) \\ \label{Z(X) e1quals Z_i(X)}
 Z(X) = Z_1(X) \quad \mbox{for} \quad X \in \mc T_1;  \qquad Z(X) = Z_2(X) \quad \mbox{for} \quad X \in \mc T_2.
\end{gather}

We will show that 
\be  \label{the glued is in st_epsilon} t \in [0,\varepsilon] \ \Rightarrow \  \mc P(t) = 0.  \ee
Indeed, let $s_{11}, s_{12}, \dots, s_{1n}$  and $s_{21}, s_{22}, \dots, s_{2m}$ be the simple objects of $\mc P_1(0,1]$ and  $\mc P_2(0,1]$, respectively.  Then  $\{s_{1 i}\}_{i=1}^n \subset \sigma^{ss}_1$,  $\{s_{2 i}\}_{i=1}^m \subset \sigma^{ss}_2$ and by \eqref{help formula 4}, \eqref{Z(X) e1quals Z_i(X)}, and \eqref{phase formula} we deduce that \begin{gather} Z(s_{1 i}),  Z(s_{2 j}) \in \RR_{>0} \exp(\ri \pi (\varepsilon,1)),\end{gather}
and on the other hand by \eqref{P01 is ext closure} it follows that $Z(X)$ is a positive linear combination of  $\{ Z(s_{1 i}) \}_{i=1}^n$,  $\{ Z(s_{2 i}) \}_{i=1}^m$  for $X\in \mc P(t)\setminus \{0\}$, $t\in (0,1]$, and therefore  $Z(X) \in \RR_{>0} \exp(\ri \pi (\varepsilon,1))$ , hence \eqref{phase formula} gives $\phi_\sigma(X) \in (\varepsilon,1) $ and \eqref{the glued is in st_epsilon} follows. This in turn implies $\exp(\ri \pi [0,\varepsilon])\cap P_\sigma^{\mc T} = \emptyset$ and then for obtaining  $\sigma \in \st_\varepsilon(\mc T)$ (recall Definition \ref{definitin of st_varepsilon}) it remains to  show that $\sigma$ is a full stability condition. We will  prove  this by showing  that $\mc P(0,1]$ is a finite length  abelian category (then it follows that $\HH^{\mc P(0,1]}\cong \HH^{n+m}$ and $\sigma$ is full, since $\sigma \in \HH^{\mc P(0,1]}$). However \cite[Proposition 3.5 (a)]{CP} claims that if $0$ is an isolated point for  $\Im \left (Z_i(\mc P_i(0,1]) \right )$ for $i=1,2$ (which is satisfied due to \eqref{are of finite length} and \eqref{help formula 4}), then $\mc P(0,1)$ is a finite length  category, and on the other hand   due to \eqref{the glued is in st_epsilon} holds   $\mc P(0,1]=\mc P(0,1)$.  Therefore  indeed $\mc P(0,1]$ is finite length category   and   $\sigma$ is a full stability condition. 

Finally,  from \eqref{P_i subset}  it follows that $P^{\mc T_i}_{\sigma_i}\subset P^{\mc T}_\sigma$,  therefore  $\ol{P^{\mc T_i}_{\sigma_i} }\subset \ol{P^{\mc T}_\sigma}$, and hence  $\vol\left (\ol{P^{\mc T_i}_{\sigma_i} } \right )\leq \vol\left ( \ol{P^{\mc T}_\sigma} \right )$ for $i=1,2$, recalling \eqref{first step in main ineq}  we derive:
\begin{gather}\label{final step for main ineq}\frac{\vol\left ( \ol{P_{\sigma}}\right )}{2}\geq \max\left \{ \norm{\mc T_1}_{\varepsilon
}, \norm{\mc T_1}_{\varepsilon
} \right \}-\mu.\end{gather}
 This ineqaullity holds for any $\mu >0$     and from the very definition \ref{main def} we deduce \eqref{main inequality}.
\epr
\begin{remark} \label{stab with infinite Psigma} Let $\mc T$, $\mc T_1$, $\mc T_2$ be as in Theorem  \ref{prop inequality} (in particular there is a SOD $\mc T=\langle \mc T_1, \mc T_2\rangle$).
From the proof of Theorem \ref{prop inequality} we see that if  for some $i=1,2$ there exists a full  $\sigma \in \st(\mc T_i)$ with infinite set of phases $P_{\sigma_i}$, then there exists a full $\sigma \in \st(\mc T_i)$ with infinite  $P_{\sigma}$ as well.
\end{remark}
\begin{coro} \label{coro for esc coll} For any exceptional collection $(E_0,E_1,\dots,E_n )$ in a proper  triangulated  category and for any $0\leq i\leq n$ we have:
\begin{gather} \label{inequality for esc coll} \norm{\scal{E_0,E_1,\dots,E_n}}_{\varepsilon}\geq  \max\left \{ \norm{\scal{E_0,E_1,\dots,E_i}}_{\varepsilon}, \norm{\scal{E_{i+1},E_{i+2},\dots,E_n}}_{\varepsilon} \right \}.\end{gather} 
\end{coro}
\bpr Due to Remark \ref{from exc coll to phase gap} the categories  $\scal{E_0,E_1,\dots,E_n}$, $\scal{E_0,E_1,\dots,E_i}$, $\scal{E_{i+1},E_{i+2},\dots,E_n}$ have phase gaps. All the conditions of Theorem \ref{prop inequality} are satisfied for the semi-orthogonal  decomposition $\scal{E_0,E_1,\dots,E_n}=\scal{ \scal{E_0,E_1,\dots,E_i}, \scal{E_{i+1},E_{i+2},\dots,E_n} }$, hence equality \eqref{main inequality}  gives rise to \eqref{inequality for esc coll}.
\epr
\begin{coro} \label{blow up} Let $X$ be a smooth algebraic variety and let $Y$ be a smooth sub-variety 
so that  $K_0(D^b(X))$, $K_0(D^b(Y))$ have finite rank and $D^b(X)$, $D^b(Y)$ have phase gaps.  Denote
 by $\wt{X}$ the smooth algebraic variety obtained by blowing up $X$ along the center $Y$. 

Then $D^b(\wt{X})$ has phase gap and $\norm{D^b(\wt{X})}_{\varepsilon}\geq \max\left \{\norm{D^b(X)}_{\varepsilon}, \norm{D^b(Y)}_{\varepsilon}\right \}$ for any $\varepsilon \in (0,1)$.  \end{coro}
\bpr  \cite[Theorem 4.2]{BondalOrlov} ensures that there is a semi-orthogonal decomposition 

 $D^b(\wt{X})=\langle \mc T_1, \mc T_2,\dots,\mc T_k,  D^b(X)\rangle $, where $\mc T_i$ is equivalent to $D^b(Y)$ for $i=1,2,\dots,k$. Now Theorem \ref{prop inequality} ensures that the inequality holds.  \epr

\section{The space of stability conditions and the norms on wild Kronecker quivers} \label{norms and stab cond on Kroneckers}

This section is devoted to our main example (namely $D^b(K(l))$), where we compute both: the stability space and $\norm{\cdot}_{\varepsilon}$.

\subsection{Recollection on the action of \texorpdfstring{${\rm SL(2,\ZZ)}$}{\space} on \texorpdfstring{$\mc H$}{\space}} \label{recollection}

 Let ${\rm GL}^+(2,\RR)$ be the group of $2\times 2$ matrices with positive determinant. Recall that (see e.g. \cite{BZ}) for any matrix $\gamma \in{\rm GL}^+(2,\RR) $ we have a biholomorphism: 
\begin{gather} \label{action} {\rm GL}^+(2,\RR) \ni \begin{bmatrix} a & b \\ c & d \end{bmatrix}= \alpha : {\mc H} \rightarrow {\mc H} \qquad \alpha(z)=\frac{a z + b}{c z + d},  \end{gather}
 and this defines an action of ${\rm GL}^+(2,\RR)$ on $\mc H$.  Let $ {\rm SL}(2,\ZZ)\subset {\rm GL}^+(2,\RR)$ be the  subgroup of matrices with integer coefficients and determinant $1$.  The action of ${\rm SL}(2,\ZZ)$  on  $\mc H$  defined by the formula \eqref{action}   is properly discontinuous   (see e.g. \cite[p.17 and p.20]{Miyake}).\footnote{The definition of properly discontinuous can be seen in the proof of Corollary \ref{properly discontinuous}.}

\subsubsection{Hyperbolic, parabolic, elliptic elements} \label{hpe} 
An element $\alpha \in {\rm GL}^+(2,\RR)$ is called \textit{ elliptic, parabolic, hyperbolic}, if $\tr(\alpha)^2<4 \det(\alpha)$, $\tr(\alpha)^2=4 \det(\alpha)$, $\tr(\alpha)^2>4 \det(\alpha)$, respectively. When $\alpha$ is parabolic or hyperbolic, then $\alpha$ has no fixed points in $\mc H$ (see \cite[p. 7]{Miyake}).  Furthermore,  if $\alpha $ is elliptic, parabolic, or hyperbolic and $\alpha \in {\rm SL}(2,\ZZ)$, then each non-trivial element in the subgroup   $\langle \alpha \rangle$  generated by $\alpha$ is elliptic, parabolic, or hyperbolic, respectively. 
 
\subsubsection{Fundamental domain of a hyperbolic element}  \label{0FD remark hyperbolic} 
The definition of a fundamental domain in $\mc H$ of a subgroup $\Gamma \subset {\rm SL}(2,\ZZ)$ which we adopt here is  in \cite[p. 20]{Miyake}. We will need to determine fundamental domains of subgroups of the form $\langle \alpha \rangle$ for a non-scalar $\alpha \in {\rm SL}(2,\ZZ)$. The following  arguments will be useful.

 Let  $\alpha \in {\rm SL}(2,\ZZ)$ be a  hyperbolic element. It is well known that $\alpha$ is conjugate (in $\rm GL^+(2,\RR)$) to a  matrix  of the form \eqref{beta 1}, i.e. $\alpha'=\beta^{-1} \alpha \beta$ for some $\beta$ (see e.g. \cite[Lemma 1.3.2]{Miyake}): \be \label{beta 1} \alpha'= \begin{bmatrix}
	a & 0 \\
	0 & d
\end{bmatrix} \quad a d >0 \qquad \alpha'(z)= (a/d) z.\ee  It is clear that  any strip $F'$ as in Figure \eqref{0FDhyperbolic}, where $\delta>0$, is a  fundamental domain for the subgroup $\langle \alpha' \rangle $. Two points in $F'$ lie  in common  orbit iff they are of the form  $z_{-}\in\mk{b}_-'$, $z_+=(a/d) z_-$. Transforming $F'$ via  $\beta$  results in  a  fundamental domain $F$ for $\langle \alpha \rangle$, i. e. $F=\beta(F')$.
 
\begin{remark} \label{circles and lines} For computing $\beta(F')$ we will use the feature of  $\beta$ that it maps circles or lines   perpendicular to the real axis to circles or lines perpendicular to the real axis (see \cite[Lemma 1.1.1]{Miyake} and recall that $\beta$ is a conformal map, which maps the real axis into itself).
 \end{remark}

\subsubsection{Fundamental domain of a parabolic element}  \label{0FD remark parabolic} 
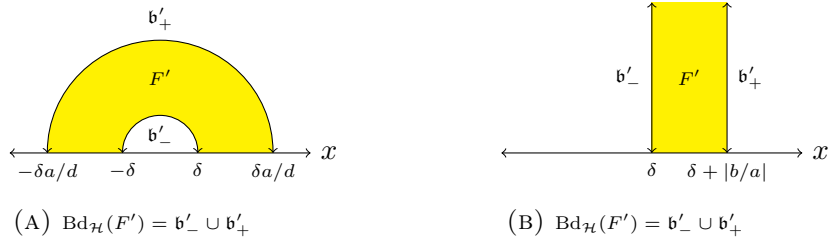
\begin{figure} \centering
   \hspace{-20mm}   \begin{subfigure}[b]{0.2\textwidth}
    \begin{tikzpicture} 
\path [fill=yellow] (-1.5,0) --(-0.5,0) to [out=90,in=180] (0,0.5) to [out=0,in=90] (0.5,0)-- (1.5,0) to [out=90,in=0] (0,1.5)to [out=180,in=90] (-1.5,0) ;	
\draw [<->] (-2,0)--(2,0);
\draw  [<->]    (0.5,0) arc [radius=0.5, start angle=0, end angle= 180];
\node [below] at (0.5,0) {\tiny $\delta$\normalsize};
\node [below] at (-0.5,0) {\tiny $-\delta$\normalsize};
\node [below] at (0,0.5) {\tiny $\mk{b}_-'$\normalsize};
\node [above] at (0,1.5) {\tiny $\mk{b}_+'$\normalsize};
\draw  [<->]    (1.5,0) arc [radius=1.5, start angle=0, end angle= 180];
\node [below] at (1.5,0) {\tiny $\delta a/d$\normalsize};
\node   at (0,1) {\tiny $F'$\normalsize};
\node [below] at (-1.5,0) {\tiny $- \delta a/d$\normalsize};
\node [right] at (2,0) {$x$};
\end{tikzpicture}
        \caption{\tiny $ {\rm Bd}_{\mc H}(F')=\mk{b}_-'\cup \mk{b}_+'$\normalsize}
        \label{0FDhyperbolic}
    \end{subfigure}
		\hspace{30mm}
		 \begin{subfigure}[b]{0.2\textwidth}
       \begin{tikzpicture} 
			\path [fill=yellow](0,0)--(1,0)-- (1,2) -- (0,2)-- (0,0) ;	
\draw [<->] (-2,0)--(2,0);
\draw [<->] (0,0)--(0,2);
\node [right] at (2,0) {$x$};
\node [below] at (0,0) {\tiny $\delta$\normalsize};
\draw [<->] (1,0)--(1,2);
\node [below] at (1,0) {\tiny $\delta+\abs{b/a}$\normalsize};
\node [left] at (0,1) {\tiny $\mk{b}_-'$\normalsize};
\node [right] at (1,1) {\tiny $\mk{b}_+'$\normalsize};
\node   at (0.5,1) {\tiny $F'$\normalsize};
\end{tikzpicture}
        \caption{\tiny $ {\rm Bd}_{\mc H}(F')=\mk{b}_-'\cup \mk{b}_+'$\normalsize}
        \label{0FDparabolic}
    \end{subfigure} 
     \caption{Fundamental domains}\label{0FD}
\end{figure}

If $\alpha \in {\rm SL}(2,\ZZ)$ is   parabolic,  then  $\alpha$ is conjugate    to a  matrix  of the form \eqref{beta 2}, i.e. $\alpha'=\beta^{-1} \alpha \beta$ for some $\beta  \in \rm GL^+(2,\RR)$ (\cite[Lemma 1.3.2]{Miyake}): \be \label{beta 2} \alpha'= \begin{bmatrix}
	a & b \\
	0 & a
\end{bmatrix} \quad a\neq 0 \qquad \alpha'(z)=z+b/a.\ee  
 Any strip $F'$ as in Figure \eqref{0FDparabolic}, where $\delta\in \RR$, is a  fundamental domain for subgroup $\langle \alpha' \rangle $. Two points in $F'$ lie in a common  orbit iff they are of the form  $z_{-}\in\mk{b}_-'$, $z_+=z_-+\abs{b/a}$. Transforming $F'$ via $\beta$ results in a  fundamental domain $F$ for $\langle \alpha \rangle$, i. e. $F=\beta(F')$. 

\begin{remark} \label{prop disc free} Let  $\alpha\in {\rm SL}(2,\ZZ)$ be a non-scalar matrix. From the arguments in Subsections \ref{0FD remark parabolic}, \ref{0FD remark hyperbolic} it follows that if $\alpha$ is either  hyperbolic or parabolic, then we have:

(a) $\langle \alpha \rangle\cong \ZZ$ acts free and properly discontinuous on $\mc H$;

(b) Let $F$ be a fundamental domain of $\langle \alpha \rangle$ obtained by the method in Subsections \ref{0FD remark hyperbolic}, \ref{0FD remark parabolic}. 

\hspace{4mm} Then  for each $u\in {\rm Bd}_{\mc H}(F)$ there exists  an open subset $U\subset \mc H$, s. t. $u\in U$ and

\hspace{4mm} $\{i\in \ZZ: \alpha^i(U)\cap F\neq \emptyset\}$ is finite (in fact contains only two elements).
  \end{remark}

\subsection{The Helix in  \texorpdfstring{$D^b(Rep_k(K(l)))$}{\space} for \texorpdfstring{$l\geq 2$}{\space} }  \label{there are no Ext-nontrivial...}
From now on we assume that $l\geq 2$ and denote $\mc T_l = D^b(Rep_k(K(l)))$. On \cite[p. 668 down]{Macri}  Macr\`i constructs  a family of exceptional objects   $\{s_i\}_{i\in \ZZ}$ and states that: if $i<j$, then  $\hom^k(s_i,s_j)\neq 0$ only if $k=0$, and  $\hom^k(s_j,s_i)\neq 0$ only if $k=1$, without giving a proof of this statement. In this section  we view  $\{s_i\}_{i\in \ZZ}$ as a helix, as defined in \cite[p. 222]{BP} (and introduced earlier in \cite{B}, \cite{GoroRuda}) and using the results about geometric helices in  \cite[p. 222]{BP}  we give a simple proof of  Lemma \ref{nonvanishings}.   

Lemma \ref{actionLEMMA} is the place where the ${\rm SL}(2,\ZZ)$-action on $\mc H$  comes into play, it is a simple but important observation  for the rest of the paper.

 We write $\ul{\dim}(X)=(n,m)$, $\ul{\dim}_0(X)=n$,  $\ul{\dim}_1(X)=m$. 
 for a representation of the form: \vspace{2mm} $$X= \begin{diagram}[w=3em]
k^n & \upperarrow{}
\lift{-2}{  \vdots }
\lowerarrow{} & k^m
\end{diagram} \in Rep_k(K(l))$$  

Recall that $Rep_k(K(l))$ is hereditary category in which for any two $X, Y \in Rep_k(Q)$ with dimension vectors $\ul{\dim}(X)=(n_x,m_x)$,  $\ul{\dim}(Y)=(n_y,m_y)$ holds the equality (the Euler Formula): 
\be \label{euler} \hom(X,Y)-\hom^1(X,Y) = n_x n_y + m_x m_y -l n_x m_y \ee

Let $s_0, s_1 \in \mc T_l$ be so that $s_0[1]$, $s_1$ are  the simple objects in $Rep_k(Q)$ with $\ul{\dim}(s_0[1])=(1,0)$ and  $\ul{\dim}(s_1)=(0,1)$. Using \eqref{euler} one easily computes $\hom(s_0,s_1)=l$, $\hom^p(s_0,s_1)=0$ for $p\neq 0$ and $\hom^*(s_1,s_0)=0$. With  the terminology from Section \ref{notations} we can say that $(s_0, s_1)$ is a strong  exceptional pair, furthermore  it is a full exceptional pair in $\mc T_l = D^b(Rep_k(K(l)))$.  

\begin{remark} \label{mutations}  Recall that (see e.g. \cite[p. 222]{BP}) for any exceptional pair $(A,B)$ in any proper triangulated category $\mc T$ one defines objects  $L_A(B)$ (left mutation) and $R_B(A)$ (right mutation) by the
 triangles \be \label{L_A and R_B} \bd L_A(B)&\rTo &\Hom^*(A,B)\otimes A & \rTo^{ev^*_{A,B}} & B \ed  \qquad   \bd A &\rTo^{coev^*_{A,B}} &\Hom^*(A,B)^{\check{}}\otimes B &\rTo & R_B(A) \ed \ee
 and   $(L_A(B),A)$, $(B,R_B(A))$ are exceptional pairs as well,  they are full  if $(A,B)$ is full. 

It follows, that for  any exceptional collection $\mc E=(E_0,E_1,\dots, E_n)$ in $\mc T$ and 
for  any $0\leq i < n$  the sequences  $R_i(\mc E)=(E_0,E_1,\dots,E_{i+1},R_{E_{i+1}}(E_i),\dots,  E_n)$,  $L_i(\mc E)=(E_0,E_1,\dots,L_{E_{i}}(E_{i+1}),E_{i},\dots,  E_n)$ are  exceptional and $\scal{R_i(\mc E)}=\scal{L_i(\mc E)} = \scal{\mc E}$.  The sequences $L_i(\mc E)$ and $R_i(\mc E)$ are called left and right mutations of $\mc E$.
\end{remark}
In particular, from the exceptional pair $(s_0,s_1)$ we get objects $L_{s_0}(s_1)$, $R_{s_1}(s_0)$,  denoted by $s_{-1}$, $s_2$, respectively, and each two adjacent elements in the sequence  $s_{-1}, s_{0} ,s_1,s_2$ form a full exceptional pair.  Applying iteratively  left/right mutations on the left/right standing exceptional pair  generates a sequence (infinite in both directions) of exceptional objects $\{s_i\}_{i\in \ZZ}$. This is the helix induced by the exceptional pair $(s_0,s_1)$, as defined in \cite[p. 222]{BP}. Any two adjacent elements in this sequence form a full exceptional pair in $\mc T$. Actually the right mutation generates an action of $\ZZ$ on the set of equivalence classes of exceptional pairs (w.r. $\sim$), (the inverse is the left mutation).    By the transitivity of this action   shown in   \cite{WCB1}, it follows that in $ \mc T_l $:
\begin{gather} \label{list of exceptional pairs} \textit{complete lists  of  exceptional  pairs and objects (up to shifts) are} \  \{(s_i,s_{i+1})\}_{i\in\ZZ} \ \textit{and} \ \{s_i\}_{i\in\ZZ}. \end{gather}
It  is shown in \cite[Example 2.7]{BP} that any strong exceptional pair  is geometric (\cite[Definition on p. 223]{BP}), which applied to our strong exceptional collection   $(s_0,s_1)$ gives the following vanishings:
\be \label{nonvanish1} i\leq j, p\neq 0 \qquad \Rightarrow \qquad  \hom^p(s_i,s_j)=0. \ee  We will show in Lemma \ref{nonvanishings}  that, when $p=0$, the dimensions above  do not vanish.  When $j=i+1$, we use the equalities $\hom^p(L_A(B),A)=\hom^{-p}(A,B)$ for any $p\in \ZZ$, and analogously for $R_B(A)$, (see again \cite[Example 2.7]{BP}, or  \cite[p. 157 down]{DHKK} for details)  to deduce that   $\hom(s_i,s_{i+1})=\hom(s_0,s_1)=l$ for any $i$.   We obtained $\{s_i\}_{i\in \ZZ}$ using  the triangles \eqref{L_A and R_B}, therefore  for any $i\in \ZZ$ exists  a triangle in $\mc T$:
\be \label{triangles}  s_{i-1}\rightarrow  s_{i}^l \rightarrow s_{i+1} \rightarrow  s_{i-1}[1]. \ee Using these triangles  we prove:

  \begin{lemma} \label{actionLEMMA} Let $i\in \ZZ$ be any integer.  Let $Z:K_0(\mc T_l) \rightarrow \CC$ be any group homomorphism, s.t.  
	$\frac{Z(s_{i})}{Z(s_{i-1})} \in \mc H$. Then for any $j \in \ZZ$ we have  $\frac{Z(s_{i+j})}{Z(s_{i-1+j})}\in \mc H$ and  the equality: 
\begin{gather} 
\label{actionCENTRAL}    \frac{Z(s_{i+j})}{Z(s_{i-1+j})}=\alpha_l^j\left (\frac{Z(s_{i})}{Z(s_{i-1})}\right ), 
  \end{gather}
where    $ \alpha_l=\begin{bmatrix}
	l & -1 \\
	1 & 0
\end{bmatrix} \in {\rm SL}(2,\ZZ) $  and  $\alpha_l^j: \mc H \rightarrow \mc H$  is the corresponding  automorphism  given in  \eqref{action}.
	\end{lemma}	
		\bpr		From 	\eqref{triangles} it follows that $Z(s_{i+1})=l Z(s_i) -  Z(s_{i-1})$ for any $i\in \ZZ$, hence we can write:
\begin{gather} 
\begin{bmatrix}
	Z(s_{i+1}) \\
	Z(s_{i})
\end{bmatrix} = \alpha_l \begin{bmatrix}
	Z(s_{i}) \\
	Z(s_{i-1})
\end{bmatrix}; \quad  \begin{bmatrix}
	Z(s_{i}) \\
	Z(s_{i-1})
\end{bmatrix} = \alpha_l^{-1} \begin{bmatrix}
	Z(s_{i+1}) \\
	Z(s_{i})
\end{bmatrix}, \ \  \textit{where} \ \  \alpha_l=\begin{bmatrix}
	l & -1 \\
	1 & 0
\end{bmatrix} \in {\rm SL}(2,\ZZ)
  \end{gather}
	hence by induction  we obtain:
	\begin{gather} \label{action_help} \begin{bmatrix}
	Z(s_{i+j}) \\
	Z(s_{i-1+j})
\end{bmatrix} = \begin{bmatrix}
	l & -1 \\
	1 & 0
\end{bmatrix}^j \begin{bmatrix}
	Z(s_{i}) \\
	Z(s_{i-1})
\end{bmatrix}  \qquad  \qquad j\in \ZZ, \end{gather} therefore $ \frac{Z(s_{i+j})}{Z(s_{i-1+j})}=\frac{(\alpha_l^j)_{11} Z(s_{i}) + (\alpha_l^j)_{12} Z(s_{i-1})}{(\alpha_l^j)_{21} Z(s_{i}) + (\alpha_l^j)_{22} Z(s_{i-1})}$ and,  recalling \eqref{action}, we deduce \eqref{actionCENTRAL}. 
\epr
 The statement of the following lemma is a slight modification of \cite[Lemma 4.1]{Macri} and we give a proof here:
\begin{lemma} \label{nonvanishings} Assume that $l\geq 2$. Then no two elements in $\{s_i\}_{i\in \ZZ}$ are isomorphic and:  $s_{\leq 0}[1], s_{\geq 1} \in Rep_k(K(l))$. Furthermore, we have: 
\begin{gather} \label{nonvanish11} i\leq j \quad \Rightarrow \quad   \left \{ \begin{array}{c c} \hom(s_i,s_j)  \neq 0 &  \\  \hom^p(s_i,s_j) = 0 & \ \ \text{for} \ \  p\neq 0  \end{array} \right . ;\\ 
\label{nonvanish12} i > j+1 \quad \Rightarrow \quad   \left \{ \begin{array}{c c} \hom^1(s_i,s_j)  \neq 0 &    \\ \hom^p(s_i,s_j) = 0 & \ \ \text{for} \ \  p\neq 1  \end{array} \right . .\end{gather}    \end{lemma} 
\bpr	The matrix $\alpha_l$ in \eqref{actionCENTRAL}  has trace $\tr(\alpha_l)= l\geq 2$ and therefore  (see Subsection \ref{hpe}) $\alpha_l^j$ is either parabolic or hyperbolic for any $j\in \ZZ$, in particular it has no fixed points in $\mc H$.  If $s_i \cong s_{i+j}$ for some $i\in \ZZ$, $j\in \NN$, then since $(s_{i-1}, s_i)$ and $(s_{i+j-1}, s_{i+j})$  are both full exceptional pairs it follows that  $s_{i-1} \cong s_{i+j-1}$ (note that $s_{i-1} \cong s_{i+j-1}[k]$ for $k\neq 0$ is impossible by the already proved  \eqref{nonvanish1}), and hence $\frac{Z(s_{i+j})}{Z(s_{i-1+j})}=\frac{Z(s_{i})}{Z(s_{i-1})}$, which contradicts \eqref{actionCENTRAL}. Therefore no two elements in $\{s_i\}_{i\in \ZZ}$ are isomorphic, 
 and  the non-vanishings in  \eqref{nonvanish11} follow. Indeed, if $\hom(s_i,s_j)=0$ for some $i<j$, then $(s_j,s_i)$ is a full exceptional pair, which contradicts  the already proven $s_{i-1} \not \cong s_j$.   
  Since $Rep_k(K(l))$ is hereditary category, the non-vanishings  in  \eqref{nonvanish11}  imply that $s_{\leq 0}[1], s_{\geq 1} \in Rep_k(K(l))$.  

 Recall \cite{BK} that there exists an exact  functor (the Serre functor)  $F:\mc T_l \rightarrow \mc T_l$, s. t. $\hom(A,B)=\hom(B,F(A))$ for any two objects $A,B \in \mc T$.   Furthermore,  the formula on \cite[p.223 above]{BP} (in our case $n=2$) says that $F(s_i)\cong s_{i-2}[1]$  for each $i\in \ZZ$ , hence for any three integers $i,j,p$  we have $\hom^p(s_i,s_j)=\hom^{1-p}(s_j,s_{i-2})$. Now \eqref{nonvanish12} follows from  \eqref{nonvanish11}.  
\epr The following corollary of \cite[Theorem 0.1]{Miyachi} ensures   existence of a  functor in
 $\Aut( \mc T_l)$, which plays the role of the functor $(\cdot)\otimes \mc O(1)$   for any $l\in \ZZ$. \begin{coro} {\rm (of \cite[Theorem 0.1]{Miyachi})} \label{the auto-eq F} There exists  $A_l\in \Aut(\mc T_l)$, s. t. $A_l(s_i)\cong s_{i+1}$ for each $i\in \ZZ$ (recall that $l\geq 2$).  \end{coro} \bpr In \cite[Definitions 2.2, 2.3]{Miyachi},   from any quiver  $Q$   they   define a quiver, denoted by ${\bf \Gamma}^{irr}$,  which is the disjoint union of all connected components of the  Auslander-Reiten quiver of $D^b(Q)$,   isomorphic to the connected component which contains  the simple representations of $Q$.    \cite[Theorem 2.4 ]{Miyachi} applied to the quiver $Q=K(l)$  gives an isomorphism of quivers $\rho: \ZZ\times (\ZZ K(l))\rightarrow {\bf \Gamma}^{irr}$, where:  

 (a) $\ZZ K(l)$ is  isomorphic to quiver with set of vertices $\ZZ$ and for any $(i,j) \in \ZZ \times \ZZ$ there are $l$ arrows from $i$ to $j$ iff $j=i+1$ and no arrows otherwise; 

   (b) $\ZZ \times (\ZZ K(l)) $ is quiver whose connected components are $\{i\}\times \ZZ K(l) $ and each of them is a labeled copy of $\ZZ K(l)$. 
	
	(c) The translation functor of $\mc T_l$ induces action on ${\bf \Gamma}^{irr}$ which by $\rho$ corresponds to an automorphism $\sigma$ of $\ZZ \times (\ZZ K(l))$ acting by increasing the first component by $1$.  The Serre functor shifted by $[-1]$ induces an automorphism  $\tau $ of $\ZZ \times (\ZZ K(l))$, which on vertices maps $(i,j)$ to $(i,j-2)$. 
		
		Taking into account \eqref{nonvanish11}, \eqref{nonvanish12} and that the Serre functor maps $s_i$ to $s_{i-2}[+1]$ (see the proof of Lemma \ref{nonvanishings}) we see that $\rho$ can be chosen so that $\rho(i,j)\cong s_j[i]$ for any $i,j\in \ZZ$. Combining \cite[Corollary 1.9; Theorem 3.7]{Miyachi} and the last paragraph of the proof of \cite[Theorem 4.3]{Miyachi} we see that  $\Aut(\mc T_l)\cong \ZZ \times (\ZZ \ltimes {\rm PGL}_n(k))$ and a generator of the factor $\ZZ$ in $\ZZ \ltimes {\rm PGL}_n(k)$ is the desired  $A_l\in \Aut(\mc T_l)$.
\epr
\subsection{The principal \texorpdfstring{$\CC$}{\space}-bundle \texorpdfstring{$  \st(\mc T_l) \rightarrow \mc X_l $}{\space}} \label{the union}
We will denote by $\mc X_l$ the set of orbits of the $\CC$-action on $\st(\mc T_l)$, i.e. $\mc X_l=\st(\mc T_l)/\CC$.   In this section we show that $\mc X_l$ is a complex manifold biholomorphic to  either $\CC$ or $\mc H$  (Corollary \ref{trvial bundle}) and that the action of $\langle A_l \rangle $ on $\st(\mc T_l)$, where $A_l\in \Aut(\mc T_l)$ is the functor from Corollary \ref{the auto-eq F},   descends to an action by biholomorphism on $\mc X_l$   (Corollary \ref{open cover of mc X}).

 In \cite[Section 3.3]{Macri} are constructed stability conditions generated by a full Ext-exceptional collection. The set of stability conditions generated by a full Ext-exceptional pair $(A,B)$ will be denoted by $\Theta'_{(A,B)}$ and for any full exceptional pair\footnote{not necessarily Ext} $(A,B)$ the notation  $\Theta_{(A,B)}$  will denote the union of the sets $\Theta'_{(A[a],B[b])}$ s.t. $(A[a],B[b])$ is an Ext-pair (see  \cite[formulas (8), (10)]{DK3}).  The  idea of  Macr\`i  in \cite[Subsection 3.3 and Section 4]{Macri}  how to prove    that $\st(\mc T_l)$ is simply connected  is to show that  $\st(\mc T_l)$ is covered by  $\{\Theta_{(s_i,s_{i+1})} \}_{i\in \ZZ}$  and then use Seifert-Van Kampen theorem. We prove here  Proposition \ref{lemma for f_E(Theta_E)}, formula \eqref{st}, and Lemma \ref{sections} following this insight of Macr\`i. 
 
 The proof of the following proposition is a simpler analogue of the proof of \cite[Proposition 2.7.]{DK3}: 
\begin{prop} \label{lemma for f_E(Theta_E)} 
For $i\in \ZZ$ the subset $ \Theta_{(s_i,s_{i+1})}\subset \st(\mc T_l)$ has the following description:  \begin{gather} \label{Theta_E n=2} \Theta_{(s_i,s_{i+1})} =\left  \{\sigma\in \st(\mc T): (s_i,s_{i+1}) \subset \sigma^{ss}  \ \mbox{and} \  \begin{array}{l} \phi_\sigma(s_i) < \phi_\sigma(s_{i+1})\end{array} \right \}.  \end{gather}  
In particular, the set $\Theta_{(s_i,s_{i+1})}$ is biholomorphic to the contractible set $\mc S =\{(z_1,z_2)\in \CC^2 ; \Im(z_1)<\Im(z_2)\}\subset \CC^2$ by the  following map $\bd \Theta_{(s_i,s_{i+1})}& \rTo^{\varphi_i} &  \mc S \ed $  \begin{gather} \label{mapping}
\Theta_{(s_i,s_{i+1})}\ni (Z, \mc P) \mapsto \left ( \log\abs{Z(s_i)}+ \ri \pi \phi_{\sigma}(s_i), \log\abs{Z(s_{i+1})}+ \ri \pi \phi_{\sigma}(s_{i+1}) \right )\in \mc S.\end{gather}   
\end{prop}
\bpr We use \cite[formulas (12), (17), (18), and Lemma 2.4]{DK3} and deduce that $\sigma \in \Theta_{(s_i,s_{i+1})}$ iff  $(s_i,s_{i+1}) \subset \sigma^{ss} $ and 
$(\phi_\sigma(s_i),\phi_\sigma(s_{i+1}) )\in \bigcup_{ \textbf{p} \in A_0 }    S^1(-\infty,1)  -  \textbf{p}$, where $ S^1(-\infty,1)=\{(x,y): x-y<1\} $ and $ A_0=\left \{(0,p) \in \ZZ^{2} :  (s_i,s_{i+1}[p]) \ \mbox{is Ext} \right \}$. From \eqref{nonvanish11} we see that $ A_0=\left \{(0,p) : p\leq -1 \right \}$ and \eqref{Theta_E n=2} follows.  Since the  map defined by    \cite[formulas (9), (11)]{DK3} is homeomorphism, it follows that   \eqref{mapping} is homeomorphism. Let us identify $\Hom(K_0(\mc T_l),\CC) \cong \CC^2$ via the basis $[s_i], [s_{i+1}]$ of $K_0(\mc T_l)$, and let $\bd \st(\mc T_l)&\rTo^{proj}&\Hom(K_0(\mc T_l),\CC)\ed$ be the projection  $proj(Z, \mc P)=Z$. Then the following diagram: \bd[size=1.5em] \Theta_{(s_i,s_{i+1})} & \rTo^{proj}                & \Hom(K_0(\mc T_l),\CC) \\
                    \dTo^{\varphi_i}               &                           &    \dTo^{\simeq}                         \\
												\mc S          &   \rTo^{\exp \times \exp} &  \CC^2  \ed   is commutative. Since the horizontal arrows are local biholomorphisms and we already showed that $\varphi_i$ is homeomorphism, it follows that   $\varphi_i$ is biholomorphic.  \epr

\cite[Lemma A.1]{DK1}  says   that for each $\sigma \in \st(\mc T_l)$ there exists a $ \sigma$-exceptional pair.  This means that (see \cite[Corollary 3.18]{DK1})  for each $\sigma$ there exists an Ext-exceptional pair $(A,B)$ with  $\sigma \in \Theta_{(A,B)}'$. From \eqref{list of exceptional pairs}  we see that  $(A,B)$ is  of the form $(s_i[a],s_{i+1}[b])$ for  some  $i,a,b \in \ZZ$. Therefore:
 \be \label{st} \st(\mc T_l) = \bigcup_{i\in \ZZ} \Theta_{(s_i,s_{i+1})}. \ee
\begin{lemma} \label{sections} Let $i,j$ be two different integers. Then the equality below holds
 \begin{gather} \label{intersection of two} \Theta_{(s_i,s_{i+1})} \cap \Theta_{(s_j,s_{j+1})} = \left  \{\sigma\in \st(\mc T_l): (s_i,s_{i+1}) \subset \sigma^{ss}  \ \mbox{and} \   \phi_\sigma(s_i) < \phi_\sigma(s_{i+1})<\phi_\sigma(s_i)+1 \right \}, \end{gather}
 and therefore $ \Theta_{(s_i,s_{i+1})} \cap \Theta_{(s_j,s_{j+1})} = \bigcap_{p\in \ZZ}  \Theta_{(s_p,s_{p+1})}$. It follows that $\st(\mc T_l)$  is contractible. \end{lemma}
\bpr We show first the inclusion $\supset$. So let   $(s_i,s_{i+1}) \subset \sigma^{ss}  \ \mbox{and} \   \phi_\sigma(s_i) < \phi_\sigma(s_{i+1})<\phi_\sigma(s_i)+1$. We will show  below that the given inequalities imply: 

\begin{gather} \label{intermediate step} \sigma \in \Theta_{(s_{i-1},s_{i})}\cap \Theta_{(s_{i+1},s_{i+2})};\begin{array}{c} \phi_\sigma(s_{i-1}) < \phi_\sigma(s_{i})<\phi_\sigma(s_{i-1})+1\\  \phi_\sigma(s_{i+1}) < \phi_\sigma(s_{i+2})<\phi_\sigma(s_{i+1})+1\end{array}. \end{gather}  Then by induction we obtain the inclusion $\subset $. We  use \cite[Proposition 2.2]{DK3} first  to show that  $s_{i-1}$ and  $s_{i+2}$ are semi-stable. More precisely, the given inequality  is the same as $\phi_\sigma(s_{i+1}[-1])<\phi_\sigma(s_i)<\phi_\sigma(s_{i+1})$, which together with \eqref{nonvanish11} imply that $s_{i},s_{i+1}[-1]$ is a $\sigma$-exceptional pair (see \cite[Definition 3.17]{DK1}).  From two consecutive  triangles  in the sequence of triangles \eqref{triangles} it follows that $s_{i-1}$ and  $s_{i+2}[-1]$ are in the extension closure of  $s_{i},s_{i+1}[-1]$. Since we have also $\hom(s_i,s_{i+1})\neq 0$, we can apply   \cite[Proposition 2.2]{DK3}  and it ensures that  $s_{i-1},s_{i+2}[-1]\in \sigma^{ss}$.  
 Now from \eqref{nonvanish11}, \eqref{nonvanish12} it follows:
\begin{gather} \label{one intermediate ineq} \phi(s_{i-1})\leq \phi(s_{i})\leq \phi(s_{i+1}) \leq \phi(s_{i+2})\leq \phi(s_{i-1})+1 \leq \phi(s_{i+1})+1.   \end{gather}
The given $\phi_\sigma(s_i) < \phi_\sigma(s_{i+1})<\phi_\sigma(s_i)+1$ amounts to $\frac{Z(s_{i+1})}{Z(s_{i})}\in \mc H$ and then  Lemma \ref{actionLEMMA} ensures that   $\frac{Z(s_{i})}{Z(s_{i-1})}, \frac{Z(s_{i+2})}{Z(s_{i+1})}\in \mc H$, which together with \eqref{one intermediate ineq} leads to  $\phi(s_{i-1})< \phi(s_{i})< \phi(s_{i-1})+1$, $  \phi(s_{i+1})< \phi(s_{i+2}) < \phi(s_{i+1})+1$. Thus we derived \eqref{intermediate step} and the inclusion  $\supset$ follows. 

To show the opposite inclusion  $\subset$  it is enough to consider the case $i<j$. So if $\sigma \in \Theta_{(s_i,s_{i+1})} \cap \Theta_{(s_j,s_{j+1})}$, then \eqref{Theta_E n=2} shows that $s_i,s_{i+1},s_j,s_{j+1} \in \sigma^{ss}$ and $\phi(s_i)<\phi(s_{i+1})$,  $\phi(s_j)<\phi(s_{j+1})$, therefore using \eqref{nonvanish11}, \eqref{nonvanish12} we get:
\begin{gather} \phi(s_i)<\phi(s_{i+1})\leq \phi(s_j) < \phi(s_{j+1})\leq\phi(s_i)+1 \end{gather} 
and the equality \eqref{intersection of two} is proved.  The map $\varphi_i$ in \eqref{mapping} restricts to a biholomorphism  \bd \Theta_{(s_i,s_{i+1})}\cap  \Theta_{(s_j,s_{j+1})} & \rTo^{{\varphi_i}_{\vert}} & \{(z_1,z_2)\in \CC^2 ; \Im(z_1)<\Im(z_2)<\Im(z_1)+1\}\subset \CC^2\ed hence $ \Theta_{(s_i,s_{i+1})} \cap \Theta_{(s_j,s_{j+1})} = \bigcap_{p\in \ZZ}  \Theta_{(s_p,s_{p+1})}$ is contractible. 
Since $\st(\mc T_l)$ is covered by the contractible sets $\{ \Theta_{(s_i,s_{i+1})} \}_{i\in \ZZ}$ (see \eqref{st}), using \cite[Remark A.6]{DK3} we deduce that $\st(\mc T_l)$ is contractible.  
   \epr 
	\begin{lemma} \label{sections1} Let  $\sigma =(Z,\mc P) \in \Theta_{(s_i,s_{i+1})}$ (in particular we have $s_i,s_{i+1} \in \sigma^{ss}$ and $\phi(s_i)<\phi(s_{i+1})$). If $\phi(s_{i+1})>\phi(s_i)+1$, then $s_j\not \in \sigma^{ss}$ for $j\neq i, j \neq i+1$.    If $\phi(s_{i+1})=\phi(s_i)+1$, then $s_j\in \sigma^{ss}$ for each $j$.
	\end{lemma}
	\bpr Assume first that  $\phi(s_{i+1})>\phi(s_i)+1$. If $s_j \in \sigma^{ss}$ for some $j>i+1$, then due to \eqref{nonvanish11} $\phi(s_j)\geq \phi(s_{i+1})$, and hence $\phi(s_j)>\phi(s_i)+1$, $\hom^1(s_j,s_i)=0$, which contradicts \eqref{nonvanish12}.  The arguments for the case $j<i$ are similar.
	
	Assume now  that  $t=\phi(s_{i+1})=\phi(s_i)+1$. Recall that $s_0[1]$ and $s_1$ are the simple ojects in $Rep_k(K(l))$ (see after \eqref{euler}). It follows that for each $j$ we  $s_j[k]$ is in the extension closure of  $s_0[1]$ and $s_1$ for some $k$. Due to Corollary \ref{the auto-eq F} it follows that  for each $j$ we  $s_j[k]$ is in the extension closure of  $s_i[1]$ and $s_{i+1}$ for some $k$, and since  $s_i[1], s_{i+1} \in \mc P(t)$ it follows that $s_j[k]\in \mc P(t)$, therefore $s_j\in \sigma^{ss}$.
	\epr
	\begin{df} \label{def of cal Z} We  will denote the quotient  $ \st(\mc T_l)/\CC$ by $\mc X_l$,   the corresponding projection by $\bd  \st(\mc T_l)& \rTo^{\pr} & \mc X_l\ed $. The intersection  $ \bigcap_{p\in \ZZ}  \Theta_{(s_p,s_{p+1})}$ will be denoted by  $\mc Z$.  Due to Lemma \ref{sections} we have $\mc Z= \Theta_{(s_i,s_{i+1})} \cap \Theta_{(s_j,s_{j+1})} $ for any $i\neq j$.
	
	From \eqref{st}  we get a disjoint union $\st(\mc T_l) = \mc Z \amalg \amalg_{i\in \ZZ} \left ( \Theta_{(s_i,s_{i+1})}\setminus \mc Z \right )$. \end{df}
	 Corollary \ref{trvial bundle 1}  and Lemma \ref{sections} imply:
	\begin{coro} \label{trvial bundle} $\mc X_l$ is  biholomorphic either to $\CC$ or to $\mc H$ and $\bd  \st(\mc T_l)& \rTo^{\pr} & \mc X_l\ed $ is trivial $\CC$-principal bundle.
\end{coro}
  The action \eqref{left action} descents to an action  by biholomorphisms on $\mc X_l$. To show this and some basic properties of this action we note first: 
	\begin{lemma} \label{actions on the blocks}  For any $i, j\in \ZZ$, any $\lambda \in \CC$, and  $A_l\in \Aut(\mc T_l)$  from Corollary \ref{the auto-eq F}  hold the equalities:
	\begin{gather}\label{F on mc Z}  \lambda \star \Theta_{(s_i,s_{i+1})} =  \Theta_{(s_i,s_{i+1})} \quad   A_l^j\cdot \Theta_{(s_i,s_{i+1})}= \Theta_{(s_{i+j},s_{i+j+1})} \qquad \lambda \star \mc Z = \mc Z \qquad  A_l^j\cdot\mc Z=A_l^j.\end{gather}
	\end{lemma}
\bpr From \eqref{star properties 1} we see that the conditions  $s_i,s_{i+1} \in \sigma^{ss}$, $\phi_{\sigma}(s_i)<\phi_{\sigma}(s_{i+1})$ are equivalent to the conditions:  $s_i,s_{i+1} \in (z\star \sigma)^{ss}$,   $\phi_{z\star \sigma}(s_i)<\phi_{z\star \sigma}(s_{i+1})$, hence by \eqref{Theta_E n=2} we obtain the first equality. 

From Lemma \ref{the auto-eq F} by induction we get $A_l^j(s_i)\cong s_{i+j}$, $A_l^j(s_{i+1})\cong s_{i+j+1}$. Now  with the help of  \eqref{aut properties 1} and \eqref{Theta_E n=2} we  establish the second equality by a sequence of equivalences:
\begin{gather} \sigma \in \Theta_{(s_{i},s_{i+1})} \iff \begin{array}{c} (s_{i},s_{i+1})\subset \sigma^{ss}\\ \phi_\sigma(s_{i}) <\phi_\sigma(s_{i+1}) \end{array} \iff \begin{array}{c} ( A_l^j(s_{i}),A_l^j(s_{i+1}) )\subset \overline{A_l^j(\sigma^{ss})}=(A_l^j\cdot\sigma)^{ss}\\ \phi_{\sigma}((A_l^j)^{-1}(A_l^j(s_{i}))) <\phi_\sigma((A_l^j)^{-1}(A_l^j(s_{i+1}))) \end{array} \nonumber  \\ \iff
\begin{array}{c} ( s_{i+j},s_{i+j+1}) )\subset (A_l^j\cdot\sigma)^{ss}\\ \phi_{A_l^j\cdot\sigma}(s_{i+j}) <\phi_{A_l^j\cdot\sigma}(s_{i+j+1}) \end{array}\iff  A_l^j\cdot\sigma \in \Theta_{(s_{i+j},s_{i+j+1})}.  \nonumber 
\end{gather}  The third and the fourth equalities in \eqref{F on mc Z} follow from the already proven first and second.\epr
\begin{coro} \label{open cover of mc X} The action \eqref{left action} descents to an action  by biholomorphisms on $\mc X_l$ by the formula $\Phi\cdot \pr(\sigma)=\pr(\Phi\cdot\sigma)$. Denote $\Theta_i = \pr(\Theta_{(s_i,s_{i+1})})$ for $i\in\ZZ$  and $\mk{Z} =\pr(\mc Z)$.
 Then for any $i,j \in \ZZ$, $j\neq 0$ we have: {\rm (a)}  $ \Theta_{(s_i,s_{i+1})}= \pr^{-1}(\Theta_i)$; {\rm (b)} $\mc Z =\pr^{-1}(\mk{Z})$; {\rm (c)} $\{ \Theta_i\}_{i\in\ZZ}$ is an open cover of $\mc X_l$;  \begin{gather} \label{(d)} {\rm (d)} \qquad  \mk{Z} = \bigcap_{p\in \ZZ}\Theta_p =  \Theta_i \cap  \Theta_{i+j}      \qquad  A_l^j\cdot \Theta_i =  \Theta_{i+j}  \qquad   A_l^j\cdot\mk{Z} = \mk{Z}. \end{gather} 
 \end{coro} 
\bpr Since the actions \eqref{star}, \eqref{left action} commute  and $\pr$ from Corollary \ref{trvial bundle} has holomorphic sections, it follows the first sentence. The rest of the corollary follows  from   the definition of  quotient topology on  $\mc X_l$  and from  Lemmas \ref{actions on the blocks}, \ref{sections}. \epr

\subsection{The action \texorpdfstring{$\langle A_l \rangle$}{\space} on \texorpdfstring{$\mc X_l$ is free and properly discontinuous}{\space} for  \texorpdfstring{$l\geq 2$}{\space}} \label{final section}

In this Section we complete the proof of Theorem \ref{main}.

Corollary \ref{open cover of mc X} gives a decomposition of $\mc X_l$ into  $\langle A_l \rangle$-invariant subsets  $\mk{Z}$ and $\mc X_l\setminus \mk{Z}$ and furthermore  it gives a decomposition $\mc X_l\setminus \mk{Z} =  \amalg_{i\in \ZZ} \left ( \Theta_{i}\setminus \mk{Z} \right )$ into subsets which   $\langle A_l \rangle$ permutes, more precisely:  
 \begin{gather} \label{subactions}  \langle A_l \rangle \cdot  \mk{Z} =\mk{Z}  \qquad    \qquad A_l^j\cdot ( \Theta_{i}\setminus \mk{Z}) =   \Theta_{i+j}\setminus \mk{Z} \qquad  \mc X_l\setminus \mk{Z} =  \amalg_{i\in \ZZ} \left ( \Theta_{i}\setminus \mk{Z} \right ).   
\end{gather} 
 Now we construct  biholomorphisms  between $\mc H$,  $\Theta_j$ and $\mk{Z}$ and  describe the action of $\langle A_l \rangle$ on $\mk{Z}$.
\begin{lemma} \label{biholomorphisms} Let $j\in \ZZ$. Let $\bd \mc H & \rTo^\gamma & \mc S \ed$ be the map $\gamma(z)=(0,z)$ and let $\varphi_j$ be as in \eqref{mapping}. Then the function \eqref{biholo} is a biholomorphism. 
\begin{gather} \label{biholo} \bd  \mc H & \rTo^\gamma & \mc S & \rTo^{\varphi_j^{-1}} &\Theta_{(s_j,s_{j+1})}& \rTo^{\pr_{\vert}}& \Theta_{j} \ed.   \end{gather}
This function, restricted to the strip $\{z\in \CC: 0<\Im(z)<\pi\}$, is a biholomorphism between this strip and $\mk{Z}=\pr(\mc Z)$.  In particular, we obtain a biholomorphism:
\begin{gather} \label{biholom} \bd \mc H & \rTo^{\psi}& \mk{Z} \ed \qquad \psi(z)=  \pr\circ \varphi_j^{-1}\circ \gamma(\log{\abs{z}}+\ri \Arg(z)). \end{gather} We claim that for any $p\in \ZZ$ and any $z\in \mc H$ we have:
\begin{gather}\label{action1}  A_l^{-p}\cdot\psi(z)=\psi(\alpha_l^p(z)), \end{gather}
where $\alpha_l$ is the matrix in Lemma \ref{actionLEMMA} and $\alpha_l^p(z)$ is defined in \eqref{action}.
\end{lemma}
\bpr Since all composing  maps in \eqref{biholo}  are holomorphic, the composite  function is also holomorphic.  We will construct an inverse holomorphic function.  Let $\bd \mc S & \rTo^{pr_i}&\CC \ed $  be the projections  $pr_i(z_1,z_2)=z_i$, $i=1,2$ and let $\kappa_1$, $\kappa_2$ be the  holomorphic functions \eqref{kappas}, which obviously satisfy \eqref{props of kappas}:
\begin{gather} \label{kappas} \bd \CC & \lTo^{\kappa_1} & \Theta_{(s_j,s_{j+1})} & \rTo^{\kappa_2} &  \CC \ed \qquad \qquad  \kappa_1=pr_1\circ \varphi_j \quad \kappa_2=pr_2\circ \varphi_j \\ 
\label{props of kappas} \kappa_2\circ \ \varphi_j^{-1}\circ\gamma = Id_{\mc H} \qquad  \kappa_1\circ \ \varphi_j^{-1}\circ\gamma = 0.\end{gather} 
 One computes by straightforward  application of  \eqref{star properties 1}, \eqref{star properties 2}, and \eqref{mapping}:
\begin{gather}\label{kappa i proper} \kappa_1(\lambda \star \sigma) = \lambda + \kappa_1( \sigma) \quad \kappa_2(\lambda \star \sigma) = \lambda + \kappa_2( \sigma) \qquad \mbox{for} \quad \lambda \in \CC, \sigma \in \Theta_{(s_j,s_{j+1})}\end{gather}
therefore    $\kappa_1(\left (-\kappa_1(\sigma) \right ) \star \sigma) =0$ and  $\varphi_j(\left (-\kappa_1(\sigma) \right ) \star \sigma )=(0,\kappa_2(\left (-\kappa_1(\sigma) \right ) \star \sigma))\in \mc S$, hence  $\kappa_2(\left (-\kappa_1(\sigma) \right ) \star \sigma)\in \mc H$ and   we can define a holomorphic function  \eqref{kappa} satisfying \eqref{properties of kappa}, \eqref{properties of kappa1} (recall also \eqref{props of kappas}):
\begin{gather} \label{kappa} \bd \Theta_{(s_j,s_{j+1})} & \rTo^{\kappa}&\mc H \ed \qquad \kappa(\sigma) =  \kappa_2\left ( \left (-\kappa_1(\sigma) \right ) \star \sigma \right ) = \kappa_2\left (  \sigma \right )-\kappa_1(\sigma)\qquad \mbox{for} \quad \sigma \in  \Theta_{(s_j,s_{j+1})} \\ \label{properties of kappa} \varphi_j^{-1}\circ \gamma \circ \kappa (\sigma) = \left (-\kappa_1(\sigma) \right ) \star \sigma \qquad \kappa(\lambda \star \sigma) = \kappa( \sigma) \qquad  \qquad \mbox{for} \quad\lambda \in \CC, \sigma \in  \Theta_{(s_j,s_{j+1})} \\
\label{properties of kappa1}  \kappa\circ\ \varphi_j^{-1}\circ \gamma = Id_{\mc H}.  \end{gather} 
Due to the  equality $\kappa(\lambda \star \sigma) = \kappa( \sigma)$ we get a well (and uniquely) defined function \bd \bd \Theta_j & \rTo^{\kappa'} & \mc H \ed \qquad \kappa'\circ \pr_{\vert} = \kappa. \ed 
Furthermore, since $\pr_{\vert}$ is a  holomorphic $\CC$-principal bundle (Corollaries \ref{open cover of mc X}, \ref{trvial bundle}), it has holomorphic sections and it follows that $\kappa'$ is holomorphic.  We claim that $\kappa'$ is the desired inverse of \eqref{biholo}. 

Indeed,  using  \eqref{properties of kappa1} we obtain  $\kappa'\circ(\pr_{\vert}\circ \varphi_j^{-1}\circ \gamma )=\kappa \circ \varphi_j^{-1}\circ \gamma = Id_{\mc H}$.  Using  \eqref{properties of kappa} we compute $(\pr_{\vert}\circ \varphi_j^{-1}\circ \gamma ) \circ \kappa' \circ \pr_{\vert} = \pr_{\vert}$, which implies $(\pr_{\vert}\circ \varphi_j^{-1}\circ \gamma ) \circ \kappa'=Id_{\Theta_j}$, since $\pr_{\vert}$ is surjective. 

To consider  the restriction of \eqref{biholo} to  the strip we take $z\in \CC$ with  $0<\Im(z)<\pi$ and let  $\sigma =\varphi_j^{-1}\circ \gamma (z)$, then  $\varphi_j(\sigma)=(0,z)$ and \eqref{mapping}  shows that $s_j,s_{j+1} \in \sigma^{ss}$,  $\phi_\sigma(s_j)=0<\Im(z)/\pi=\phi_\sigma(s_{j+1})<1$,  therefore $\sigma \in \mc Z$ (see \eqref{intersection of two}), and hence $\pr(\sigma)\in \mk{Z}$.  Conversely, take any $\pr( \sigma) \in \mk{Z}$, $\sigma \in \mc Z$, then  $\phi_\sigma(s_j)<\phi_\sigma(s_{j+1})<\phi_\sigma(s_j)+1$ and by \eqref{kappa i proper}, \eqref{kappa}  we get $\kappa'(\pr(\sigma)) = \kappa(\sigma)=-\kappa_1(\sigma)+\kappa_2\left (  \sigma \right )$, therefore  \eqref{mapping}: $\Im(\kappa'(\pr(\sigma)))=\Im(-\kappa_1(\sigma)+\kappa_2(\sigma))=\pi(\phi_\sigma(s_{j+1})-\phi_\sigma(s_j))<\pi$. Thus we get   \eqref{biholom}.

 Let  $\sigma \in\mc Z$,  recalling that $\kappa'$ is the inverse of \eqref{biholo} and the definition of $\psi$ in \eqref{biholom},   we compute
\begin{gather} \psi^{-1}( \pr(\sigma))=\exp\circ\kappa'( \pr(\sigma))= \exp\circ\kappa( \sigma)=\exp\left ( \kappa_2(\sigma)-\kappa_1\left ( \sigma \right ) \right )\nonumber \\ 
=\mbox{by \eqref{mapping} } =
\exp\left (\log\abs{Z_{\sigma}(s_{j+1})} +\ri \pi \phi_{\sigma}(s_{j+1})-\log\abs{Z_{\sigma}(s_{j})} -
\ri \pi \phi_{\sigma}(s_{j})  \right ) =\nonumber \\ = \mbox{by \eqref{phase formula}} 
=\frac{Z_{\sigma}(s_{j+1})}{Z_{\sigma}(s_{j})}.\nonumber 
\end{gather}

 To show \eqref{action1} we first recall that  for $\sigma\in \mc Z$ we have $\{s_i\}_{i\in \ZZ}\subset \sigma^{ss}$, $\left \{ \frac{Z_\sigma(s_{i+1})}{Z_\sigma(s_i)} \right \}_{i\in \ZZ} \subset \mc H$ (Lemma \ref{sections}) and $A_l^p\cdot\sigma \in \mc Z$ for any $p\in\ZZ $  (see   \eqref{F on mc Z}), then using what we just computed,  we get (here $\sigma \in \mc Z$): 
\begin{gather} \psi^{-1}(A_l^{-p}\cdot \pr(\sigma))=\psi^{-1} (\pr(A_l^{-p}\cdot\sigma)))=\frac{Z_{A_l^{-p}\cdot\sigma}(s_{j+1})}{Z_{A_l^{-p}\cdot\sigma}(s_{j})}=
 \mbox{by   \eqref{aut properties 2}} 
= \frac{Z_{\sigma}(A_l^{p}(s_{j+1}))}{Z_{\sigma}(A_l^{p}(s_{j}))}\nonumber  \\ 
=\mbox{by Corollary \ref{the auto-eq F}} = \frac{Z_{\sigma}(s_{p+j+1})}{Z_{\sigma}(s_{p+j})}=
\mbox{by Lemma \ref{actionLEMMA} } = \alpha_l^p\left (\frac{Z_{\sigma}(s_{j+1})}{Z_{\sigma}(s_{j})}\right) = \alpha_l^p(\psi^{-1}(\pr(\sigma)) ).\nonumber
\end{gather} Now \eqref{action1} follows from the obtained $\psi^{-1}(A_l^{-p}\cdot \pr(\sigma))= \alpha_l^p(\psi^{-1}(\pr(\sigma)) )$,  since $\psi$ is biholomorphism.
\epr 
 The matrix  $ \alpha_l =\begin{bmatrix}
	l & -1 \\
	1 & 0
\end{bmatrix} $   has $\tr(\alpha_l)=l$.  Since $l\geq 2$, then $\alpha_l$ is either parabolic or hyperbolic (see Section \ref{hpe}).  As a   corollary following from \eqref{subactions}, Remark \ref{prop disc free} (a) and \eqref{action1} we get:
\begin{coro} \label{partially prop disc} The action of $ \langle A_l \rangle $ on $\mc X_l$ is free and its restriction to $\mk{Z}$ is properly discontinuous. 
\end{coro} 

\begin{figure} \centering
   \hspace{-40mm}   \begin{subfigure}[b]{0.2\textwidth}
       \begin{tikzpicture}
\path [fill=yellow] (-1,0) --(0,0) to [out=90,in=180] (1.5,1.5) to [out=0,in=90] (3,0) -- (4.5,0) -- (4.5,3)--(-1,3)--(-1,0);
\draw [<-] (-1,0)--(0,0);
\draw [->] (0,0)--(5.3,0);
\draw [<->] (4.5,0)--(4.5,3);
\draw [dashed, ->] (0, 0) -- (0,3);
\draw  [<->]    (3,0) arc [radius=1.5, start angle=0, end angle= 180];
\node [below] at (0,0) {$0$};
\node [below] at (3.8,0) {$1$};
\draw (3.8,-.05) -- (3.8, .05);
\node [below] at (3,0) {$\frac{2 a_l}{a_l^2 +1}$};
\node [below] at (4.5,0) { $\frac{a_l^2 +1}{2a_l}$ };
\node [right] at (5.3,0) {$x$};
\node [left] at (0,3) {$y$};
\node [right] at (4.5,1.5) {$\mk{b}_+$};
\node [below] at (1.5,1.5) {$\mk{b}_-$};
\node      at (1.5,2.3) {$F_l$};
\node [below] at (1.5,0) {$\frac{a_l}{a_l^2 +1}$};
\draw (1.5,-.05) -- (1.5, .05);
\end{tikzpicture}
        \caption{$l\geq3$}
        \label{FDhyperbolic}
    \end{subfigure} 
		\hspace{50mm} \begin{subfigure}[b]{0.2\textwidth}
   \begin{tikzpicture}
\path [fill=yellow] (-1,0) --(0,0) to [out=90,in=180] (2,2) to [out=0,in=90] (4,0) --(4,3)--(-1,3)--(-1,0);
\draw [<-] (-1,0)--(0,0);
\draw [->] (0,0)--(5.3,0);
\draw [->] (4,0)--(4,3);
\draw [dashed, ->] (0, 0) -- (0,3);
\draw   [<->]    (4,0) arc [radius=2, start angle=0, end angle= 180];
\node [below] at (0,0) {$0$};
\node [below] at (4,0) {$1=\frac{a_2^2 +1}{2a_2}=\frac{2 a_2}{a_2^2 +1}$};
\node  at (2,2.5) {$F_l$};
\node [right] at (4,2) {$\mk{b}_-$};
\node [below] at (2,2) {$\mk{b}_+$};
\node [right] at (5.3,0) {$x$};
\node [left] at (0,3) {$y$};
\end{tikzpicture}
        \caption{$l=2$}
        \label{FDparabolic}
    \end{subfigure}
     \caption{Fundamental domains of $\langle \alpha_l \rangle $}\label{1FD}
\end{figure}
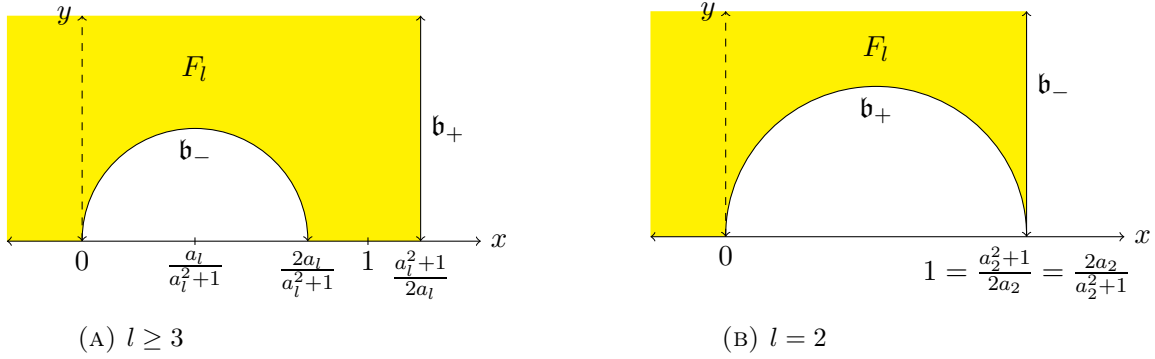

The next step is to find a fundamental domain of the  action of $ \langle A_l \rangle $ on $\mk{Z}\cong \mc H$. We choose $j\in \ZZ$ and use  the biholomorphism $\psi:\mc H \rightarrow  \mk{Z}$ \eqref{biholom} (depending on $j$) constructed in  Lemma \ref{biholomorphisms}. By \eqref{action1} we see that the fundamental domain we want is of the form $\psi(F_l)$, where  $F_l\subset \mc H$ is a fundamental domain of  the action of $\langle \alpha_l \rangle$ on $\mc H$ discussed in Section \ref{recollection}.
 We  use now  Sections \ref{0FD remark hyperbolic} and  \ref{0FD remark parabolic} to find such a $F_l$.
\begin{lemma}  \label{lemma FDhyperbolic} Let us denote  $a_l=\frac{l+\sqrt{l^2-4}}{2}$  for $l\geq 2$, $a_2=1$.  A fundamental domain of $\langle \alpha_l \rangle$ is: 
\begin{gather} \label{F_l} F_l=\left \{x+\ri y \in \mc H : \begin{array}{c}  \frac{a_l^2+1}{2 a_l}\left ( x^2+y^2 \right ) \geq  x \leq \frac{a_l^2+1}{2 a_l}\end{array} \right \} \nonumber \\[-2mm]  \\[-2mm] \nonumber {\rm Bd}_{\mc H}(F_l)=\left \{x+\ri y \in \mc H : \begin{array}{c}  \frac{a_l^2+1}{2 a_l}\left ( x^2+y^2 \right ) =  x  \ \text{or} \  x= \frac{a_l^2+1}{2 a_l}\end{array} \right \}.
\end{gather}
For $l\geq 3$ the fundamental domain is  shown in Figure \eqref{FDhyperbolic} and for  $l=2$ in Figure  \eqref{FDparabolic}.
 
Two points in $F_l$   lie  in a  common orbit iff they satisfy $z_{\pm} \in {\rm Bd}_{\mc H}(F_l)$,   $\Arg(z_+) =\Arg(z_-)$.

 For each $u\in {\rm Bd}_{\mc H}(F_l)$ there exists  an open subset $U\subset \mc H$, s. t. $u\in U$ and
 $\{i\in \ZZ: \alpha^i(U)\cap F_l\neq \emptyset\}$ is finite (in fact contains only two elements).
\end{lemma}
\bpr  We need to consider two cases $l\geq 3$ and $l=2$, because this determines the type of $\alpha_l$ according to the classification recalled  in Subsection \ref{hpe}.

Assume first that $l\geq 3$.  The  feature of $a_l$ we need, which we have  for all $l\geq 3$, is: $a_l >1$, that's why we will omit the subscript $l$ and write just $a$, remembering that $a>1$.  Now $\alpha_l$ is  a  hyperbolic element in ${\rm SL(2,\ZZ)}$ (see Subsection \ref{hpe}) and we can use the method  described in Subsection \ref{0FD remark hyperbolic}.  We note first that
$ \beta^{-1} \  \alpha_l \ \beta=\begin{bmatrix}
	a & 0 \\
	0 & a^{-1}
\end{bmatrix},\nonumber $
where 
\be a=\frac{l+\sqrt{l^2-4}}{2}  \qquad \beta^{-1} =\frac{1}{a-a^{-1}}\begin{bmatrix}
	1 & -a^{-1}  \\
	-1 & a
\end{bmatrix} \quad \beta= \begin{bmatrix}
	a & a^{-1} \\
	1 & 1
\end{bmatrix}. \ee
Since $a>1$, it follows that $\det(\beta)>0$ and  $\beta$ determines a biholomorphism  
of $\mc H$ (recall \eqref{action}):
\be \beta(z)=\frac{a^2 z +1}{a (z+1)}  \qquad z \in \mc H. \ee 

We choose the parameters  of the strip  $F'$ in Figure \eqref{0FDhyperbolic} to be $\delta=1/a^2$, $\delta a^2=1$. Since $\beta(1)=\frac{a^2+1}{2 a}$, $\beta(-1)=\infty$, by Remark \ref{circles and lines} we see that $\beta(\mk{b}_+')=\mk{b}_+$ (see Figures \eqref{0FDhyperbolic} and \eqref{FDhyperbolic}). Since $\beta(\frac{-1}{a^2})=0$, $\beta(\frac{1}{a^2})=\frac{2 a}{a^2 +1}$ and using  Remark \ref{circles and lines} again,  we deduce that $\beta(\mk{b}_-')=\mk{b}_-$. Finally, since $\beta(1/a)=1$ we deduce that $\beta(F')=F_l$, where $F'$, $F_l$ are depicted in Figures \eqref{0FDhyperbolic} and \eqref{FDhyperbolic}. The first part of  the lemma is proved for $l\geq 3$.  Two points in $F_l$ lie in a  common  orbit iff they are of the form $\beta(a^2 z), \beta(z), z \in  \mk{b}_-'$.
  For $z\in \mk{b}_-'$ we have $\abs{z}^2=1/a^4$ and we compute:
\begin{gather} \frac{\beta(a^2 z)}{\beta(z)}=\frac{(a^4 z +1)(z+1)}{(a^2 z +1)^2}=\frac{(a^4 z +1)(z+1)}{(a^2 z +1)^2}\frac{\ol{z}^2}{\ol{z}^2}=\frac{(1 +\ol{z})(\frac{1}{a^4}+\ol{z})}{(\frac{1}{a^2} +\ol{z})^2} 
=\frac{(1 +\ol{z})(1+a^4 \ol{z})}{(1 + a^2 \ol{z})^2}=\frac{\ol{\beta(a^2 z)}}{\ol{\beta(z)}}\nonumber \end{gather}
hence   $\beta(a^2 z) , \beta( z) $ are parallel   for $z\in \mk{b}_-'$ and, being in a common quadrant, $\beta(a^2 z)\in \RR^{>0} \ \beta( z)$. 

Assume  that $l=2$.  In this case  $\alpha_2$ is parabolic and   we can use the method  described in Subsection \ref{0FD remark parabolic}.  We note first that
$ \beta^{-1} \  \alpha_2 \ \beta=\begin{bmatrix}
	1 & -1 \\
	0 & 1
\end{bmatrix},\nonumber $
where 
$  \beta^{-1} =\begin{bmatrix}
	\frac{1}{2} & 	\frac{1}{2}  \\
	-1 & 1
\end{bmatrix} \quad \beta= \begin{bmatrix}
	1 & -	\frac{1}{2} \\
	1 & 	\frac{1}{2}
\end{bmatrix}. $
And now  $\beta(F')$ is a fundamental domain, where $F'$ is of the form in Figure \eqref{0FDparabolic} and $\beta$ is:
\be \beta(z)=\frac{2 z-1}{2 z +1}  \qquad z \in \mc H. \ee 

We choose the parameters  of the strip  $F'$ in Figure \eqref{0FDhyperbolic} to be $\delta=-1/2$, $\delta +\abs{b/a}=1/2$. Since $\beta(1/2)=0$, $\beta(\infty)=1$, $\beta(-1/2)=\infty$, by Remark \ref{circles and lines} we see that $\beta(\mk{b}_\pm')=\mk{b}_\pm$ (see Figures \eqref{0FDparabolic} and \eqref{FDparabolic}).  Finally, since $\beta(0)=-1$, we deduce that $\beta(F')=F_2$,  where $F'$, $F_2$ are depicted in Figures \eqref{0FDparabolic} and \eqref{FDparabolic}. It remains to describe the orbits in  $F_2$. Two points in $F_2$  lie in a common orbit   iff they are of the form $\beta(z-1), \beta(z), z\in \mk{b}_+'$.
  For $z\in \mk{b}_+'$ we have $\Re(z) =1/2$ and with the  computation below the lemma is proved (recall also Remark \ref{prop disc free} (b)):
\begin{gather} \frac{\beta( z)}{\beta(z-1)}=\frac{(2 z-1)^2}{(2 z+1) (2 z -3)}=\frac{(\ri 2 \Im(z))^2}{(2+2 \ri \Im(z)) (2 \ri \Im(z)-2)}=\frac{\Im(z)^2}{1+\Im(z)^2}>0.\nonumber \end{gather}
\epr

\begin{coro} \label{intermediate coro} A fundamental domain $\mc F_l'$ of the action $\langle A_l \rangle$ on $\mk{Z}$ is the image of the set: 
\begin{gather}  \left \{ (u+\ri v)\in  \CC : 0 < v < \pi \ \ \text{\rm and} \ \   \re ^{u+\Delta_l} \geq  \cos(v) \leq  \re^{-u+\Delta_l}   \right \} \nonumber \\[-3mm] \label{mc Fprim}\\[-3mm]
=\left \{(u+\ri v)\in  \CC: v\in (0,\pi ) \ \ \text{\rm and}  \begin{array}{c}  v\geq \arccos\left (\min\left \{1,\re^{\Delta_l-\abs{u}} \right \} \right )\end{array} \right \} \nonumber\end{gather}
by the biholomprhism \eqref{biholo},  where $\Delta_l=\log\left (  \frac{a_l^2+1}{2 a_l}\right )$.  The following specifications hold: 

{\rm (a)} The boundary  ${\rm Bd}_{\mk{Z}}(\mc F_l')$ is the image of $\{ \pm u + \ri \arccos\left (\re^{\Delta_l-\abs{u}} \right ) : \abs{u}>\Delta_l \} $ by  \eqref{biholo}.

{\rm (b)} Two points in $\mc F_l'$ lie in the same orbit of the action iff they are in  ${\rm Bd}_{\mk{Z}}(\mc F_l')$ and have the same
  imaginary parts via the inverse of \eqref{biholo}. 

{\rm (c)} For any $q\in {\rm Bd}_{\mk{Z}}(\mc F_l')$ there exists  an open $U\subset \mk{Z}$, s. t. $q\in U$ and
 $\{i\in \ZZ: A^i(U)\cap \mc F_l' \neq \emptyset\}$ is

  \hspace{4mm}  finite (in fact contains only two elements).
\end{coro}
\bpr Recalling  \eqref{biholo},  \eqref{biholom} we see that $\mc F_l'=\psi(F_l)=\pr_{\vert}\circ \varphi_j^{-1} \circ \gamma \left ( \exp_{\vert \{z\in\mc H:\Im(z)<\pi\}\rightarrow \mc H}^{-1}(F_l)\right ) $. So we have to determine $\exp_{\vert \{z\in\mc H:\Im(z)<\pi\}\rightarrow \mc H}^{-1}(F_l)$. An element $(u+\ri v)\in \CC$ is in the latter  set iff $0<v<\pi$ and $\exp(u+\ri v)\in F_l$, which by Lemma\ref{F_l} is the same as: $\frac{a_l^2+1}{2 a_l}\re ^{2 u} \geq  \cos(v) \re^u \leq \frac{a_l^2+1}{2 a_l}$, so we get
\begin{gather} \mc F_l'= \pr_{\vert}\circ \varphi_j^{-1} \circ \gamma\left ( \left \{ (u+\ri v)\in \mc H : v < \pi \ \  \frac{a_l^2+1}{2 a_l}\re ^{u} \geq  \cos(v) \leq  \re^{-u}  \frac{a_l^2+1}{2 a_l} \right \} \right ) \end{gather}
and the corollary follows from \eqref{action1}, the fact that $\psi$ is biholomorphism,  and the properties of $F_l$ obtained in   Lemma \ref{lemma FDhyperbolic}. 
\epr

\begin{coro} \label{fundamental domain}
Let $j\in \ZZ$. Let $\bd \Theta_j & \rTo^{u+\ri v}& \mc H\ed $ be the inverse  of \eqref{biholo}.  Let $\Delta_l=\log\left (  \frac{a_l^2+1}{2 a_l}\right )$. The subset   $\mc F_l =\left \{ q\in \Theta_j :  v(q)\geq \arccos\left (\min\left \{1,\re^{\Delta_l-\abs{u(q)}} \right \} \right ) \right \}$ has the following properties:

{\rm (a)} ${\rm Bd}_{\mk{Z}}(\mc F_l')={\rm Bd}_{\Theta_j}(\mc F_l) =  \left \{ q\in   \Theta_j :    v(q)= \arccos\left (\min\left \{1,\re^{\Delta_l-\abs{u(q)}} \right \} \right )\right \}\subset \mk{Z}$. 

	{\rm (b)} The interior of $\mc F_l$ w.r. to $\mc X_l$ is   $\mc F_l^o=  \left \{ q\in   \Theta_j :    v(q)> \arccos\left (\min\left \{1,\re^{\Delta_l-\abs{u(q)}} \right \} \right )\right \}$. 
	
	{\rm (c)} $\mc X_l = \cup_{i\in \ZZ} A_l^i(\mc F_l) $
	
	{\rm (d)}  $A_l^m(\mc F_l^o)\cap A_l^n(\mc F_l^o)=\emptyset$ for $m\neq n$.

  {\rm (e)}  Two points in $\mc F_l$ lie in the same orbit of the $\langle A_l \rangle$-action on $\mc X_l$  iff they are in  ${\rm Bd}_{\Theta_j}(\mc F_l)$ and have the   same imaginary parts (i.e. $v$ has the same values on them). 
	
	 {\rm (f)}  For each $q\in {\rm Bd}_{\Theta_j}(\mc F_l)$ there exists  an open subset $U\subset \mc X_l$, s. t. $q\in U$ and
 $\{i\in \ZZ: A^i(U)\cap \mc F_l \neq \emptyset\}$ is finite (in fact contains only two elements).
\end{coro}
\bpr (a): It is clear that 
 ${\rm Bd}_{\Theta_j}(\mc F_l) = \left \{ q\in   \Theta_j :    v(q)= \arccos\left (\min\left \{1,\re^{\Delta_l-\abs{u(q)}} \right \} \right )\right \}$,  since \eqref{biholo} is biholomorphism.  ${\rm Bd}_{\mk{Z}}(\mc F_l')$ is the same set by Corollary \ref{intermediate coro} (a). 

(b): Since $\Theta_j$ is open subset in $\mc X_l$, the interiors  of $\mc F_l$ w.r. to $\Theta_j$ and w.r. to $\mc F_l$ coincide. Hence (b) is due to the fact that \eqref{biholo} is biholomorphism. 

(c) and (d): From Lemma \ref{biholomorphisms} and the definition of $\mc F_l'$  in Corollary \ref{intermediate coro} we see that \begin{gather} \label{intermidiate help} \mc F_l=\mc F_l' \cup (\Theta_j\setminus \mk{Z}) \qquad \qquad  \mc F_l^o=\mc F_l'^o \cup (\Theta_j\setminus \mk{Z}). \end{gather} Therefore by \eqref{subactions}  and since  $\mc F_l'$ is fundamental domain of the action of $\langle A_l \rangle$ on $\mk{Z}$ (as defined on  \cite[p. 20]{Miyake}) we get: $\cup_{i\in \ZZ} A_l^i(\mc F_l) = \cup_{i\in \ZZ} A_l^i(\mc F_l') \cup \cup_{i\in \ZZ} A_l^i(\Theta_j\setminus \mk{Z})   = \mk{Z} \cup (\mc X_l \setminus \mk{Z})=\mc X_l$ and $A_l^m(\mc F_l^o)\cap A_l^n(\mc F_l^o)  =  (A_l^m(\mc F_l'^o) \cup A_l^m(\Theta_j\setminus \mk{Z}) )\cap  (A_l^n(\mc F_l'^o) \cup A_l^n(\Theta_j\setminus \mk{Z}) )  =\emptyset$ for $m\neq n$.

(e): 
From \eqref{intermidiate help} and \eqref{subactions} we see that if two points in $\mc F_l$ lie in a common orbit, then they lie in $\mc F_l'$, and then we use the already proven (a) here and  Corollary \ref{intermediate coro}  (b) to obtain (e). 

  (f): For $q\in  {\rm Bd}_{\Theta_j}(\mc F_l)$ the same neighborhood $U\ni q$ ensured by Corollary \ref{intermediate coro}  (c) satisfies the required property in (f) due to  \eqref{subactions}, \eqref{intermidiate help},  and the already proven (a) here. \epr 
\begin{figure} \centering
   \hspace{-40mm}   \begin{subfigure}[b]{0.2\textwidth}
       \begin{tikzpicture}
\path [fill=yellow] (-3,0.9)  to [out=350,in=135] (-0.9,0) -- (0.9,0) to [out=45,in=190] (3,0.9)--(3,2)--(-3,2)--(-3,0.9);
\path [fill=purple] (-3,2)  -- (3,2) -- (3,3)--(-3,3)--(-3,2);
\draw [<->] (-3,0.9)  to [out=350,in=135] (-0.9,0) ;
\draw [<->] (0.9,0)  to [out=45,in=190] (3,0.9);
\draw [dashed, ->] (0, 0) -- (0,3);
\node [left] at (0,3) {$v$};
\draw[<->] (-3,0)--(3,0);
\draw[dashed, <->] (-3,1)--(3,1);
\draw[dashed, <->] (-3,2)--(3,2);
\node [right] at (3,1) {$v=\frac{\pi}{2}$};
\node [right] at (3,2) {$v=\pi$};
\node [right] at (3,0) {$u$};
\node [below] at (0,0) {$0$};
\node [below] at (0.9,0) {$\Delta_l$};
\node [below] at (-0.9,0) {$-\Delta_l$};
\end{tikzpicture}
        \caption{$l\geq3$}
        \label{2FDhyperbolic}
    \end{subfigure} 
		\hspace{50mm} \begin{subfigure}[b]{0.2\textwidth}
    \begin{tikzpicture}
\path [fill=yellow] (-3,0.9)  to [out=350,in=90] (0,0)  to [out=90,in=190] (3,0.9)--(3,2)--(-3,2)--(-3,0.9);
\path [fill=purple] (-3,2)  -- (3,2) -- (3,3)--(-3,3)--(-3,2);
\draw [<->] (-3,0.9)  to [out=350,in=90] (0,0) ;
\draw [<->] (0,0)  to [out=90,in=190] (3,0.9);
\draw [dashed, ->] (0, 0) -- (0,3);
\node [left] at (0,3) {$v$};
\draw[<->] (-3,0)--(3,0);
\draw[dashed, <->] (-3,1)--(3,1);
\draw[dashed, <->] (-3,2)--(3,2);
\node [right] at (3,1) {$v=\frac{\pi}{2}$};
\node [right] at (3,2) {$v=\pi$};
\node [right] at (3,0) {$u$};
\node [below] at (0,0) {$0$};
\end{tikzpicture}
        \caption{$l=2$}
        \label{2FDparabolic}
    \end{subfigure}
     \caption{$(u+\ri v)(\mc F_l)=\mk{F}_l$}\label{FD}
\end{figure}

\begin{coro} \label{properly discontinuous} The action of $\langle A_l \rangle$ on $\mc X_l$ is free  and properly discontinuous. 
\end{coro}
\bpr In Corollary \ref{partially prop disc} we already showed that the action is free. Take any two $q_1, q_2 \in \mc X_l$.  We need to find open subsets $U_i\ni q_i$ in $\mc X_l$, $i=1,2$, such that  the set $\{i\in \ZZ: A_l^i(U_1)\cap U_2\neq \emptyset\}$ is finite (see e.g. \cite[p. 17]{Miyake}). From  (c) in Corollary \ref{fundamental domain} it follows that it is enough to consider the case $q_1, q_2 \in \mc F_l$.\footnote{because the set $\{i\in \ZZ: A_l^i(A_l^m U_1)\cap A_l^n U_2\}$ is $n-m+\{i\in \ZZ: A_l^i( U_1)\cap U_2\}$} If $q_1, q_2 \in {\rm Bd}_{\Theta_j}(\mc F_l)$, then $q_1, q_2 \in \mk{Z}$ (see Corollary \ref{fundamental domain} (a))  and the neighborhoods $U_i\ni q_i$ we  need  exist since $\mk{Z}$ is an $\langle A_l \rangle$-invariant open subset and  by Corollary \ref{partially prop disc}. If $q_1,q_2\in \mc F_l^o$, then $U_1=U_2= F_l^o$ satisfy the condition we need by Corollary \ref{fundamental domain} (d). Since $\mc F_l = F_l^o \cup {\rm Bd}_{\Theta_j}(\mc F_l) $, it remains to consider the case $q_1 \in  {\rm Bd}_{\Theta_j}(\mc F_l)$, $q_2\in  F_l^o$. In this case we take $U\ni q_1$ as in Corollary \ref{fundamental domain} (f) and put $U_1=U$, $U_2=\mc F_l^o$, and the corollary is proved.   \epr 

\begin{coro} \label{covering} The orbit-space $\mc X_l/\langle A_l \rangle$ with the quotient topology carries a structure of a one dimensional complex manifold, s. t. the projection $\wt{\pr}  : \mc X_l \rightarrow \mc X_l/\langle A_l \rangle$ is a holomorphic covering map: a universal cover of $\mc X_l/\langle A_l \rangle$. In particular $\pi_1(\mc X_l/\langle A_l \rangle)=\ZZ$.
\end{coro}
\bpr  In the previous corollary we showed that the action  is free and proper. From these properties    \cite[Proposition 1.2]{Huck} ensures that  $\mc X_l/\langle A_l \rangle$ has the structure of a  one dimensional complex manifold, s. t. $\wt{\pr}  : \mc X_l \rightarrow \mc X_l/\langle A_l \rangle$ is a holomorphic  principal $\ZZ$-bundle and recalling that  $ \mc X_l$ is contractible (see Corollary \ref{trvial bundle})  the corollary follows from the long exact sequence for the homotopy groups associated to $\wt{\pr}$.  
\epr
 From now on we fix  $j\in \ZZ$ and let $\mc F_l\subset \Theta_j$ be the closed in  $\Theta_j$ subset obtained in Corollary \ref{fundamental domain}. By (c) in Corollary \ref{fundamental domain} it follows that $\wt{\pr}(\mc F_l)=\wt{\pr}(\Theta_j)=\mc X_l/\langle A_l \rangle$, i. e. we have the surjectivity of the restriction of the map  $\wt{\pr}$ from Corollary \ref{covering}:

\begin{lemma} \label{proper map}The restriction $\wt{\pr}_{\vert \mc F_l }: \mc F_l \rightarrow \mc X_l/\langle A_l \rangle$  is a proper surjective map.
\end{lemma}
\bpr It remains to show the properness. Recall that  $\Theta_j$ is an open subset in $\mc X_l$, and $\mc F_l$ is a closed subset in  $\Theta_j$. Thus, the restriction $\bd \Theta_j & \rTo^{\wt{\pr}_{\vert \Theta_j }} & \mc X_l/\langle A_l \rangle \ed $ is a local biholomorphism between locally compact spaces. Let $K\subset X_l/\langle A_l \rangle$ be  compact. 

For any $q\in K$ we fix $q'\in \mc F_l$, s. t. $\wt{\pr}(q')=q$.  Now we will choose an open $U_q\ni q'$ subset of $\Theta_j$ with certain properties for each $q\in K$.  If $q'\in \mc F_l^o$ we choose  $U_q\ni q'$ whose closure w. r. to $\Theta_j$ is compact and contained in $\mc F_l^o$ (in particular $A_l^j(U_q)\cap \mc F_l=\emptyset$ for $j\neq 0$). If $q'\in {\rm Bd}_{\Theta_j}({\mc F}_l)$ we choose  $U_q\ni q'$ with compact closure  w. r. to $\Theta_j$, which is contained in $\mk{Z}$ and such that  $\{i\in \ZZ: A_l^i(U_q)\cap \mc F_l \neq \emptyset\}$ is finite (we can do this by Corollary \ref{fundamental domain} (f) and since ${\rm Bd}_{\Theta_j}({\mc F}_l)\subset \mk{Z}$). Since $K$ is compact and $\wt{\pr}$ is open, we have $K\subset \cup_{i=1}^n \wt{\pr}(U_{q_i})$ for a finite family $\{q_i\}_{i=1}^n\subset K$, i.e.
\begin{gather}(\wt{\pr}_{\vert \mc F_l})^{-1}(K)\subset \cup_{i=1}^n\wt{\pr}^{-1}(\wt{\pr}(U_{q_i}))\cap \mc F_l=  \cup_{i=1}^n\cup_{m\in\ZZ} A_l^m(U_{q_i})\cap \mc F_l. \nonumber \end{gather} 
By our choice of $U_q$ the union on the right hand side is finite and each element $A_l^m(U_{q_i})$  in it is contained in a compact subset of $\Theta_j$ (recall also \eqref{subactions}). Since $\mc F_l$ is closed in $\Theta_j$, we deduce that $(\wt{\pr}_{\vert})^{-1}(K)$ is  contained in a compact subset of $\mc F_l$, and  therefore, being closed subset of compact, $(\wt{\pr}_{\vert})^{-1}(K)$ is compact. \epr
Now we can prove:

\begin{prop} \label{final proof} If  $l\geq 3$, then  $\mc X_l$ is biholomorphic to $\mc H$, and  Corollary \ref{trvial bundle} implies Theorem \ref{main}.
\end{prop}
\bpr Suppose that $\mc X_l$ is not biholomorphic to $\mc H$. We will obtain a contradiction. 

By  Corollary \ref{trvial bundle} we see that we have a biholomorphism $\mc X_l \cong \CC$, and  we showed  in Corollary \ref{covering}   that $\wt{\pr}: \mc X_l \rightarrow \mc X_l / \langle A_l \rangle$ is a universal covering of $\mc X_l / \langle A_l \rangle$ and $\pi_1(\mc X_l / \langle A_l \rangle)=\ZZ$.
   \cite[Prop. 27.12]{Otto} ensures that $\mc X_l / \langle A_l \rangle$ is biholomorphic to one of these: $\CC,\CC^{\star}$, or a torus. However  $\pi_1(\mc X_l / \langle A_l \rangle)=\ZZ$ and we deduce  that  $\theta: \mc X_l / \langle A_l \rangle \cong  \CC^{\star}$ for some biholomorphism $\theta$. 	Let us denote by $f$ the following composition: \begin{gather} \label{f} f: \bd  \mc H & \rTo^\gamma & \mc S & \rTo^{\varphi_j^{-1}} &\Theta_{(s_j,s_{j+1})}& \rTo^{\pr_{\vert}}& \Theta_{j} & \rTo^{\wt{\pr}_{\vert}} & \mc X_l / \langle A_l \rangle & \rTo^{\theta} & \CC^{\star}. \ed \end{gather} As in Corollary \ref{fundamental domain} we denote  by $\bd \Theta_j & \rTo^{u+\ri v}& \mc H\ed $  the inverse  of \eqref{biholo}. 
	Let us denote $\mk{F}_l=(u+\ri v)(\mc F_l)$ (see Figure \eqref{2FDhyperbolic}).  From Corollary \ref{fundamental domain}  (or Corollary \ref{intermediate coro}) (a) we see that:
\begin{gather}{\rm Bd}_{\mc H}(\mk{F}_l)=\{\pm \delta(t)+\ri t: t\in (0,\pi/2) \}, \qquad \text{where} \qquad\delta(t)=(\Delta_l-\log(\cos(t))).  \end{gather} Then from Lemma \ref{biholomorphisms}, Corollary \ref{fundamental domain}, and Lemma \ref{proper map}  follow (a), (b), (c) below: 

(a) $f:\mc H \rightarrow \CC^{\star}$ is a local biholomorphism;

(b) $f(\mk{F}_l)=\CC^{\star}$  and $\bd \mk{F}_l& \rTo^{f_{\vert \mc F_l}} & \CC^{\star}\ed$ is proper;

(c)  $q_\pm \in \mk{F}_l$ and $f(q_+)=f(q_-)$ implies $q_+ = q_-$ or $q_{\pm}=\pm \delta(t)+\ri t$ for some $t\in (0,\pi/2)$.

For any $t \in (0,\pi/2)$ let us  denote  the segment  $\gamma_t$ in $\mc F_l$ defined by $\gamma_t(s) =s+\ri t$ for $s\in [-\delta(t),\delta(t)]$. Then  (see Figure \eqref{3FDhyperbolic})
\begin{gather} \label{F_l^tpm} \mk{F}_l = \mk{F}_l^{t-} \amalg \im(\gamma_t) \amalg \mk{F}_l^{t+}, \ \ \mbox{where} \ \ \mk{F}_l^{t-}=\{z\in \mk{F}_l:\Im(z)<t\},   \ \ \mk{F}_l^{t+}=\{z\in \mk{F}_l:\Im(z)>t\}.  \end{gather}    
From (c) above  it follows that $f(\mc  F_l^{t+})$, $f(\im(\gamma_t))$, $f(\mk{F}_l^{t-})$ are pairwise disjoint  and we can write : \begin{gather}\label{a disj union} \CC^{\star}=f( \mk{F}_l^{t+}) \amalg f(\im(\gamma_t)) \amalg f( \mk{F}_l^{t-}) \qquad f_{\vert \mk{F}_l}^{-1}(f( \mk{F}_l^{t\pm}) )=\mk{F}_l^{t\pm} \qquad f_{\vert \mk{F}_l}^{-1}(f( \im(\gamma_t)) )=\im(\gamma_t). \end{gather}
 
Due to (c) above $f\circ \gamma_t$ is a closed Jordan curve in $\CC^\star$. By Jordan curve theorem (see e.g. \cite[Theorem 4.13 on p. $70_1$]{Mark}) we have $\CC \setminus f(\im(\gamma_t)) = U_t\amalg V_t$, where $U_t, V_t$ are the connected components of the complement of $\im(f\circ \gamma_t)$ in $\CC$.  In particular, since    $f(\mk{F}_l^{t\pm})$ are connected, by \eqref{a disj union} we get either  $f(\mk{F}_l^{t+})\subset U_t\setminus \{0\}$, $f(\mk{F}_l^{t-})\subset V_t\setminus\{0\}$  or $f(\mk{F}_l^{t-})\subset U_t\setminus\{0\}$, $f(\mk{F}_l^{t+})\subset V_t\setminus\{0\}$.  Due to the first equality in  \eqref{a disj union} we can write:
$ f( \mk{F}_l^{t-}) \amalg  f( \mk{F}_l^{t+}) =(U_t \setminus \{0\}) \amalg  (V_t\setminus \{0\}) 
 $  and we obtain:
\begin{gather} \label{alternative equalities} \left \{f(\mk{F}_l^{t+}) = U_t\setminus \{0\} \ \text{and} \  f(\mk{F}_l^{t-})= V_t\setminus\{0\} \right \} \ \text{or} \   \left \{  f(\mk{F}_l^{t-})= U_t\setminus\{0\} \ \text{and} \ f(\mk{F}_l^{t+})= V_t\setminus\{0\} \right \}. \end{gather}

   Furthermore, by Jordan theorem  we can assume that   $U_t \cup f(\im (\gamma_t))$ is compact in $\CC$, $V_t \cup \im (\gamma_t)$ is closed,  non-compact in $\CC$.  We claim that $0\in U_t$. Indeed, if $0\not \in U_t$, then  $U_t\setminus \{0\}=U_t$ and by \eqref{alternative equalities} either $f(\mk{F}_l^{t+}\cup \im(\gamma_t))$ or  $f(\mk{F}_l^{t-}\cup \im(\gamma_t))$ is compact, which in turn by \eqref{a disj union} and (b) above leads to  either $\mk{F}_l^{t+}\cup \im(\gamma_t)$ or  $\mk{F}_l^{t-}\cup \im(\gamma_t)$ is compact, which is a contradiction (see \eqref{F_l^tpm}). 
	So we see that $0\in U_t$, $V_t= V_t\setminus\{0\}$, and
	$ \left \{f(\mk{F}_l^{t+}) = U_t\setminus \{0\} \ \text{and} \  f(\mk{F}_l^{t-})= V_t \right \} \ \text{or} \   \left \{  f(\mk{F}_l^{t-})= U_t\setminus\{0\} \ \text{and} \ f(\mk{F}_l^{t+})= V_t \right \} $  for any $t\in (0,\pi/2)$.  We will show that we can ensure: 
	\begin{gather} \label{property} \forall t\in (0,\pi/2) \qquad \ \ 0\in U_t \ \text{and} \ f(\mk{F}_l^{t-})= U_t\setminus\{0\} \ \text{and} \   f(\mk{F}_l^{t+})= V_t. \end{gather}
	Indeed, if for some $t\in (0,\pi/2)$ this holds,   then  $f(\mk{F}_l^{t'-})= U_{t'}\setminus\{0\}$,  $ f(\mk{F}_l^{t'+})= V_{t'}$ for any other $t'\in (0,\pi/2)$.\footnote{Otherwise we would have  $f(\mk{F}_l^{t'-})= V_{t'}$,  $ f(\mk{F}_l^{t'+})= U_{t'}\setminus\{0\}$, which implies  $V_{t'}\subset U_t$ for $t'<t$, which is impossible or $U_t\setminus\{0\}\subset V_{t'}$ for $t<t'$. The latter in turn implies the contradiction  $U_{t'}\cap  V_{t'}\neq \emptyset$, since $(U_t\cap U_{t'})\setminus\{0\}\neq \emptyset$.} On the other hand by composing $\theta$ in \eqref{f} with $\frac{1}{z}$ we can ensure that the equalities in \eqref{property}  hold for some $t$ and  then they will hold for any $t\in (0,\pi/2)$. Next step is  to show:	
	\begin{gather}\label{convergences} \{z_i\}_{i\in \NN} \subset \mk{F}_l \ \text{and} \ \lim_{i\rightarrow \infty} \Im(z_i)=0 \qquad \Rightarrow \qquad    \lim_{i\rightarrow \infty} \abs{f(z_i)}=0. \end{gather}
We can assume that $\{z_i\}_{i\in \NN} \subset \mk{F}_l^{t-}$ for some $t\in (0,\pi/2)$ and therefore by \eqref{property} $\{f(z_i)\}_{i\in \NN} \subset U_t\setminus\{0\}$ is a bounded sequence. If  $\lim_{i\rightarrow \infty} \abs{f(z_i)}=0$ does not hold, then for some subsequence $\{z_{i_m}\}_{m\in \NN}$ holds  $\lim_{i\rightarrow \infty} f(z_{i_m})=
 q \in \CC^{\star}$. Now $\{f(z_{i_m})\}_{m\in \NN} \cup \{q\}$ is a compact subset of $\CC^{\star}$, but $f_{\vert \mk{F}_l}^{-1}(\{f(z_{i_m})\}_{m\in \NN} \cup \{q\}) \supset \{z_{i_m}\}_{m\in \NN}$ is not compact, since $\lim_{i\rightarrow \infty} \Im(z_{i_m})=0$ and $f_{\vert \mk{F}_l}^{-1}(\{f(z_{i_m})\}_{m\in \NN} \cup \{q\})$ does not contain any point with zero imaginary part, so we get a contradiction to (b) and proved \eqref{convergences}.

Since $l\geq 3$, then $\Delta_l=\log\left (\frac{a_l^2+1}{2  a_l} \right )>0$ (recall that $a_l=\frac{l+\sqrt{l^2-4}}{2}$ ). Therefore we can choose $0<\epsilon< \min\{\Delta_l,\pi/4\}$ and denote (see figure \eqref{3FDhyperbolic}):
\begin{gather} Q=\{u+\ri v \in \mc H: \max\{u,v\}<  \epsilon\}.  \end{gather} 
Since $Q\subset \mk{F}^o_l$  is an open subset and $Q\subset \mk{F}_l^{\frac{\pi}{4}-}$, by (a), (c) and \eqref{property}  it follows that $ \label{f vert} \bd Q & \rTo^{f_{\vert Q}} & f(Q) \ed$ is a biholomoprphism between open simply-connected  bounded subsets in $\CC$.   Obviously ${\rm Bd}_\CC(Q)={\rm Bd}_{\mc H}(Q)\amalg [-\epsilon, \epsilon] $ is   (image of) a closed Jordan curve. 
 On the other hand from \eqref{convergences} we see that, if exists, any continuous extension of $f_{\vert Q}$ to  ${\rm Cl}_{\CC}(Q)$ must map the entire real segment $[-\epsilon, \epsilon]$ to $\{0\}$ and thus it cannot be one-to-one. Therefore by  boundary correspondence theorem \cite[Theorem 2.24 on p. $70^3$]{Mark}   ${\rm Bd}_{\CC}(f(Q))$ cannot  be (image of) a closed Jordan  curve, otherwise  a one-to-one extension of  $f_{\vert Q}$ to a homeomorphism between ${\rm Cl}_\CC(Q)$ and  ${\rm Cl}_\CC(f(Q))$  must exist. We will  prove the theorem by deriving   the contrary:    ${\rm Bd}_{\CC}(f(Q))$ is an image  of a closed Jordan curve.

Since $f_{\vert \mk{F}_l}$ is proper, it follows that it  is a closed map (see e.g. \cite[p. 35]{Otto}), therefore $f(Q\amalg{\rm Bd}_{\mc H}(Q))$ is closed in $\CC^*$, and it follows that  
$f(Q)\cup f({\rm Bd}_{\mc H}(Q)) \cup\{0\}$ is closed in $\CC$. Hence  
${\rm Cl}_{\CC}(f(Q))\subset f(Q)\cup f({\rm Bd}_{\mc H}(Q))\cup \{0\}$. On the other hand \eqref{convergences} shows that $\{0\}\in{\rm Cl}_{\CC}(f(Q))$ and $f({\rm Bd}_{\mc H}(Q))\subset {\rm Cl}_{\CC}(f(Q))$ by continuity, hence
$ {\rm Cl}_{\CC}(f(Q)) = f(Q)\amalg f({\rm Bd}_{\mc H}(Q))\amalg \{0\}. $
Since $f(Q)$ is open in $\CC$, we obtain:
  \begin{gather} {\rm Bd}_{\CC}(f(Q))=f({\rm Bd}_{\mc H}(Q))\amalg \{0\}. \end{gather} Take any homeomorphism  $\bd & (0,2\pi) \rTo^{\kappa}& {\rm Bd}_{\mc H}(Q)\ed $, s. t. $\lim_{t\rightarrow 0}\kappa(t)=-\epsilon$,  $\lim_{t\rightarrow 2\pi}\kappa(t)=+\epsilon$, then $\bd (0,2\pi) & \rTo^{f\circ \kappa}& f({\rm Bd}_{\mc H}(Q))\ed$ is a homeomorphism and by \eqref{convergences} we have $\lim_{t\rightarrow 0} f(\kappa(t))=\lim_{t\rightarrow 2\pi} f(\kappa(t))=0$, then it is clear that the following  function \eqref{final function} is continuous and bijective, hence by compactness of $\SS^1$ it is homeomorphism and  the theorem is proved.
\begin{gather}\label{final function} \bd \SS^1 & \rTo^{\zeta} & f({\rm Bd}_{\mc H}(Q))\amalg \{0\}\ed \quad  \qquad\zeta(\exp(\ri t))=\left \{ \begin{array}{c c}  f(\kappa(t)) & t\in(0,2\pi) \\ 0 & t=0 \end{array}\right.\end{gather}

\begin{figure} \centering
        \begin{tikzpicture}
\path [fill=yellow] (-0.7,0)  --  (-0.7,0.7) --  (-0.7,0.7)  --  (0.7,0.7) --  (0.7,0.7)  --  (0.7,0);
\draw [<->] (-5.9,1.8)  to [out=350,in=135] (-2.3,0) ;
\draw [<->] (2.3,0)  to [out=45,in=190] (5.9,1.8);
\draw [dashed, ->] (0, 0) -- (0,3);
\draw [ -] (-4.7, 1.5) -- (4.7,1.5);
\node [left] at (0,3) {$v$};
\draw[<->] (-5.9,0)--(5.9,0);
\draw[dashed, <->] (-5.9,2)--(5.9,2);
\node [below] at (-0.7,0) {$-\epsilon$};
\node [below] at (0.7,0) {$\epsilon$};
\node [right] at (5.9,2) {$v=\frac{\pi}{2}$};
\node [below] at (-1.5,1.5) {$\im(\gamma_t)$};
\node [right] at (5.9,0) {$u$};
\node [below] at (0,0) {$0$};
\node [below] at (2.3,0) {$\Delta_l$};
\node [below] at (-2.3,0) {$-\Delta_l$};
\node [above] at (0,0) {$Q$};
\end{tikzpicture}
        \caption{$l\geq3$}
        \label{3FDhyperbolic}
 		     \end{figure} \epr

\subsection{The set of phases} \label{the set of phases}
Let $l\geq 2$ be an integer. Recall that  we denote \begin{gather} \label{al} a_l=\frac{l+\sqrt{l^2-4}}{2} \ \ \Rightarrow \ \  a_l^{-1}=\frac{l-\sqrt{l^2-4}}{2}; \quad  a_l^{-1}+a_l=\frac{a_l^2+1}{a_l}=l. \end{gather} 

In this section we write  for short $P_\sigma^l$ instead of $P_\sigma^{\mc T_l}$, and  will determine $P_\sigma^l$ (Proposition \ref{lemma for P_sigma}).

We start by some comments on the root system of $K(l)$.  The root system of $K(l)$ is $ \Delta_{l+}=\Delta_+(K(l)) = \Delta_{l+}^{re} \cup  \Delta_{l+}^{im}, $
where $\Delta_{l+}^{re }=\{(n,m)\in \NN^2  : n^2 + m^2 - l m n = 1 \}$ and  $\Delta_{l+}^{im }=\{(n,m)\in \NN^2  : n^2 + m^2 - l m n \leq  0 \}\setminus \{(0,0)\}$. It is well known that the real roots  $\Delta_{l+}^{re}$ are exactly the dimension vectors of the exceptional representations in $Rep_k(K(l))$ and for the imaginary roots  $\Delta_{l+}^{im}$ we have  formula (74) in \cite{DHKK}. 
 From \eqref{list of exceptional pairs} and Lemma \ref{nonvanishings} we have the complete list $\{s_{\leq 0}[1], s_{\geq 1}\}$ of exceptional representations in $Rep_k(K(l))$. Let us denote the corresponding dimension vectors as follows:
\be \label{nmi} (n_i,m_i)=\left \{ \begin{array}{c c} \ul{\dim}(s_i) & i\geq 1 \\  \ul{\dim}(s_i[1]) & i\leq 0 \end{array} \right .\ee
  Therefore     we can write: 
\begin{gather} \label{roots for K(l)}  \Delta_{l+}^{re} = \{(n_i,m_i): i\in \ZZ\}  \quad    \Delta_{l+}=\{(n_i,m_i): i\in \ZZ\} \cup \left \{a_l^{-1} \leq \frac{n}{m} \leq a_l : n\in \NN_{\geq 1}, m \in \NN_{\geq 1}\right \}. \end{gather} 

We will need later the following facts for the real roots $\{(n_i,m_i): i\in \ZZ\}$:
\begin{lemma} \label{lemma for nmi} {\rm (a)} $(n_{-1},m_{-1})=(l,1)$,  $(n_0,m_0)=(1,0)$, $(n_1,m_1)=(0,1)$,   $(n_2,m_2)=(1,l)$

 {\rm (b)} $(m_{-i},n_{-i})=(n_{i+1},m_{i+1})$ for $i\geq 0$.
 
{\rm (c)}  $n_{-i}>m_{-i}$  and $n_{i+1}<m_{i+1}$ for $i\geq 0$; $n_{i+1}>0$ and  $m_{-i}>0$ for $i\geq 1$.

{\rm (d)} $\frac{n_i}{m_i} = \frac{l}{2}-\sqrt{\frac{l^2}{4}-1+\frac{1}{m_i^2}}$   and  $\frac{n_{-i}}{m_{-i}} = \frac{l}{2}+\sqrt{\frac{l^2}{4}-1+\frac{1}{m_{-i}^2}}$ for $i\geq 1$. 

{\rm (e)} $\frac{n_{-1}}{m_{-1}}>\frac{n_{-2}}{m_{-2}}>...>\frac{n_{-i}}{m_{-i}}\stackrel{i\rightarrow\infty}{\rightarrow} a_l$ and  $0=\frac{n_{1}}{m_{1}}<\frac{n_{2}}{m_{2}}<...<\frac{n_{i}}{m_{i}}\stackrel{i\rightarrow\infty}{\rightarrow} a_l^{-1}$.
\end{lemma}
\bpr (a)  $(n_0,m_0)=(1,0)$, $(n_1,m_1)=(0,1)$ follow from the definition.  The triangles \eqref{triangles} for $i=1$,  and $i=0$ amount to short exact sequences: $\bd[1em] s_1^l&\rTo & s_2 & \rTo & s_0[1]\ed$ and $\bd[1em] s_1&\rTo & s_{-1}[1] & \rTo & s_0[1]^l\ed$ in $Rep_k(K(l))$, and it follows that $\ul{\dim}(s_2)= (1,l)$,  $\ul{\dim}(s_{-1}[1])= (l,1)$

(b) The equality for $0\leq i \leq 1$ follows from (a).  We make the induction assumption  that for some $p\geq 1$ the equality holds for any $0\leq i\leq p$, we will make the induction step, namely that  the equality  for $i=p+1$ follows from this induction assumption. Indeed,   for $i\geq 1$ from \eqref{triangles} we obtain the following short exact sequences in $Rep_k(K(l))$:  $\bd s_{-i-1}[1]&\rTo & s_{-i}^l[1] & \rTo & s_{-i+1}[1]\ed$, $\bd s_{i}&\rTo & s_{i+1}^l & \rTo & s_{i+2}\ed$ therefore for $i\geq 1$ we obtain:
\begin{gather} \label{recursive 1} n_{-i-1}= l n_{-i} - n_{-i+1} \qquad   m_{-i-1}= l m_{-i} - m_{-i+1} \\ 
\label{recursive 2}  n_{i+2}= l n_{i+1} - n_{i} \qquad   m_{i+2}= l m_{i+1} - m_{i}\end{gather} having these recursive formulas one easily caries out  the inductive step. 

(c) Due to (b) it is enough to show that $n_{i+1}<m_{i+1}$. For $i=0$ this is shown in (a).  For $i\geq 1$  we have $\hom(s_{i+1}, s_{-i}[1])> 0$, $\hom^1(s_{i+1}, s_{-i}[1])=0 $ (recall Lemma \ref{nonvanishings} ), hence the Euler formula amounts to:
\begin{gather} \scal{\ul{\dim}(s_{i+1}),\ul{\dim}(s_{-i}[1])}= \hom(s_{i+1}, s_{-i}[1])-\hom^1(s_{i+1}, s_{-i}[1])>0 \nonumber \\
\Rightarrow \scal{(n_{i+1},m_{i+1}), (n_{-i},m_{-i})}=n_{i+1} n_{-i} + m_{i+1} m_{-i} - l n_{i+1}m_{-i} >0. \nonumber   \end{gather}
Putting the equality from (b) in the last inequality we get $n_{i+1} m_{i+1} + m_{i+1} n_{i+1} - l n_{i+1}n_{i+1} >0$. Therefore $n_{i+1}(2 m_{i+1}  - l n_{i+1}) >0$, hence $n_{i+1}>0$, $ m_{i+1}   >\frac{l}{2} n_{i+1}\geq  n_{i+1}$ and (c) is proved.

(d) Take any $i\in \ZZ$, $i\neq 0$. From (c) we know that $m_i\neq 0$. From \eqref{roots for K(l)} we know that $n_i^2 +m_i^2 - l n_i m_i =1 $,  hence via the quadratic equation   $\left (\frac{n_i}{m_i} \right )^2- l \frac{n_i}{m_i}+1 -\frac{1}{m_i^2}  =0 $ we get $\frac{n_i}{m_i}=\frac{1}{2}\left ( l\pm \sqrt{l^2 - 4 + \frac{4}{m_i^2}} \right )$. One  checks that $\frac{1}{2}\left ( l+ \sqrt{l^2 - 4 + \frac{4}{m_i^2}} \right )>1$,  $\frac{1}{2}\left ( l- \sqrt{l^2 - 4 + \frac{4}{m_i^2}} \right )<1$ and  then from (c) we deduce (d).

(e) Using \eqref{recursive 1}, \eqref{recursive 2}, and induction one shows that $ m_{i}<m_{i+1}$, $m_{-i-1}>m_{-i}$ for $i\geq 0$ and then (e) follows from (d).  
\epr

Now we  have the necessary notations to describe $P_\sigma^l$: 

\begin{prop} \label{lemma for P_sigma} Let $\sigma \in \st(D^b(K(l)))$. 

{\rm (a)} If $\sigma \not \in \mc Z$  (see also Definition \ref{def of cal Z}), then  the  set of phases $P_{\sigma}^l$ is  finite (has up to 4 elements).

{\rm (b)} If $\sigma  \in \mc Z$, then    for any $j\in \ZZ$ we have the following formulas: 
 \begin{gather}\label{one inequality} 0<\phi_\sigma(s_{j+1})-\phi_\sigma(s_{j})<1 \\ \label{density P_sigma} \exp(\ri  \pi (1-\phi_{\sigma}(s_j)) )  \cdot P_\sigma^l =  \{ \pm 1\} \cup \{\pm\exp \left ( \ri f\left (n_i/m_i \right ) \right ):i\neq 0   \} \cup\pm\exp \left ( \ri f\left (\QQ \cap \left  [a_l^{-1},  a_l\right ] \right ) \right ),  \end{gather}   where  $x$, $y$, and   the  function (strictly increasing smooth) $f:[0,\infty)\rightarrow [\pi(\phi_\sigma(s_{j+1})-\phi_\sigma(s_{j})), \pi )$ are: 
\begin{gather} \label{f of tminuss} f(t)=\arccos\left( \frac{ x y-t}{\sqrt{t^2 + x^2 - 2 t  x y}} \right ) \quad x =\frac{\abs{Z(s_{j+1})}}{\abs{Z(s_j)}} \quad y= \cos\left (\pi  (\phi_\sigma(s_{j+1})-\phi_\sigma(s_{j})) \right ).\end{gather}
 {\rm (c)} For $\sigma \in \mc Z $ holds the equality $\{\pm \exp(\ri \pi \phi_\sigma(s_i))\}_{i\in \ZZ}= P_\sigma^l\setminus L(P_\sigma^l)$ (recall that by $L(P_\sigma^l)$ we denote the set of limit points in the circle of $P_\sigma^l$).

{\rm (d)} For any $\sigma \in \mc Z$ and any $j\in \ZZ$  hold:
\begin{gather} \label{e1}  \lim_{k\rightarrow +\infty} \pi \phi_\sigma(s_{k}[-1]) = u_\sigma\leq  v_\sigma =  \lim_{k\rightarrow -\infty} \pi \phi_\sigma(s_{k})   \\
    \pi (\phi_\sigma(s_j)-1)<\pi \phi_\sigma(s_{j+1}[-1])<\pi \phi_\sigma(s_{j+2}[-1])<\cdots<u_\sigma 
\leq  \nonumber \\[-2mm]   \label{e2}  \\[-2mm] \leq  v_\sigma <\cdots  <\pi \phi_\sigma (s_{j-2}) < \pi \phi_\sigma (s_{j-1}) <\pi \phi_\sigma(s_j)   \nonumber \\ \label{e3} 
\ol{P_\sigma^l} =\pm \exp \left (\ri \{\pi\phi_\sigma(s_{j+k}[-1])\}_{k\geq 1} \cup\ri [u_\sigma,v_\sigma]\cup \ri \{\pi \phi_\sigma(s_{j-k})\}_{k\geq 0} \right )\\
\label{e4} \frac{v_\sigma - u_\sigma}{u_\sigma-\pi \phi_\sigma(s_{j+1}[-1])}=\frac{f(a_l)-f(a_{l}^{-1})}{f(a_{l}^{-1})-\arccos(y)} \qquad \frac{v_\sigma - u_\sigma}{\pi \phi_\sigma(s_{j})-v_\sigma}=\frac{f(a_l)-f(a_{l}^{-1})}{\pi-f(a_l)}, 
 \end{gather}
where $f$, $x$, $y$ are as in \eqref{f of tminuss} and $ u_\sigma=f(a_{l}^{-1})+\pi (1-\phi_{\sigma}(s_j))$, $ v_\sigma=f(a_{l})+\pi (1-\phi_{\sigma}(s_j))$.

{\rm (e)} Let  $\sigma \in \mc Z $ and $0<\varepsilon<1$. Then ${\mathbb S}^1\setminus P_{\sigma}^l$ contains a closed $\varepsilon$-arc iff there exists $i\in \ZZ$ such that $\phi_\sigma(s_{i+1})-\phi_\sigma(s_{i})>\varepsilon$.
   \end{prop}
Before giving the proof of this proposition we make some preparatory  steps.

 For a pair of complex numbers $v=(z_1,z_2)$   we discussed in \cite{DHKK} (see  \cite[Lemma 3.18]{DHKK} and the first row of the proof) the following subset of the circle  \be R_{{v},\Delta_{l+}}=\left \{ \pm \frac{n z_1 + m z_2}{\abs{n z_1 + m z_2}} \vert (n,m) \in \Delta_{l+} \right \} \subset {\mathbb S}^1.\ee From \cite[Remark 3.16]{DHKK} and \eqref{roots for K(l)} we deduce that:

\begin{lemma} \label{lemma on R_v Deltal} For any pair of complex numbers $v=(z_1,z_2)$ of the form $ z_i= r_i \exp(\ri \phi_i)$, $r_i>0$, $i=1,2$, $0 < \phi_2 < \phi_1 \leq \pi$ holds:
\begin{gather} \label{lemma on R_v Deltal1}  R_{v,\Delta_{l+}}= \{ \pm \exp \left ( \ri \phi_1 \right )\} \cup \{\pm\exp \left ( \ri f\left (n_i/m_i \right ) \right ):i\neq 0   \} \cup \left \{\pm\exp \left ( \ri f\left (n/m \right ) \right ): n/m\in\left [a_l^{-1},  a_l\right ] \right \} \end{gather}
where $f:[0,\infty)\rightarrow [\phi_2, \phi_1 )\subset  (0, \pi )$ is the strictly increasing smooth function:
\begin{gather} \label{f of t} f(t)=\arccos\left( \frac{t r_1 \cos(\phi_1)+r_2 \cos(\phi_2)}{\sqrt{t^2 r_1^2 + r_2^2 + 2 t r_1 r_2 \cos(\phi_1-\phi_2)}} \right ), \ \ f(0)=\phi_2, \ \lim_{t \rightarrow \infty }f(t) = \phi_1.\end{gather}
\end{lemma}

	\textit{ Proof of Proposition \ref{lemma for P_sigma}.}  Let $\sigma = (Z, \mc P) \in \st(D^b(K(l)))$. From Lemma \ref{sections} and equality \eqref{st} we have either $\sigma \in \Theta_{(s_j,s_{j+1})}\setminus \mathcal Z$ for some $j\in \ZZ$ or  $\sigma \in  \mathcal Z$ (see also Definition \ref{def of cal Z}).
	
 (a)	Assume first that $\sigma \in \Theta_{(s_j,s_{j+1})}\setminus \mathcal Z$ for some $j\in \ZZ$.  Then  by  Lemmas \ref{lemma for f_E(Theta_E)} and \ref{sections}   we see that \be \label{case in a lemma} s_j, s_{j+1} \in \sigma^{ss}  \quad  \phi( s_j)+1\leq   \phi( s_{j+1}). \ee We will show that in this case 
	 $P_\sigma^l=\{\pm \exp(\ri \pi \phi_\sigma(s_j)), \pm \exp(\ri \pi \phi_\sigma(s_{j+1}))\}$.
	Indeed, \eqref{case in a lemma} implies that there exists $k\geq 1$ such that $\phi(s_j)\leq \phi(s_{j+1}[-k])< \phi(s_j)+1$. From Lemma \ref{nonvanishings} it follows that $(s_j,s_{j+1}[-k])$ is a $\sigma$-exceptional pair (as defined in \cite[Definition 3.17]{DK1}). From \cite[Corollary 3.18]{DK1} (and its proof) it follows that the extension closure of 
 $(s_j,s_{j+1}[-k])$	equals $\mc P(t,t+1]$ for some $t\in \RR$. Since $(s_j, s_{j+1}[-k])$ is an exceptional pair, each element $Y$ in the extension closure of $(s_j, s_{j+1}[-k])$ can be put in a triangle of the form $\bd  s_{j+1}[-k]^a & \rTo & Y& \rTo & s_j^b & \rTo & s_{j+1}[-k+1]^b \ed$ for some $a,b\in \NN$. Take any $X\in \sigma^{ss}$, then  for some $i\in \ZZ$ we have $\phi_\sigma(X[i]) \in (t,t+1]$ and therefore we have a triangle:
\begin{gather} \label{triangle11}  \bd s_{j+1}[-k]^a & \rTo^\alpha & X[i] & \rTo^\beta & s_j^b & \rTo & s_{j+1}[-k+1]^a.  \ed  \end{gather}

If $a=0$ or $b=0$, then $X[i]\cong s_j^b $ or $X[i]\cong s_{j+1}[-k]^a $ and hence  $\phi_\sigma(X[i])=\phi(s_j)$ or $\phi_\sigma(X[i])=\phi(s_{j+1}[-k])$ and the $\exp(\ri \pi \phi_\sigma(X)) \in \{\pm \exp(\ri \pi \phi_\sigma(s_j)), \pm \exp(\ri \pi \phi_\sigma(s_{j+1}))\}$.

Next assume that  $a\neq 0$ and  $b\neq 0$. If $\phi (s_{j})=\phi (s_{j+1}[-k])$,  then we obtain   $\exp(\ri \pi \phi_\sigma(X)) =  \exp(\ri \pi \phi_\sigma(s_j))$using \eqref{triangle11}. Thus, we reduce to $\phi (s_{j})<\phi (s_{j+1}[-k])$, which in turn by \eqref{triangle11}, $X[i]\in \sigma^{ss}$, and \eqref{vanishing formula}   implies that $\alpha=0$  or $\beta=0$. If $\alpha =0$, then $s_j^b \cong X[i] \oplus s_{j+1}[-k+1]^b$ and by \cite[Lemma 3.7]{DK1} it follows that $\phi(s_j) = \phi(X[i])=\phi(s_{j+1}[-k+1])$; if $\beta =0$, then $s_{j+1}[-k]^a \cong X[i] \oplus s_{j}[-1]^b$ and by \cite[Lemma 3.7]{DK1} it follows that $\phi(s_{j+1}[-k]) = \phi(X[i])=\phi(s_{j}[-1])$. Thus we see that  \eqref{case in a lemma} implies  $P_\sigma^l=\{\pm \exp(\ri \pi \phi_\sigma(s_j)), \pm \exp(\ri \pi \phi_\sigma(s_{j+1}))\}$,  and (a) is proved.

(b) If  $\sigma \in \mathcal Z$,   then  Lemma \ref{sections} shows that  for any $j\in \ZZ$ holds  $ s_j, s_{j+1} \in \sigma^{ss} $,  $ \phi_\sigma( s_j)< \phi_\sigma( s_{j+1}) < \phi_\sigma( s_j) +1  . $ Choosing one $j\in \ZZ$, denoting $\sigma'=(Z',\mc P')= (-\log\abs{Z(s_j)}+ \ri \pi (1-\phi_{\sigma}(s_j)))\star \sigma$ and  using \eqref{star properties 1}, \eqref{star properties 2}, we get: 
\begin{gather}\label{another case in a lemma1} Z'(s_j)=-1, \ \ \abs{Z'(s_{j+1})}=\frac{\abs{Z(s_{j+1})}}{\abs{Z(s_{j})}}, \ \ \phi_{\sigma'}(s_{j+1})=\phi_{\sigma}(s_{j+1})+1-\phi_{\sigma}(s_{j})   \\  \label{another case in a lemma}1=\phi_{\sigma'}( s_j)<   \phi_{\sigma'}(s_{j+1}) < \phi_{\sigma'}( s_j) +1=2 \  \Rightarrow \  0< \phi_{\sigma'}( s_{j+1}[-1]) < \phi_{\sigma'}( s_j) =1 . \end{gather}
Let $\mc A$ be the  extension closure of $(s_j,s_{j+1}[-1])$. 
Utilizing Lemma \ref{nonvanishings} and recalling that $\hom(s_j,s_{j+1})=l\geq 2$ (see for example the arguments before \eqref{triangles}) we see that $(s_j,s_{j+1}[-1])$ is an $l$-Kronecker pair \cite[Definition 3.20]{Dimitrov},  and by \eqref{another case in a lemma} it is a $\sigma'$-exceptional pair as well. From \cite[Corollary 3.18]{DK1} (and its proof) we see that the extension closure $\mc A$  of  $(s_j,s_{j+1}[-1])$ coincides with $\mc P'(0,1]$. Applying \cite[Lemma 3.19]{Dimitrov}  to $(s_j,s_{j+1}[-1])$ we see that $\mc A$ is the heart of a bounded t-structure  of $\mc T_l$ and due to the equality  $\mc A=\mc P'(0,1]$ we have actually  $\sigma' \in \HH^{\mc A}$ (see \cite[Definition 2.28]{Dimitrov}). Now all the conditions of \cite[Corollary 3.21]{Dimitrov} with the exceptional pair $(s_j,s_{j+1}[-1])$ hold and we deduce that $P_{\sigma'}^l=R_{v,\Delta_{l+}}$, where $v=(Z'(s_j),Z'(s_{j+1}[-1]))$. On the other hand  \eqref{rotating P_sigma} shows that $\exp(\ri  \pi (1-\phi_{\sigma}(s_j)) ) \cdot P_{\sigma}^l=R_{v,\Delta_l}$

To determine  the set $R_{v,\Delta_l}$ we use  Lemma \ref{lemma on R_v Deltal} and observe  that now (see \eqref{another case in a lemma1}) \\ $v=\left (-1,\frac{\abs{Z(s_{j+1})}}{\abs{Z(s_{j})}}  \exp(\ri \pi (\phi_{\sigma}(s_{j+1})-\phi_{\sigma}(s_{j})) \right )$, $0 < \pi  (\phi_{\sigma}(s_{j+1})-\phi_{\sigma}(s_{j})) <  \pi$,  in particular the equality \eqref{lemma on R_v Deltal1}  yields \eqref{density P_sigma} and the function \eqref{f of t} has the form \eqref{f of tminuss}.

(c) Let $\sigma \in \mc Z$. In (b) $j$ was any integer, here we  choose  $j=0$.  Now formulas  \eqref{another case in a lemma1} and \eqref{another case in a lemma} give:   \begin{gather} Z'(s_0)=-1 \qquad Z'(s_1[-1])=\abs{Z'(s_1)} \exp(\ri \pi \phi_{\sigma'}(s_1[-1])),   \ \ 0< \phi_{\sigma'}( s_{1}[-1]) < \phi_{\sigma'}( s_0)=1. \end{gather}
  Since $s_0[1]$, $s_1$ are the simple representations and since $s_{\geq 1}, s_{\leq 0}[1] \in Rep_k(K(l))$ (Lemma \ref{nonvanishings}), it follows that (see also \eqref{nmi}) for any $i\geq 1$ $Z'(s_i)=n_i Z'(s_0[1]) + m_i  Z'(s_1)$, and  for any $i\leq 0$ $Z'(s_i[1])=n_i Z'(s_0[1]) + m_i  Z'(s_1)$, and now using   \cite[Remark 3.16]{DHKK} (in particular  $f$ is as in \eqref{f of tminuss}) we obtain :
\begin{gather}\label{formula for the accumulation 1} \pm \frac{Z'(s_i)}{\abs{Z'(s_i)}}= \left \{ \begin{array}{c c} \mp \frac{n_i Z'(s_0) + m_i  Z'(s_1[-1])}{\abs{n_i Z'(s_0) + m_i  Z'(s_1[-1])}} = \mp \exp(\ri f(n_i/m_i)) & i\geq 1\\
\mp 1 & i=0  \\ \pm \frac{n_i Z'(s_0) + m_i  Z'(s_1[-1])}{\abs{n_i Z'(s_0) + m_i  Z'(s_1[-1])}} = \pm \exp(\ri f(n_i/m_i)) & i \leq -1 \end{array} \right. .
\end{gather}
 In (b) we showed that $P_{\sigma'}^l$ equals the set on the RHS of \eqref{density P_sigma}. Due to  Lemma \ref{lemma for nmi} we get that $L(P_{\sigma'}^l)=\pm\exp \left ( \ri f\left ( \left  [a_l^{-1},  a_l\right ] \right ) \right )$, and therefore 
    \eqref{formula for the accumulation 1} and \eqref{phase formula} imply  that $P_{\sigma'}^l\setminus L(P_{\sigma'}^l)=\{\pm \exp(\ri \pi \phi_{\sigma'}(s_i))\}_{i\in \ZZ}$. Recalling that $\sigma' = \lambda \star\sigma  $ for certain $\lambda \in \CC$ with the help of formulas \eqref{star properties 1} and \eqref{rotating P_sigma} we deduce the desired $P_{\sigma}^l\setminus L(P_{\sigma}^l)=\{\pm \exp(\ri \pi \phi_{\sigma}(s_i))\}_{i\in \ZZ}$.

(d) In (b) we showed that $P_{\sigma'}^l$ for $\sigma'=(Z',\mc P')= (-\log\abs{Z(s_j)}+ \ri \pi (1-\phi_{\sigma}(s_j)))\star \sigma$ equals the RHS of \eqref{density P_sigma} and taking into account Lemma \ref{lemma for nmi}  we deduce that $P_{\sigma'}^l\setminus L(P_{\sigma'}^l)=\pm 1 \cup \{\pm\exp \left ( \ri f\left (n_i/m_i \right ) \right ):i\neq 0 \}$, which combined with (c) yields:
\begin{gather}\label{equality for e1}  \pm 1 \cup \{\pm\exp \left ( \ri f\left (n_i/m_i \right ) \right ):i\neq 0 \} = \{\pm \exp(\ri \pi \phi_{\sigma'}(s_i))\}_{i\in \ZZ}.  \end{gather} 
Recalling that \eqref{one inequality}  holds for any $j\in \ZZ$ and also \eqref{nonvanish12}, \eqref{another case in a lemma} we derive:
\begin{gather}\label{equality for e2}  0=\phi_{\sigma'}(s_j)-1< \phi_{\sigma'}(s_{j+1}[-1])<\phi_{\sigma'}(s_{j+2}[-1])<\cdots  < \phi_{\sigma'} (s_{j-2}) <  \phi_{\sigma'} (s_{j-1}) < \phi_{\sigma'}(s_j) =1 \end{gather}
We already know that (see (b) of the proposition and (e) in Lemma \ref{lemma for nmi}) 
 \begin{gather}\label{f(0)<dots} f(0)=f(n_1/m_1)=\pi (\phi_{\sigma}(s_{j+1})-\phi_{\sigma}(s_{j}))=\pi \phi_{\sigma'}(s_{j+1}[-1])=\arccos(y). \end{gather} Furthermore from (e) in Lemma \ref{lemma for nmi} we deduce:
	\begin{gather}\label{equality for e3}  0<f\left (\frac{n_1}{m_1}\right )<f\left (\frac{n_2}{m_2}\right )<\cdots<f(a_l^{-1})\leq f(a_l)<\cdots<f\left(\frac{n_{-2}}{m_{-2}}\right)<f\left(\frac{n_{-1}}{m_{-1}}\right)<\pi. \end{gather}
	By induction the equalities \eqref{equality for e1},  \eqref{equality for e2}, \eqref{equality for e3} imply: 
	\begin{gather}\label{k geq 1 Rightarrow} k\geq 1  \ \ \Rightarrow \ \ f\left (\frac{n_k}{m_k}\right )=\pi\phi_{\sigma'}(s_{j+k}[-1]) \quad f\left (\frac{n_{-k}}{m_{-k}}\right )=\pi\phi_{\sigma'}(s_{j-k}). \end{gather} 
	Now recalling that (see \eqref{star properties 2})  \begin{gather} \label{phiprim and phi} \forall i \in \ZZ \ \ \  \phi_{\sigma'}(s_{i}) = \phi_{\sigma}(s_{i})+1-\phi_\sigma(s_{j})\end{gather}	
we deduce \eqref{e1}, \eqref{e2}, \eqref{e4} from \eqref{k geq 1 Rightarrow}, \eqref{f(0)<dots}, and \eqref{equality for e2}. The equality \eqref{e3} in turn follows from \eqref{density P_sigma}, \eqref{k geq 1 Rightarrow}, \eqref{phiprim and phi}. 
		
(e) Follows easily from the already proven (d).
		
	 \subsection{Computing the norms of Kornecker quivers, i. e.  \texorpdfstring{$\norm{D^b(K(l))}_\varepsilon$}{\space}} \label{norms on kroneckers}. 

If we define the function: 
\begin{gather} \label{Function} F: (0,+\infty)\times (-1,+1)\times  (0,+\infty) \rightarrow (0,\pi) \qquad F(x,y,t)=\arccos\left( \frac{ x y-t}{\sqrt{t^2 + x^2 - 2 t  x y}} \right )\end{gather}
 then using Proposition  \ref{lemma for P_sigma}  (a), Lemma \ref{lemma for nmi} (e), and   formulas  \eqref{density P_sigma},  \eqref{f of tminuss} one concludes that:  
\begin{prop} \label{formula for the volume}
Let $\sigma = (Z, \mc P)\in \st(\mc T_l)$.

 If $\sigma \not \in \mc Z$, then $\vol\left (\ol{P_\sigma^l}\right )=0$.  If $\sigma \in \mc Z$, then for any $j\in \ZZ$    holds: \begin{gather}   \frac{1}{2} \vol\left (\ol{P_\sigma^l}\right ) =  F\left (x_j(\sigma),y_j(\sigma),a_l\right )- F\left (x_j(\sigma),y_j(\sigma),a_l^{-1}\right ), \nonumber \\[-2mm] \label{closure length} \\[-2mm] \mbox{where} \qquad  x_j(\sigma)=\frac{\abs{Z(s_{j+1})}}{\abs{Z(s_{j})}} \qquad y_j(\sigma)=\cos\left (\pi \left (\phi_\sigma(s_{j+1})-\phi_\sigma(s_{j})  \right )\right ). \nonumber \end{gather}
\end{prop}
One computes \begin{gather} \label{partial xt F} \frac{\partial}{\partial x} F(x,y,t)=\frac{-t \sqrt{1-y^2}}{t^2+x^2- 2 t x y} \qquad \frac{\partial}{\partial t} F(x,y,t)=\frac{x \sqrt{1-y^2}}{t^2+x^2- 2 t x y} 
\end{gather}
 and therefore:  \begin{gather} \frac{\partial}{\partial x}\left ( F(x,y,a_l)-F(x,y,a_l^{-1}) \right )=\frac{a_l (a_l^2-1) \sqrt{1-y^2}}{(1+(a_l x)^2- 2 a_l x y)(a_l^2+x^2- 2 a_l x y)}(1-x^2), \end{gather} which implies that for any $x>0$, $y\in (-1,+1)$ we have:
\begin{gather} F(x,y,a_l)-F(x,y,a_l^{-1})\leq  F(1,y,a_l)-F(1,y,a_l^{-1})   \label{going to x=1}  \end{gather}
On the other hand one computes that for any $y\in (-1,+1), t \in (0,+\infty)$ holds:
\begin{gather} F(1,y,a_l)-F(1,y,a_l^{-1}) = \arccos\left( \frac{  y-a_l}{\sqrt{a_l^2 + 1 - 2 a_l   y}} \right )-\arccos\left( \frac{ a_l y-1}{\sqrt{a_l^2 + 1 - 2 a_l   y}} \right )  \\ \label{derivative} \frac{\partial}{\partial y}\left ( F(1,y,a_l)-F(1,y,a_l^{-1}) \right )=\left \{ \begin{array}{c c} \frac{a_l^2-1}{\sqrt{1-y^2} (a_l^2+1- 2 a_l y)}>0 & l\geq 3  \\ 0 & l=2 \end{array} \right . \\
  \frac{\partial}{\partial t}\left ( F(1,y,t)-F(1,y,t^{-1}) \right )= \frac{2 \sqrt{1-y^2}}{t^2 +1 - 2 t y}>0.
\end{gather}
Therefore the numbers \eqref{Ky(l)} depending on $\varepsilon \in (0,1)$ and $l\geq 2$  satisfy \eqref{props of Ky(l) 0}, \eqref{props of Ky(l) 1}, \eqref{props of Ky(l) 3}: 
\begin{gather} \label{Ky(l)} K_\varepsilon(l)= \arccos\left( \frac{  \cos(\pi \varepsilon) -a_l}{\sqrt{a_l^2 + 1 - 2 a_l   \cos(\pi \varepsilon)}} \right )-\arccos\left( \frac{ a_l \cos(\pi \varepsilon)-1}{\sqrt{a_l^2 + 1 - 2 a_l   \cos(\pi \varepsilon)}} \right )   \\
\label{props of Ky(l) 0}   0<\varepsilon<1  \qquad \Rightarrow \qquad  K_\varepsilon (2)=0 \\ 
 \label{props of Ky(l) 1} l\in \NN_{\geq 3} \ \   0<u<v<+1  \qquad \Rightarrow \qquad  K_u(l)> K_v(l) \\
\label{props of Ky(l) 3} 0<\varepsilon<1  \ \ 2\leq l_1 < l_2 \in \NN_{\geq 2} \qquad \Rightarrow \qquad  K_\varepsilon (l_1)< K_\varepsilon (l_2)\\
\label{props of Ky(l) 4} \lim_{l\rightarrow +\infty}  K_\varepsilon (l)=\pi(1-\varepsilon). \end{gather}
The inequality \eqref{going to x=1} and the  derivative \eqref{derivative}    imply that for $\varepsilon \in (0,+1)$ and $ l\geq 2$ holds:
\begin{gather} \label{the supremum 1} \sup_{(x,y)\in (0,+\infty)\times (-1,\cos(\pi \varepsilon))}\left \{ F(x,y,a_l)- F(x,y,a_l^{-1})\right \} =  K_\varepsilon(l), 
\end{gather} Note that $\sup_{(x,y)\in (0,+\infty)\times (-1,1)}\left \{ F(x,y,a_l)- F(x,y,a_l^{-1})\right \}$ is always equal to $\pi$ independently on $l\geq 3$ as opposed to $K_\varepsilon(l)$, which is strictly increasing on $l$.

Finally we note that for $\varepsilon=1/2$  the expression \eqref{Ky(l)} takes a simple form (recall that $l=\frac{a_l^2+1}{a_l}$):\footnote{one   shows this using the equality $\arccos(x)-\arccos\left (\sqrt{1-x^2}\right )=\arccos\left (2 x \sqrt{1-x^2}\right )$, which holds for $0\leq x \leq \frac{1}{\sqrt{2}}$} 
\be \label{props of Ky(l) 2} 
K_{\frac{1}{2}}(l)=\arccos\left (\frac{2}{l}\right ).
\ee

Now we can compute  $ \norm{D^b(K(l))}_{\varepsilon}$.
\begin{prop} \label{norm of Kroneckers} Let $\varepsilon \in (0,1)$,  $l\geq 2$,  and let  $ K_{\varepsilon}(l)$ be as  in \eqref{Ky(l)}. Then 
$ \norm{D^b(K(l))}_{\varepsilon} = K_{\varepsilon}(l)$. 
\end{prop}
\bpr  From Proposition \ref{lemma for P_sigma} (a), (b), (e)  we see that  $P_\sigma$ is not dense in $\mathbb S^1$ for all $\sigma$, and   \eqref{main def eq} reduces to the following formula:
\begin{gather} \norm{D^b(K(l))}_{\varepsilon} =   \sup \left  \{ \frac{ \vol\left (\ol{P_\sigma^l}\right )}{2}: \sigma \in  \mc{Z} \ \mbox{and   there exists} \  j\in \ZZ \ \mbox{ such that}  \  \phi_\sigma(s_{j+1})-\phi_\sigma(s_{j})>\varepsilon  \right  \}\nonumber \\
= \sup \left  \{ \frac{ \vol\left (\ol{P_\sigma^l}\right )}{2}: \sigma \in \bigcup_{j\in \ZZ} \left \{ x \in \mc{Z} :   \phi_x(s_{j+1})-\phi_x(s_{j})>\varepsilon  \right \} \right  \}\nonumber \\
=\sup \left \{ \sup\left \{\frac{ \vol\left (\ol{P_\sigma^l}\right )}{2}:\sigma \in \mc{Z}  \ \mbox{and} \   \phi_\sigma(s_{j+1})-\phi_\sigma(s_{j})>\varepsilon \right \} : j \in \ZZ \right \}.\nonumber  \end{gather}
 By using \eqref{closure length} and \eqref{the supremum 1} we will show that for all $j \in \ZZ$ holds: 
 \be\label{step to prove the norm of Kronecker} \sup\left \{\frac{ \vol\left (\ol{P_\sigma^l}\right )}{2}:\sigma \in \mc{Z}  \ \mbox{and} \   \phi_\sigma(s_{j+1})-\phi_\sigma(s_{j})>\varepsilon \right \} = K_{\varepsilon}(l) \ee and then the proposition follows.   Recalling Lemma \ref{sections} we see that $\sigma \in \mc{Z}$ and  $\phi_\sigma(s_{j+1})-\phi_\sigma(s_{j})>\varepsilon $ iff $s_{j}, s_{j+1} \in \sigma^{ss}$ and $\varepsilon <\phi_\sigma(s_{j+1})-\phi_\sigma(s_{j})< 1$, furthermore restricting the map  \eqref{mapping} to the set of stability conditions $\sigma$ with $s_{j}, s_{j+1} \in \sigma^{ss}$ and $\varepsilon <\phi_\sigma(s_{j+1})-\phi_\sigma(s_{j})< 1$ we see that the set of pairs $(x_j(\sigma)$, $y_j(\sigma))$ from \eqref{closure length} for these $\sigma$ is:
\begin{gather} \left \{(x_j(\sigma), y_j(\sigma) ):\sigma \in \mc{Z}  \ \mbox{and} \   \phi_\sigma(s_{j+1})-\phi_\sigma(s_{j})>\varepsilon\right \} = (0,+\infty)\times (-1,\cos(\pi \varepsilon)).\end{gather}    Combining the latter equality with \eqref{closure length} and \eqref{the supremum 1}  leads to \eqref{step to prove the norm of Kronecker}. \epr
\begin{coro} Let $\varepsilon \in (0,1)$. Let $2\leq l_1 $ and   $ 2 \leq l_2$. Then:
  \begin{gather} \norm{D^b(K(l_1))}_\varepsilon < \norm{D^b(K(l_2))}_\varepsilon  \iff  l_1<l_2 \end{gather}
\end{coro}
\bpr Follows from the previous proposition and  \eqref{props of Ky(l) 3}. \epr

\begin{coro} Let $\varepsilon \in (0,1)$ and  let $l\in \NN_{\geq 1} $.

 Then $\norm{D^b(K(l))}_\varepsilon = 0$ iff $\st(D^b(K(l)))$ is affine (biholomorphic to $\CC^2$).

\end{coro}
\bpr Propositions \ref{norm of Dynkin} and   \ref{norm of Kroneckers} imply that $\norm{D^b(K(l))}_\varepsilon = 0$ iff $l\leq 2$, and table \eqref{table_intro1} shows that $\st(D^b(K(l)))$ is biholomorphis to $\CC^2$ iff $l\leq 2$. 
\epr

\section{ The inequality \texorpdfstring{$\norm{\scal{E_1,E_2}}_{\varepsilon}\geq K_\varepsilon(\hom^{min}(E_1,E_2))$}{\space} } \label{criteria}

 In this section we dervie a formula, which will help us to compute other norms. To that end it is useful to extend the definition of $K_\varepsilon(l)$ in \eqref{Ky(l)} by postulating  $K_\varepsilon(0)=K_\varepsilon(1)=0$.  Recall that the notation $\hom^{min}(E_1,E_2)$ is explained in \eqref{hom min}. 
\begin{prop} \label{inequality for any exc pair} Let $\mc T$ be a proper category, and let $(E_1,E_2)$ be any exceptional pair in it. Then 
\begin{gather} \norm{\scal{E_1,E_2}}_{\varepsilon}\geq K_\varepsilon\left (\hom^{min}(E_1,E_2)\right ) \ \ \mbox{for} \ \ \varepsilon \in (0,1). \end{gather}
\end{prop}
\bpr We can assume that $\hom^{\leq 0}(E_1,E_2)=0$ and  $l=\hom^1(E_1,E_2)\neq 0$, and under these assumption we have to show that 
\begin{gather} \norm{\scal{E_1,E_2}}_{\varepsilon}\geq K_\varepsilon\left (l\right ). \end{gather}
 Let $\mc D$ be the triangulated  subcategory $\scal{E_1,E_2}$. The assumptions on $(E_1,E_2)$ are the same  as in the definition of an $l$-Kronecker pair,  \cite[Definition 3.20]{Dimitrov}, and we can apply  \cite[Lemma 3.19, Corollary 3.21]{Dimitrov} to it. In particular the extension closure $\mc A$ of $(E_1,E_2)$ is a heart of a bounded t-structure in $\mc D$ with simple objects $E_1,E_2$, and any stability condition $\sigma = (Z, \mc P) \in \HH^{\mc A}\subset \st(\mc D)$ with $\arg(Z(E_1))>\arg(Z(E_2))$ satisfies $P_\sigma^{\mc D}=R_{v,\Delta_{l+}}$, where $v=(Z(E_1),Z(E_2))$. The arguments in the beginning of the proof of Lemma \ref{criteria for phase gap}  show that for each $v\in \HH^2$ there exists unique  $\sigma=(Z, \mc P) \in \HH^{\mc A}$ with $v=(Z(E_1),Z(E_2))$ and that $\sigma$ is  full. For any  $0<\mu$ such that $\mu+\varepsilon <1$ choose the vector $(-1,\exp(\ri \pi (\varepsilon +\mu)))=v_{\mu}$ and    denote by  $\sigma_{\mu}$ the stability condition  $\sigma_{\mu}=(\mc P_{\mu}, Z_{\mu}) \in \HH^{\mc A}$ with  $( Z_{\mu}(E_1), Z_{\mu}(E_2))=v_{\mu}$. The given arguments ensure that $\sigma_{\mu}$ is full and  $P_{\sigma_{\mu}}^{\mc D}= R_{v_{\mu},\Delta_{l+}}$. Using the formula for $R_{v_{\mu},\Delta_{l+}}$  in Lemma \ref{lemma on R_v Deltal}   for the given $v_\mu$ one  derives:
\be \frac{\vol\left (\ol{P_{\sigma_{\mu}}^{\mc D}} \right )}{2}=\frac{\vol\left (\ol{R_{v_{\mu},\Delta_{l+}}} \right )}{2}=K_{\varepsilon + \mu}(l) ,\ee
where  $K_{\varepsilon + \mu}(l) $ is in \eqref{Ky(l)}. Note that the arc $\exp\left (\ri \pi [\mu/2,\varepsilon+\mu/2]\right )$ is in the complement of $P^{\mc D}_{\sigma_\mu}$ and therefore $\sigma_\mu \in \st_{\varepsilon}(\mc D)$.    Now from the very  Definition \ref{main def} we see that $\norm{\mc D}_{\varepsilon}\geq K_{\varepsilon + \mu}(l)$ for any small enough positive $\mu$, letting $\mu \rightarrow 0$ we derive the desired $\norm{\mc D}_{\varepsilon}\geq K_{\varepsilon}(l)$.
\epr

\begin{coro} \label{inequality for sequences} Let $\mc E = (E_0, E_1,\dots,E_n)$ be an exceptional collection in a proper triangulated category $\mc T$. Then for any $0\leq i<j\leq n$ we have $\norm{\scal{\mc E}}_{\varepsilon} \geq K_\varepsilon \left (\hom^{min}(E_i,E_j) \right )$. 
\end{coro}
\bpr  Take $0\leq i < j \leq n$. 
 By mutating the sequence $\mc E$ (see Remark \ref{mutations}) one can get a sequence $\mc E'$ of the form $\mc E'=(E_i,E_j,C_2,\dots,C_n)$ such that $\scal{\mc E}=\scal{\mc E'}$. Corollary \ref{coro for esc coll} implies $\norm{\scal{\mc E}}_{\varepsilon}=\norm{\scal{\mc E'}}_{\varepsilon}\geq \norm{\scal{E_i,E_j}}_{\varepsilon}$, and due to Proposition \ref{inequality for any exc pair}  we get $\norm{\scal{E_i,E_j}}_{\varepsilon}\geq K_\varepsilon(\hom^{min}(E_i,E_j))$.
\epr
\begin{coro} \label{from exc collections to maximal norm} Let $\mc T$ be a proper triangulated category  such that for each $l\in \NN$ there exists a full exceptional collection $(E_0,E_1,\dots, E_n)$ and integers  $0\leq i<j\leq n$ for which $\hom^{min}(E_i,E_j)\geq l$. Then $\norm{\mc T}_{\varepsilon}= \pi (1-\varepsilon)$ for any $\varepsilon \in (0,1)$. 
\end{coro}
\bpr   The given property of $\mc T$ combined with   Corollary \ref{inequality for sequences} amounts to   $\norm{\mc T}_{\varepsilon}\geq K_{\varepsilon}(l)$ for each $l\geq \NN$ (Recall also \eqref{props of Ky(l) 3}). Now from \eqref{props of Ky(l) 4} and Remark \ref{remark leq pi(1-varepsilon)}  we obtain $\norm{\mc T}_{\varepsilon}=\pi (1-\varepsilon)$.
\epr
\begin{coro} \label{some limitations} Let $\mc T$ be a proper category, and let $0<\varepsilon <1$. 

{\rm (a)} If $\norm{\mc T}_{\varepsilon}=0$, then for any full exceptional collection $\mc E=(E_0,E_1,\dots, E_n)$ and for  any $0\leq i<j\leq n$ we have $\hom^{min}(E_i,E_j)\leq 2$.

{\rm (b)} If $\norm{\mc T}_{\varepsilon}\leq K_\varepsilon (l)$, $l\geq 2$, then for any full exceptional collection $\mc E=(E_0,E_1,\dots, E_n)$ and for  any $0\leq i<j\leq n$ we have $\hom^{min}(E_i,E_j)\leq l$.

{\rm (c)}  If  $\norm{\mc T}_{\varepsilon}<\pi (1-\varepsilon)$, then there exists $l\in \NN$ such that for any full exceptional collection $\mc E=(E_0,E_1,\dots, E_n)$ and for  any $0\leq i<j\leq n$ we have $\hom^{min}(E_i,E_j)\leq l$.
\end{coro}
We will aply Corollary \ref{from exc collections to maximal norm} to various examples. More precisely we will show that 
\begin{prop} \label{prop maximal norms}
	In the following examples of triangulated categories are satisfied the  conditions of Corollary \ref{from exc collections to maximal norm}. In particular  $\norm{\mc T}_{\varepsilon}=(1-\varepsilon) \pi$ for any $\mc T$ in this list of examples.
	
	{\rm (a)}  $D^b(Q)$, where  $Q$ is an acyclic quiver, s.t. there exists a subset $A\subset V(Q)$ such that the quiver $Q_A$ is affine and there exists a vertex $v \in V(Q)$ such that $v$ is a source or a sink in $Q_{A\cup\{v\}}$ (see Definition \ref{adjacent} for the terminology)
	
	{\rm (b)} $D^b(\PP^n), n\geq 2$; 	
	{\rm (c) }  $D^b(\PP^1\times \PP^1)$;	
	{\rm (d)\footnote{where  $\mathbb F_m$ is the $m$-th Hirzebruch surface}} $D^b({\mathbb F}_m), m\geq 0$
	
	{\rm (e)} $D^b(X)$, where $X$ is a smooth algebraic variety obtained from $\PP^n$, $n\geq 2$, or from $\PP^1\times \PP^1$, or from ${\mathbb F}_m), m\geq 0$ by a sequence of blow ups at finite  finite  number  of points;
	
	{\rm (f )}$D^b(S)$, where $S$ is any smooth complete rational surface\footnote{in particular any somoth projective surface.}
\end{prop}
\begin{df} \label{adjacent} For any quiver $Q$ and any subset $A\subset V(Q)$ we denote by  $Q_A$ the quiver whose  vertices are  $A$ and whose  arrows  are those arrows of $Q$ whose initial and final vertex is in $A$.  A vertex $v \in V(Q)$ is called adjacent to $A$ if there exists an arrow in $Q$ starting at $v$ and ending at a vertex of $A$ or  an arrow starting at a vertex of $A$ and ending at a $v$.
\end{df}
 \textit{Proof of Proposition \ref{prop maximal norms} (a):} Let $l\geq 3$.  By \cite[Corollary 3.36]{Dimitrov}  for any $l\geq 3$ there exists a an exceptional pair $(E_0,E_1)$ in $D^b(Q)$ such that $\hom^{min}(E_0,E_1)\geq l$.  In \cite{WCB1} is shown that $(E_0,E_1)$ can be extended to a full exceptional collection. Therefore we can apply Corollary \ref{from exc collections to maximal norm}  to $D^b(Q)$.
$\square$\vspace{3mm}

Now we present one  method (Lemma \ref{criteria for maximal norm}) to obtain $l$-Kronecker pairs  with arbitrary big $l$ as part of full exceptional collections, i.e. method to obtain the conditions of Corollary \ref{from exc collections to maximal norm}.  This method relies on  full exceptional  collections in which a triple remains strong after certain mutations (see (c) in the statement of Lemma \ref{criteria for maximal norm}).  In \cite{BP} a strong exceptional collection  $\mc E$ which  remains strong under all mutations is called \textit{non-degenerate}.   Furthermore  in \cite{BP} are  defined  so called \textit{ geometric}  exceptional collections and \cite[Corollary 2.4]{BP} says that geometricity implies non-degeneracy.  Furthermore, \cite[Proposition 3.3]{BP} 
 claims that a full exceptional collection of length $m$ of coherent sheaves on a smooth projective variety $X$ of dimension $n$ is geometric if and only if $m=n+1$.     In particular it follows:  
	\begin{remark} \label{one geometric and one non-geometric} The full exceptional collection $\mc E=\{\mc O, \mc O(1),\dots, \mc O(n)\}$ in $D^b(\PP^n)$ introduced by Beilinson   \cite{Beili} is geometric and therefore non-degenerate, whereas the  well known (see \cite{Rudakov}, \cite{Gorod}) strong  full exceptional collection   of line bundles
		$(\mc O(0,0), \mc O(0,1), \mc O(1,0), \mc O(1,1) )$  in $D^b(\PP^1\times \PP^1)$ is not geometric.  \end{remark}
That's why the method of Lemma \ref{criteria for maximal norm} is readily applied to $D^b(\PP^n)$, whereas applying it to $D^b(\PP^1\times \PP^1)$  requires some additional 	 arguments to ensure (c) in Lemma \ref{criteria for maximal norm}.
	
\begin{lemma} \label{criteria for maximal norm} Let $\mc T$ be a proper triangulated category and  $\varepsilon \in (0,1)$.  Let  $\mc E=(F_0,F_1,F_2,E_3,\dots,E_n)$ be   a   full exceptional collection with $n\geq 3$. Let  $\{F_i\}_{i\in \NN}$ be a sequence starting with $F_0,F_1,F_2$ and $F_{i+1}=R_{F_{i}}(F_{i-1})$ for $i\geq 2$.  If the following three properties hold:

{\rm (a)}  $\hom(F_0,F_1)<\hom(F_0,F_2)$;  {\rm (b)} $l=\hom(F_1,F_2)\geq 2$;
 {\rm (c)} $(F_0,F_i,F_{i+1})$ is strong for all $i\geq 1$,
then $\mc T$ satisfies  the condition of Corollary \ref{from exc collections to maximal norm} and $\norm{\mc T}_{\varepsilon}=\pi (1-\varepsilon)$.
\end{lemma}
\bpr 
 Now \eqref{L_A and R_B} becomes  \be \label{L_A and R_B op}    \bd F_{i-1} &\rTo^{coev^*_{F_{i-1},F_{i}}} &\Hom^*(F_{i-1},F_{i})^{\check{}}\otimes F_{i} &\rTo & R_{F_{i}}( F_{i-1})=F_{i+1} \qquad i\geq 2. \ed \ee 
Since the property of being full is preserved under mutations, it follows that  $(F_0,F_{i-1},F_i,E_3\dots,E_n)$ is   full  for each $i\geq2$.   We will show that \eqref{increasing sequence}  holds, and then our $\mc T$ satisfies the conditions of  Corollary \ref{from exc collections to maximal norm}, hence $\norm{\mc T}_{\varepsilon}=\pi (1-\varepsilon)$. 
\begin{gather}\label{increasing sequence}  i\in \NN_{\geq 2} \ \ \Rightarrow \ \  \hom(F_{0},F_{i-1})< \hom(F_{0},F_{i})  \end{gather} 
To show \eqref{increasing sequence} we first note that due to  (c) we have $\hom^k(F_{i-1},F_i)=0$ for each  $k\neq 0$ and each $i\geq 2$ and  it follows that (see e.g. \cite[Example 2.7]{BP}) $l=\hom(F_{1},F_{2})=\hom(F_{i-1},F_{i})=\hom(F_{i},F_{i+1})$ for each $i\geq 2$ and then \eqref{L_A and R_B op} has the form:
  \be \label{L_A and R_B opa}    \bd F_{i-1} &\rTo^{coev^*_{F_{i-1},F_{i}}} & F_{i}^{\oplus l} &\rTo & F_{i+1} \ed \qquad \qquad i\geq 2.\ee 
	In (a)  we are given $\hom(F_{0},F_{i-1})<\hom(F_{0},F_{i})$ for $i=2$ and we will show \eqref{increasing sequence} by induction. Indeed, since $(F_0,F_{i-1},F_i)$ is a strong exceptional collection for each $i\geq 2$,  applying $\Hom(F_0,\_)$ to \eqref{L_A and R_B opa}   yields  short exact sequences between finite dimensional vector spaces:
	\be \label{L_A and R_B vector spaces}    \bd 0 &\rTo  &\Hom(F_0,F_{i-1}) &\rTo & \Hom(F_0,F_{i})^{\oplus l} &\rTo &  \Hom(F_0,F_{i+1}) & \rTo &  0, \ed \qquad i\geq 2. \ee 
 The obtained exact sequences and  $l \geq 2$ imply: 
	\begin{gather}\nonumber  \hom(F_0,F_{i+1})=l \hom(F_0,F_{i}) - \hom(F_0,F_{i-1})\geq 2 \hom(F_0,F_{i}) - \hom(F_0,F_{i-1}) \\[-2mm]
	\label{computation} \\[-2mm] \nonumber  =\hom(F_0,F_{i}) + (\hom(F_0,F_{i}) - \hom(F_0,F_{i-1})), \end{gather}
	hence for  $i\geq 2$ the inequality $\hom(F_0,F_{i}) > \hom(F_0,F_{i-1})$ implies $\hom(F_0,F_{i+1}) > \hom(F_0,F_{i})$. 
The lemma is proved.\epr

\textit{Proof of Proposition \ref{prop maximal norms} (b).} In Remark \ref{one geometric and one non-geometric} is given a full strong  exceptional collection $\mc E$ on $D^b(\PP^n)$ which remains strong under all mutations. Using \cite[Example 2.9]{BP} one computes  $\hom(\mc O,\mc O(1))=\hom(\mc O(1),\mc O(2))=n+1<\hom(\mc O,\mc O(2))=\frac{(n+1)(n+2)}{2}$.  Therefore we can apply Lemma \ref{criteria for maximal norm} and the corollary follows. $\square$ \vspace{2mm}

\textit{Proof of Proposition \ref{prop maximal norms} (c).} Let us denote here $\mc T=D^b(\PP^1 \times \PP^1)$. Exceptional collections in $\mc T$ have been studied in  \cite{Rudakov} and \cite{Gorod}. In particular  the full strong exceptional collection \\ $(\mc O(0,0), \mc O(0,1), \mc O(1,0), \mc O(1,1) )$ mentioned in Remark \ref{one geometric and one non-geometric}  satisfies $\hom(\mc O(0,0),\mc O(0,1))=$\\ $\hom(\mc O(0,1),\mc O(1,1))=2$ and \\  $\hom(\mc O(0,0),\mc O(1,1))=4$ (see \cite[p. 3]{perling} or \cite[Example 6.5]{B}). After one mutation we get a full exceptional collection $(F_0,F_1,F_2,E_3)$ in which $(F_0,F_1,F_2)$ is strong,   $\hom(F_0,F_1)<\hom(F_0,F_2)$, and   $\hom(F_1,F_2)=2$.  Let  $\{F_i\}_{i\in \NN}$ be a sequence starting with $F_0,F_1,F_2$ and $F_{i+1}=R_{F_{i}}(F_{i-1})$ for $i\geq 2$. To apply Lemma \ref{criteria for maximal norm}   and deduce that $\norm{\mc T}_{\varepsilon}=\pi (1-\varepsilon)$ we need  to show that $(F_0,F_i,F_{i+1})$ is strong for all $i\geq 1$. 

 From \cite[Proposition 5.3.1, Theorem 3.3.1.]{Gorod} it follows that:
\begin{gather} \label{For each exceptional pair} \mbox{\textit{For each exceptional pair}}  \  (E,F) \ \mbox{\textit{in}} \  \mc T \  \mbox{\textit{there is at most one}} \  i\in \ZZ \mbox{\textit{ with}} \ \hom^i(E,F)\neq 0.\end{gather}
From the way we defined $\{F_i\}_{i\in \NN}$ it follows (see e.g. \cite[Example 2.7]{BP}) $2=\hom(F_{1},F_{2})=\hom(F_{i-1},F_{i})=\hom(F_{i},F_{i+1})$ for all $i\geq 2$, hence taking into account \eqref{For each exceptional pair}, to show that $(F_0,F_i,F_{i+1})$ is strong for all $i\geq 1$ suffices to show that $\hom(F_0,F_i)\neq 0$ for each $i\geq 1$.  Now \eqref{L_A and R_B} becomes distinguished triangle 
  \be \label{L_A and R_B opala}    \bd F_{i-1} &\rTo^{coev^*_{F_{i-1},F_{i}}} & F_{i}^{\oplus 2} &\rTo & F_{i+1} &\rTo & F_{i-1}[1] \ed \qquad  i\geq 2.\ee 
	We have $0<\hom(F_0,F_1)<\hom(F_0,F_2)$. Assume that for some $i\geq 2$ holds  \begin{gather} \label{something for induction} 0<\hom(F_0,F_1)<\dots<\hom(F_0,F_{i-1})<\hom(F_0,F_{i})\end{gather} we will show that this implies $\hom(F_0,F_{i})<\hom(F_0,F_{i+1})$ and by induction the corollary follows. Applying $\Hom(F_0,\_)$ to \eqref{L_A and R_B opala} and since $\hom^k(F_0,F_{i-1})=\hom^k(F_0,F_{i})=0$ for $k\neq 0$ one easily deduces that $\hom^k(F_0,F_{i+1})=0$ for $k\not \in \{ -1, 0\}$. If $\hom^{-1}(F_0,F_{i+1})\neq 0$, then by \eqref{For each exceptional pair} it follows that $\hom(F_0,F_{i+1})=0$ and applying $\Hom(F_0,\_)$ to \eqref{L_A and R_B opala} yields an exact sequence of vector spaces:
 \begin{gather}  \label{L_A and R_B opalala}    \bd 0 & \rTo &  \Hom^{-1}(F_0,F_{i+1}) & \rTo &  \Hom(F_0,F_{i-1}) &\rTo & \Hom(F_0,F_{i})^{\oplus 2} &\rTo & \Hom(F_0,F_{i+1})=0  \ed ,\end{gather}
which contradicts \eqref{something for induction}. Therefore $\hom^{-1}(F_0,F_{i+1})=0$ and  $\hom^k(F_0,F_{i+1})=0$ for $k\neq 0$. Now we apply  $\Hom(F_0,\_)$ to \eqref{L_A and R_B opala} again and get a short exact sequence as in \eqref{L_A and R_B opa} 
which by the same computation as in \eqref{computation}  implies  
$ \hom(F_0,F_{i+1})>  \hom(F_0,F_{i})$, 
thus we proved the corollary.
$\square$ \vspace{2 mm}

\begin{lemma} \label{help lemma for rational surfaces}
	Let $X$ be a smooth algebraic variety s. t. $D^b(X)$ satisfies the conditions of Corollary \ref{from exc collections to maximal norm}. Let $\wt{X}$ be obtained from $X$ by blowing up a point. Then  $D^b(\wt{X})$ satisfies the conditions of Corollary \ref{from exc collections to maximal norm} as well.
\end{lemma}
\bpr  \cite[Theorem 4.2]{BondalOrlov} ensures that there is a semi-orthogonal decomposition 

$D^b(\wt{X})=\langle \mc T_1, \mc T_2,\dots,\mc T_k,  D^b(X)\rangle $, where $\mc T_i$ is equivalent to $D^b(point)$ for $i=1,2,\dots,k$, which implies that $\mc T_i$ is generated by an exceptional object for each $i$. Now it is clear that the full exceptional sequences of  $D^b(X)$ ensuring the conditions of Corollary \ref{from exc collections to maximal norm}  extend to  full exceptional collections on $D^b(\wt{X})$, so these conditions are satisfied in  $D^b(\wt{X})$ as well. \epr

\textit{Proof of Proposition \ref{prop maximal norms} (d), (e), and (f).} Since ${\mathbb F}_0=\PP^1\times \PP^1$ and ${\mathbb F}_1$ is $\PP^2$ blown up at a point, then the cases $a=0,1$ are contained in Proposition  \ref{prop maximal norms} (a), (b), and Lemma \ref{help lemma for rational surfaces}. In  \cite{HP} they construct families of  full exceptional collections of invertible sheaves on $D^b( {\mathbb F}_a)$ for any $a$. To show that their exceptional  collections furnish the conditions of Corollary \ref{from exc collections to maximal norm} we just need to combine some results in  \cite{HP}. First adopt here some notations and terminology from \cite{HP}: $P,Q$ denotes basis of $Pic({\mathbb F}_a)$ (see \cite[Section 4, p.1224]{HP}) s.t. \begin{gather}\label{hom in hirzebruch 0} P.Q=1, Q^2=a, P^2=0. \end{gather} 
Hille and Perling study certain sequences of Cartier divisors on a rational surfaces $X$ which they call toric systems. Furthermore a toric system $A_1,A_2,\dots, A_n$ is called strongly  exceptional (see \cite[Definition 3.6]{HP}), if they generate a  sequence of invertible sheaves $$\mc O_X, \mc O_X(A_1),\mc O_X(A_1+A_2), \dots ,  \mc O_X(\sum_{i=1}^{n-1} A_i)$$, which is strong exceptional. For such a toric system each divisor $A_j$ is numerically left orthogonal \cite[Definition 3.1 (a)]{HP}, which means that $\chi(-A_j)=0$. Indeed, we have $${\rm Ext}^k\left (O_X(\sum_{i=1}^{j}A_i), O_X(\sum_{i=1}^{j-1}A_i) \right )=\{0\}$$ for each $k\in \ZZ$, since the sequence is exceptional, and on the other hand $${\rm Ext}^k\left (O_X(\sum_{i=1}^{j}A_i), O_X(\sum_{i=1}^{j-1}A_i)\right ) \cong H^k(\mc O_X(-A_j))=\{0\}$$ (see e.g.   \cite[the beginning of Section 3]{HP}), and hence $\chi(-A_j)=0$. Note also that  since the sequence is strong it follows that $\chi (A_j)=\dim(H^0(\mc O_X(A_j)))=\hom\left (O_X(\sum_{i=1}^{j-1}A_i), O_X(\sum_{i=1}^{j}A_i) \right )$.  On the other hand, having that $A_j$ is numerically left orthogonal and using \cite[Lemma 3.3 (i)]{HP} we derive:
\begin{gather} \label{hom in hirzebruch}
\hom\left (O_X\left(\sum_{i=1}^{j-1}A_i\right ), O_X\left (\sum_{i=1}^{j}A_i\right ) \right )=\chi(A_j)=-K_X.A_j
\end{gather}
\cite[Proposition 5.2]{HP} proves that $P, s P +Q, P, -(a+s)P+Q$ is a strongly exceptional toric system  on ${\mathbb F}_a$  when $s\geq -1$. If we denote by $E_1^{s},E_2^{s},E_3^{s},E_4^{s}$  the corresponding strong exceptional collection, then using the formula \eqref{hom in hirzebruch} and the property of toric system, that $\sum_{i=1}^n A_i = -K_X$ (see \cite[p. 1233 down]{HP}), and also the equalities \eqref{hom in hirzebruch 0} we compute:
\begin{gather} \hom(E_2^s,E_3^s)=-K_X.(s P + Q)=(2 (P +Q)-a P).(s P + Q)=2 s = a+ 2 + 2 s\end{gather}
Thus we see that $\hom(E_2^s,E_3^s)$ can be done arbitrary big.  The sequence $E_1^{s},E_2^{s},E_3^{s},E_4^{s}$ is already shown to be full  (see \cite[Theorem 5.8.]{HP}, also \cite[the beginning of the proof of Theorem 8.6.]{HP} or \cite[Proposition 2.1]{HP1} ). Part (d) is proved. Part (e) follows by recursively applying  Lemma \ref{help lemma for rational surfaces} and the already proven cases. Part (f) reduces to part (e),  since any smooth complete rational surface $S$ can be constructed after applying a finite sequence of blow ups starting with $\PP^2$ of $\mathbb F_a$, $a\geq 0$ (see e.g. \cite[the beginning of Section 4, p. 1243]{HP}).
 $\square$\vspace{2 mm}

\section{The inequality \texorpdfstring{$\norm{{\mc T}_{l_1}\oplus\dots\oplus{\mc T}_{l_n}}_{\varepsilon}<\pi (1-\varepsilon)$}{\space} } \label{section sum of Kroneckers}

The goal of this seciton is to prove the following:
\begin{prop} \label{sum of Kroneckers} Let $n\geq 1$, let $l_i\geq 1$, $i=1,2,\dots,n$ be a sequence of integers, and let $0<\varepsilon<1$.  Then for any orthogonal decomposition of the form $\mc T ={\mc T}_{l_1}\oplus{\mc T}_{l_2}\oplus\dots\oplus{\mc T}_{l_n}$,  where $\mc T_{l_i} \cong D^b(K(l_i))$, holds $\norm{\mc T}_{\varepsilon}<\pi (1-\varepsilon)$. Furthermore  $\norm{\mc T}_{\varepsilon}>0$ iff $l_i\geq 3$ for some $1\leq i\leq n$. 
\end{prop}
Before going to the proof of this proposition we prove some facts for the case   $l\geq 3$ and denote  $\mc T_l=D^b(K(l))$.  We will use notations and results from  Section \ref{norms and stab cond on Kroneckers}. 
The first step  is:
\begin{lemma}\label{u_sigma - v_sigma over}   For $\sigma\not \in \mc Z$ the set  $P_\sigma^l=\ol{P_\sigma^l}$ is finite. Otherwise, for $\sigma \in \mc Z$, we use the description of  the set $\ol{P_\sigma^l}$  as  in  Proposition \ref{lemma for P_sigma}  \eqref{e1}, \eqref{e2}, \eqref{e3}.

  For any $0<\varepsilon<1$ there exists $M_{l,\varepsilon} >0$ such that for any $\sigma \in \mc Z$  and  for any  $j\in \ZZ$ : 
\begin{gather}\label{formula for u_sigma - v_sigma over} \phi_\sigma(s_{j+1})-\phi_\sigma(s_{j})>\varepsilon \ \ \ \Rightarrow \ \ \  \frac{v_\sigma - u_\sigma}{u_\sigma-\pi \phi_\sigma(s_{j+1}[-1])}\leq M_{\varepsilon,l} \qquad \frac{v_\sigma - u_\sigma}{\pi \phi_\sigma(s_{j})-v_\sigma}\leq M_{\varepsilon,l}.  \end{gather} 
\end{lemma}
\bpr The part of the lemma, which is not contained  in Proposition \ref{lemma for P_sigma}  are the inequalities \eqref{formula for u_sigma - v_sigma over}.  So, let us chose  $\sigma\in \mc Z$, $j\in \ZZ$ and $0<\varepsilon<1$ and assume that $\phi_\sigma(s_{j+1})-\phi_\sigma(s_{j})>\varepsilon$.
 In terms of the  function \eqref{Function} we can rewrite \eqref{e4} as follows:
\begin{gather} \label{v_sigma-u_sigma..} \frac{v_\sigma - u_\sigma}{u_\sigma-\pi \phi_\sigma(s_{j+1}[-1])}=\frac{F(x,y,a_l)-F(x,y,a_l^{-1})}{F(x,y,a_l^{-1})-\arccos(y)}; \frac{v_\sigma - u_\sigma}{\pi \phi_\sigma(s_{j})-v_{\sigma}}= \frac{F(x,y,a_l)-F(x,y,a_l^{-1})}{\pi -F(x,y,a_l)}
\end{gather}
where (recall that $\sigma \in \mc Z$ implies $\phi_\sigma(s_{j+1})-\phi_\sigma(s_{j})<1$): \begin{gather} \label{-1<y<cospie}  0<x =\frac{\abs{Z(s_{j+1})}}{\abs{Z(s_j)}} \quad -1<y= \cos\left (\pi  (\phi_\sigma(s_{j+1})-\phi_\sigma(s_{j}))\right )<\cos(\pi \varepsilon): 
\end{gather}
For any   $a\in (0,+\infty)$, $b\in (-1,+1)$ the differentiable  functions $(0,+\infty)\ni t \mapsto F(a,b,t)$ and $(0,+\infty)\ni t \mapsto F(t,b,a)$ can be extended uniquely to continuous functions in $[0,+\infty)$ having values $\arccos(b)$ and  $\pi$ at $0$, respectively,  and therefore we can apply the mean value theorem to these functions. More precisely, if $h:[0,+\infty)\rightarrow \RR $ is a function obtained in such a way, then for any $0\leq \alpha < \beta < +\infty$ there exists $\alpha < t < \beta$, such that $h(\beta)-h(\alpha)=(\beta - \alpha) h'(t)$.  In particular, for any $x, y$  as in \eqref{-1<y<cospie}  we can represent all the differences in \eqref{v_sigma-u_sigma..}  as follows (recall \eqref{partial xt F}):
\begin{gather}\label{0<x'<x} \frac{1}{\pi-F(x,y,a_l)}= \frac{-1}{F(x,y,a_l)-\pi}=\frac{1}{ x\frac{a_l \sqrt{1-y^2}}{a_l^2+x'^2- 2 a_l x' y} } =\frac{a_l^2+x'^2- 2 a_l x' y}{x a_l \sqrt{1-y^2} }  \ \mbox{for some} \ 0<x'<x \\ \label{ 0<t<a}
 \frac{1}{F(x,y,a_l^{-1})-\arccos(y)}=\frac{1}{a_{l}^{-1} \frac{x \sqrt{1-y^2}}{t^2+x^2- 2 t x y}}=\frac{a_{l} (t^2+x^2- 2 t x y)}{ x \sqrt{1-y^2} }  \ \mbox{for some} \ 0<t<a_{l}^{-1} \\  
F(x,y,a_l)-F(x,y,a_l^{-1})=\frac{(a_l-a_l^{-1}) x \sqrt{1-y^2}}{t'^2+x^2- 2 t' x y}\leq \frac{(a_l-a_l^{-1})x \sqrt{1-y^2}}{t'^2+x^2- 2 t' x \cos(\pi \varepsilon)}\ \mbox{for some} \ a_l^{-1}<t'<a_l. \end{gather}
And now looking back at \eqref{v_sigma-u_sigma..}  we deduce:
\begin{gather} \label{v_sigma-u_sigma...} \frac{v_\sigma - u_\sigma}{u_\sigma-\pi \phi_\sigma(s_{j+1}[-1])}\leq 
\frac{(a_l^2-1)(t^2+x^2- 2 t x y) }{t'^2+x^2- 2 t' x \cos(\pi \varepsilon)} \\
  \label{v_sigma-u_sigma...1} \frac{v_\sigma - u_\sigma}{\pi \phi_\sigma(s_{j})-v_\sigma}\leq \frac{(1-a_l^{-2}) (a_l^2+x'^2- 2 a_l x' y) }{t'^2+x^2- 2 t' x \cos(\pi \varepsilon)}.
\end{gather}
Now since $t'^2+x^2- 2 t' x \cos(\pi \varepsilon)$ gets minimal values for $ t'= x \cos(\pi \varepsilon)$ (with respect to the variable $t'$) and for  $x=t' \cos(\pi \varepsilon)$ (with respect to the variable $x$) we have  $t'^2+x^2- 2 t' x \cos(\pi \varepsilon)\geq  x^2(1-\cos^2(\pi \varepsilon))= x^2\sin^2(\pi \varepsilon)$ and   $t'^2+x^2- 2 t' x \cos(\pi \varepsilon)\geq  t'^2 \sin^2(\pi \varepsilon) \geq   a_l^{-2} \sin^2(\pi \varepsilon) $, therefore:
{\small \begin{gather}   t'^2+x^2- 2 t' x \cos(\pi \varepsilon)\geq \max\{ a_l^{-2} , x^2\}  \sin^2(\pi \varepsilon) \end{gather} and \eqref{v_sigma-u_sigma...}, \eqref{v_sigma-u_sigma...1} can be continued  (recall that  $0<t<a_l^{-1}$ in  \eqref{ 0<t<a}  and $0<x'<x$ in \eqref{0<x'<x}): 
\begin{gather} \nonumber \frac{v_\sigma - u_\sigma}{u_\sigma-\pi \phi_\sigma(s_{j+1}[-1])}\leq 
\frac{(a_l^2-1)(t^2+x^2- 2 t x y) }{\max\{ a_l^{-2} , x^2\}  \sin^2(\pi \varepsilon)}\leq \frac{(a_l^2-1) }{  \sin^2(\pi \varepsilon)}  \sup \left \{  \frac{(t^2+x^2- 2 t x y) }{\max\{ a_l^{-2} , x^2\} }: \begin{array}{c}  t\in (0,a_l^{-1}) \\ x \in (0,+\infty)\\ y\in (-1,\cos(\pi \varepsilon))\end{array} \right \}\\ \nonumber 
 \frac{v_\sigma - u_\sigma}{\pi \phi_\sigma(s_{j})-v_\sigma}\leq \frac{(1-a_l^{-2}) (a_l^2+x'^2- 2 a_l x' y) }{\max\{ a_l^{-2} , x^2\}  \sin^2(\pi \varepsilon)} \leq  \frac{(1-a_l^{-2}) }{  \sin^2(\pi \varepsilon)}  \sup \left \{  \frac{(a_l^2+x'^2- 2 a_l x' y) }{\max\{ a_l^{-2} , x^2\} }: \begin{array}{c}  x\in (0,+\infty) \\ x' \in (0,x)\\ y\in (-1,\cos(\pi \varepsilon))\end{array} \right \} 
\end{gather} 
}
hence \eqref{formula for u_sigma - v_sigma over} follows.
\epr
\begin{coro} \label{the closed arcs} For any $\sigma \in \mc Z$ there is  closed $\frac{\vol(\ol{P_\sigma^l})}{2}$-arc $p_\sigma^l\subset \ol{P_\sigma^l}$ s.t.   $ \ol{P_\sigma^l}\setminus ( p_\sigma^l \cup - p_\sigma^l ) $ is countable.

Let $0<\varepsilon<1$. For any closed $\varepsilon$-arc $\gamma$ satisfying $P_\sigma^l \cap \gamma = \emptyset$ hold  $(p_\sigma^l\cup-p_\sigma^l )\cap(\gamma \cup-\gamma)=\emptyset$ and  any (of the four) connected component $c$ of $\SS^1\setminus\{p_\sigma^l\cup-p_\sigma^l\cup\gamma \cup-\gamma\}$ restricts $\vol(\ol{P_\sigma^l})$ as follows: \begin{gather}  c\subset \SS^1\setminus\{p_\sigma^l\cup-p_\sigma^l\cup\gamma \cup-\gamma\} \ \ \pi_0(c)=\{0\} \qquad \Rightarrow \qquad  \frac{\vol(\ol{P_\sigma^l})}{2} = \vol(p_\sigma^l)\leq M_{l,\varepsilon} \vol(c)\end{gather}  where $M_{l,\varepsilon}$ is  as in Lemma  \ref{u_sigma - v_sigma over}.
\end{coro}
\bpr For $\sigma \in \mc Z$, the set $\ol{P_\sigma^l}$ is as described in \eqref{e1}, \eqref{e2}, \eqref{e3} and then we can choose $p_\sigma^l$ to be $\exp(\ri[u_\sigma,v_\sigma])$ and $\SS^1$ can be divided  as follows (for any $j\in \ZZ$):
\begin{gather}\nonumber \SS^1 = {\rm e}^{\ri \pi [\phi_\sigma(s_{j}[-1]), \phi_\sigma(s_{j+1}[-1]))}\cup 
{\rm e}^{\ri [\pi  \phi_\sigma(s_{j+1}[-1]), u_\sigma)}\cup p_\sigma^l\cup {\rm e}^{\ri(v_\sigma,\pi \phi_\sigma(s_{j}))} \\[-2mm]\label{division of circle} \\[-2mm] \nonumber  \cup -{\rm e}^{\ri \pi [\phi_\sigma(s_{j}[-1]), \phi_\sigma(s_{j+1}[-1]))}\cup 
-{\rm e}^{\ri [\pi  \phi_\sigma(s_{j+1}[-1]), u_\sigma)}\cup \left ( -p_\sigma^l \right )\cup -{\rm e}^{\ri(v_\sigma,\pi \phi_\sigma(s_{j}))} \end{gather}

 Furthermore, let  $\gamma$ be a  closed  $\varepsilon$-arc with  $P_\sigma^l\cap \gamma=\emptyset$, then  using \eqref{e2} one easily sees that  $\gamma \subset \exp(\ri \pi (\phi_\sigma(s_j),\phi_\sigma(s_{j+1})))$ or $-\gamma \subset \exp(\ri \pi (\phi_\sigma(s_j),\phi_\sigma(s_{j+1})))$ for some $j\in \ZZ$   and therefore $\pi (\phi_\sigma(s_{j+1})-\phi_\sigma(s_{j}))>\vol(\gamma)=\pi \varepsilon$, hence by Lemma \ref{u_sigma - v_sigma over} follow the inequalities \eqref{formula for u_sigma - v_sigma over} and 
\begin{gather} \gamma \subset {\rm e}^{\ri \pi [\phi_\sigma(s_{j}[-1]), \phi_\sigma(s_{j+1}[-1]))} \qquad \mbox{or} \qquad  -\gamma \subset  {\rm e}^{\ri \pi [\phi_\sigma(s_{j}[-1]), \phi_\sigma(s_{j+1}[-1]))}. \end{gather} 
Therefore, taking into account the disjoint union \eqref{division of circle} we see that  the four components of $\SS^1\setminus \{ \gamma \cup -\gamma \cup p_\sigma^l \cup -p_\sigma^l \}$ can be ordered as $c_1$, $c_2$, $-c_1$, $-c_2$ so that: $c_1 \supset {\rm e}^{\ri [\pi  \phi_\sigma(s_{j+1}[-1]), u_\sigma)}$, $c_2\supset  {\rm e}^{\ri(v_\sigma,\pi \phi_\sigma(s_{j}))}$, in particular: 
\begin{gather} \vol(\pm c_1) \geq  u_\sigma-\pi  \phi_\sigma(s_{j+1}[-1]) \qquad  \vol(\pm c_2) \geq \pi \phi_\sigma(s_{j})-v_\sigma\end{gather} and the corollary follows from  \eqref{formula for u_sigma - v_sigma over}.
\epr
\textit{ Proof of Proposition  \ref{sum of Kroneckers} } From Remark \ref{equal norms for isomorphic} and Subsection  \ref{norms on kroneckers} we see that: 
\begin{gather}  \norm{{\mc T}_{l_i}}_\varepsilon=\norm{D^b(K(l_i))}_{\varepsilon}=K_\varepsilon(l_i) \end{gather} hence   the proposition follows  for $n=1$.

 Assume that we have already proved the proposition for  $1\leq n \leq N$.  And assume that  $\mc T ={\mc T}_{l_1}\oplus{\mc T}_{l_2}\oplus\dots\oplus{\mc T}_{l_N}\oplus {\mc T}_{l_{N+1}}$,  where $\mc T_{l_i} \cong D^b(K(l_i))$ and denote by $L$ the set $L=\{l_1,l_2,\dots,l_N,l_{N+1}\}$.  

If  $1\leq l_j\leq 2$ for some $j$, then $\norm{\mc T_{l_j}}_\varepsilon=\norm{D^b(K(l_j))}_\varepsilon=0$, and  the statement follows from the induction assumption, Corollary \ref{sum of several},  and $\norm{{\mc T}_{l_j}}_\varepsilon = 0$. Therefore we can assume that all integers in $L$ are at least $3$.  From the induction assumption there exists $\delta >0$ such that:
\begin{gather} \delta + X = \pi(1-\varepsilon), \ \ \mbox{where } \nonumber \\[-2mm]\label{sum of kroneckers 1} \\[-2mm]
X={\rm max}\left  \{\norm {\mc T_{x_1} \oplus \mc T_{x_2} \oplus \cdots \oplus \mc T_{x_j}}_\varepsilon :   j<N+1, x_i \in L \ \mbox{for} \ 1\leq i \leq j\right \} \nonumber
\end{gather}
Note that due to Remark \ref{phases and equivalences},  Proposition \ref{lemma for orthogonal composition 1} (d), and Corollary \ref{gaps in products} for any sequence $x_1,x_2,\dots,x_j$ in $L$ holds: 
\begin{gather}\norm {\mc T_{x_1} \oplus \mc T_{x_2} \oplus \cdots \oplus \mc T_{x_j}}_\varepsilon = \nonumber \\[-2mm] \label{sum of kroneckers 2} \\[-2mm]
= \sup \left \{\frac{\vol\left (\bigcup_{i=1}^j\ol{P_{\sigma_i}^{x_i}}\right )}{2}: \exists \mbox{closed} \ \varepsilon\mbox{-arc} \ \gamma \ \mbox{s.t.} \ \forall i \ \sigma_i \in \st(D^b(K(x_i)))  \ \mbox{and} \ \emptyset = P_{\sigma_i}^{x_i} \cap \gamma  \right \} \nonumber  \end{gather}
Assume now that $\sigma_i \in \st(D^b(K(l_i)))$ for $i=1,\dots,N+1$ and that there exists a closed $\varepsilon$-arc $\gamma$ satisfying $ \emptyset = P_{\sigma_i}^{l_i} \cap \gamma =\emptyset $ for $i=1,\dots,N+1$.  In particular we can represent the circle $\SS^1$: 
\begin{gather} \label{split the circle} \SS^1 = \exp(\ri (\alpha,\beta)) \cup \gamma \cup  -\exp(\ri (\alpha,\beta)) \cup-\gamma  \qquad \mbox{disjoint union}  \end{gather}
where $\alpha\in \RR$ and $\beta=\alpha+\pi(1-\varepsilon)$. 
 If for some $k$ the corresponding $\sigma_k \not \in \mc Z_{l_k} \subset \st(D^b(K(l_k)))$, then by Lemma \ref{u_sigma - v_sigma over} $\ol{P_{\sigma_k}^{l_k}}$ is finite and taking into account  \eqref{sum of kroneckers 1},  \eqref{sum of kroneckers 2} we derive: \begin{gather} \label{help for a sum of kroneckers}  \frac{\vol\left (\bigcup_{i=1}^{N+1}\ol{P_{\sigma_i}^{l_i}}\right )}{2}= \frac{\vol\left (\bigcup_{i=1,i\neq k}^{N+1}\ol{P_{\sigma_i}^{l_i}}\right )}{2}\leq X, \end{gather} 
otherwise for all $i$ we have $\sigma_i  \in \mc Z_{l_i}$, and then by Corollary \ref{the closed arcs}  $ \frac{\vol\left (\bigcup_{i=1}^{N+1}\ol{P_{\sigma_i}^{l_i}}\right )}{2}= \vol\left (\bigcup_{i=1}^{N+1}p_{\sigma_i}^{l_i}\right )$, where $p_{\sigma_i}^{l_i}$ is a closed arc as explained in  Corollary \ref{the closed arcs} and we can assume that  $p_{\sigma_i}^{l_i}\subset \exp(\ri (\alpha,\beta))$ for all $i$ (see \eqref{split the circle}). There exist $\delta_->0$, $\delta_+>0$ such that    $\bigcup_{i=1}^{N+1}p_{\sigma_i}^{l_i}\subset \exp(\ri [\alpha+\delta_-,\beta-\delta_+])$,  $\exp(\ri (\alpha+\delta_-), 
\exp(\ri (\beta-\delta_+))$ $\in \bigcup_{i=1}^{N+1}p_{\sigma_i}^{l_i}$. Let $j,k$ be such that  $\exp(\ri (\alpha+\delta_-) \in p_{\sigma_j}^{l_j}$ and  $\exp(\ri (\beta-\delta_+) \in p_{\sigma_k}^{l_k}$ If we denote $M={\rm max}\{ M_{l_i,\varepsilon}:1\leq i \leq N+1 \}$, then from Corollary \ref{the closed arcs} we obtain 
$\vol (p_{\sigma_j}^{l_j})+\vol (p_{\sigma_k}^{l_k})\leq M(\delta_++\delta_-)$. Since $\bigcup_{i=1}^{N+1}p_{\sigma_i}^{l_i}\subset \exp(\ri [\alpha+\delta_-,\beta-\delta_+])$ it follows that $\delta_++\delta_- \leq \pi(1-\varepsilon)- \vol\left (\bigcup_{i=1}^{N+1}p_{\sigma_i}^{l_i}\right )$, therefore we can write:
\begin{gather}  \frac{\vol\left (\bigcup_{i=1}^{N+1}\ol{P_{\sigma_i}^{l_i}}\right )}{2} = \vol\left (\bigcup_{i=1}^{N+1}p_{\sigma_i}^{l_i}\right )\leq \vol (p_{\sigma_j}^{l_j})+\vol (p_{\sigma_k}^{l_k}) + \vol \left ( \bigcup_{i=1,i\neq j,i\neq k}^{N+1}p_{\sigma_i}^{l_i} \right )\nonumber   \\ 
\label{help for a sum of kroneckers 2}  \leq M\left (\pi (1-\varepsilon) - \vol\left (\bigcup_{i=1}^{N+1}p_{\sigma_i}^{l_i}\right ) \right )+  \frac{\vol\left (\bigcup_{i=1,i\neq j,i\neq k}^{N+1}\ol{P_{\sigma_i}^{l_i}}\right )}{2}  \\
\leq M\left (\pi (1-\varepsilon) - \vol\left (\bigcup_{i=1}^{N+1}p_{\sigma_i}^{l_i}\right ) \right ) + X=M\left (\pi (1-\varepsilon) - \frac{\vol\left (\bigcup_{i=1}^{N+1}\ol{P_{\sigma_i}^{l_i}}\right )}{2} \right ) + X \nonumber \end{gather}
The obtained inequalities \eqref{help for a sum of kroneckers}, \eqref{help for a sum of kroneckers 2}, and the formula \eqref{sum of kroneckers 2} with $x_i = l_i$, for $i=1,2,\dots,N+1$ show that for a certain set $Y$ and a real function $G$ on $Y$ we have:
\begin{gather}\norm {\mc T_{l_1} \oplus \mc T_{l_2} \oplus \cdots \oplus \mc T_{l_{N+1}}}_\varepsilon =  \sup  \{G(y): y \in Y \}  \nonumber \\
\forall y \in Y \qquad 0\leq  G(y) \leq \pi(1-\varepsilon); \quad G(y)\leq M(\pi(1-\varepsilon)-G(y)) + X. \nonumber \end{gather}
Now recalling \eqref{sum of kroneckers 1} we get   $G(y)\leq M(\pi(1-\varepsilon)-G(y)) + \pi(1-\varepsilon)-\delta$ for any $y\in Y$, which is the same as $G(y)\leq \pi(1-\varepsilon)-\frac{\delta}{M+1}$. Therefore the proof completes with the following inequality:
\begin{gather}0<\norm {\mc T_{l_1} \oplus \mc T_{l_2} \oplus \cdots \oplus \mc T_{l_{N+1}}}_\varepsilon \leq \pi(1-\varepsilon)-\frac{\delta}{M+1}.\end{gather}

\section{Discrete derived categories and their norms} \label{discrete derived categories}

There are categories, in which every heart of a bounded t-structure has finitely many indecomposable objects up to isomorphism. Due to the following lemma the norm of these categories vanishes:

\begin{lemma} \label{lemma the set of phases with indec}  For any triangulated category $\mc T$ and any  $a \in \RR$ we have:
\begin{equation} \label{the set of phases with indec} P_\sigma^{\mc T}= \{\pm \exp(\ri \pi \phi_\sigma(I)) :  I \in \sigma^{ss} \cap \mc P(a,a+1] \ \mbox{and} \ I \ \mbox{is} \ \mc P(a,a+1]\mbox{-indecomposable} \}.\end{equation} 
\end{lemma}
\bpr
From \cite[Lemma 3.9]{DHKK} we know that
\begin{equation} \label{the set of phases with indec1} P_\sigma^{\mc T}= \{\exp(\ri \pi \phi_\sigma(I)) : I \ \mbox{is} \ \mc T-\mbox{indecomposable and} \ I \in \sigma^{ss}\}.\end{equation} 
Furthermore, the properties that for any $j\in \ZZ$ holds $I\in \sigma^{ss}$ iff $I[j]\in \sigma^{ss}$,  and that  $\phi_\sigma(I[j])=\phi_\sigma(I)+j$ whenever $I \in \sigma^{ss}$ are  axioms of  Bridgeland , which together with \eqref{the set of phases with indec1}  imply 
\begin{equation} \label{the set of phases with indec12} P_\sigma^{\mc T}= \{\pm \exp(\ri \pi \phi_\sigma(I)) : I \ \mbox{is} \ \mc T-\mbox{indecomposable and} \ I \in \sigma^{ss} \mbox{and} \ \phi_\sigma(I)\in (a,a+1]\}.\end{equation} 

 From \cite[Lemma 3.7]{DK1} it follows  that an object  $I\in \mc P(a,a+1]$  is $\mc T$-indecomposable iff it is  $\mc P(a,a+1]$-indecomposable, hence the lemma follows.
\epr

\begin{coro} \label{class of vanishing norm} Let $\mc T$ be a  category with phase gap, s.t. every heart of a bounded t-structure has finitely many indecomposable objects up to isomorphism. Then $\norm{\mc T}_\varepsilon = 0$ for every $\varepsilon\in (0,1)$. 
\end{coro}
\bpr First recall that for each $\sigma = (Z,\mc P)\in \st(\mc T)$ and for any $a\in \RR$ the subcategory $\mc P(a,a+1]$ is a heart of a bounded t-structure.   From the previous lemma $P_\sigma^{\ms T}$ is finite for each $\sigma \in \st(\mc T)$. Therefore from Definition \ref{main def eq int} it follows that $\norm{\mc T}_\varepsilon = 0$.
\epr

In representation theory was introduced  a class of triangulated categories with a particularly discrete structure, called Discrete derived categories (Vossieck \cite{Vossieck}), they were classified in \cite{BGS} and thoroughly studied in \cite{BPP1}, whereas the topology of the stability spaces on them were studied in \cite{BPP2}, \cite{Woolf}, in particular it was shown that these spaces are all contractible.  This class contains the categories $\{D^b(Q): Q \ \mbox{is Dynkin}\}$, and the discrete derived categories not contained in this list are of the form $D^b(\Lambda(r,n,m))$ for $n\geq r \geq 1$ and $m\geq 0$, where $\Lambda(r,n,m)$ is the path algebra of the quiver with relations shown on \cite[Section 4.3, Figure 1]{Woolf}.

\begin{prop} \label{norm of discrete categories}  For any discrete derived category $\mc T$ (in the sense of \cite{Vossieck}, \cite{BGS}) and any $\varepsilon \in (0,1)$ holds $\norm{\mc T}_\varepsilon = 0$.
\end{prop}
\bpr  \cite[Proposition 7.1]{BPP1} says that each heart of a bounded t-structure  in $\mc T$ has finitely many indecomposable objects and is of finite  length. In articular (see Lemma \ref{criteria for phase gap})  $\mc T$ has a phase gap and it satisfies the conditions of Lemma \ref{class of vanishing norm}, therefore   $\norm{\mc T}_\varepsilon = 0$.
\epr

\section{Topology on the class of triangulated categories with a phase gap} \label{topology section}
In this section we denote by $\mathfrak{T}'$ the set of all small triangulated  categories within a certain universe (a universe which contains the derived categories or representations of algebras) and by $\mathfrak{PG}'\subset \mathfrak{T}'$ we denote the subset of proper  categories with finite rank Grothendieck group and with a phase gap. From Proposition \ref{norm of discrete categories} (see also its proof) it follows that each discrete category is in $\frak{PG}'$. Furthermore, from \cite[Proposition 7.6]{BPP1} it follows that each discrete derived category has a full Exceptional collection.  Thus if we denote  by $\frak{DDK}'$ the subset in $\frak{T}'$ of discrete derived categories, and by  $\frak{E}'$ the subset of proper  categories with a full exceptional collection, then we have the inclusions: 
$ \frak{DDK}' \subset \frak{E}' \subset \frak{PG}' \subset \frak{T}'. $ 
   When we write $\mc A \cong\mc B$ for $\mc A$, $\mc B \in \mathfrak{T}'$,  we mean an equivalence between triangulated categories, and  by $\mathfrak{T}=\mathfrak{T}'/\cong$ ,  $\mathfrak{PG}=\mathfrak{PG}'/\cong$, $\mathfrak{DDK}=\mathfrak{DDK}'/\cong$ we denote the corresponding sets of equivalence classes and then we have inclusions:
   \begin{gather} \frak{DDK} \subset \frak{E} \subset \frak{PG} \subset \frak{T}. \end{gather}
   We  give first an example of a topology on the largest class $\mc \frak{T}$ and give evidence that  this topology is too coarse.  
  \begin{df} \label{def of topology} For any   $\mc T\in \mathfrak{T}'$ we denote a subset of $\mathfrak{T}'$ as follows :
  \begin{gather} B(\mc T)=\{ \mc T' \in \mathfrak{T}' : \mc T'\cong \mc T  \ \mbox{or there is a SOD} \ \mc T'=\langle \mc A, \mc B \rangle \ \mbox{with} \ \mc A \cong  \mc T \}.
  \end{gather}
  By definition we have  $B(\mc T_1)=B(\mc T_2)$, if $\mc T_1 \cong \mc T_2$.
  \end{df}
\begin{lemma} \label{topology} Let $\mc T'$, $\mc T$ be triangulated categories. If $\mc T'\in B(\mc T)$, then $B(\mc T')\subset B(\mc T)$. In particular, the family of sets $\{ B(\mc T) \}_{\mc T \in \mathfrak{T}'}$ is a base of a topology on $\mathfrak{T}'$, and the family of sets $\{ B(\mc T)/\cong \}_{\mc T \in \mathfrak{T}'}$ is a base of  a topology on $\mathfrak{T}$.
\end{lemma}
\bpr Since $\mc T'\in B(\mc T)$, by definition  there is a SOD $\mc T'=\langle \mc A, \mc B \rangle$ with $\mc A \cong \mc T$. 

Let $\mc T_1 \in B(\mc T')$, therefore  there is a SOD $\mc T_1=\langle \mc C, \mc D \rangle$ with $\mc C \cong \mc T'$, now the SOD $\mc T'=\langle \mc A, \mc B \rangle$ implies a SOD  $\mc C = \langle \mc A', \mc B' \rangle$, where $ \mc A' \cong \mc A \cong \mc T$. Therefore we obtain a SOD  $\mc T_1=\langle \langle \mc A', \mc B' \rangle, \mc D \rangle=\langle  \mc A', \langle \mc B',  \mc D \rangle \rangle$ with $\mc A' \cong \mc T$. i. e.  $\mc T_1 \in B(\mc T) $.
\epr

\begin{lemma} \label{lemma for closed points} If  $\mc T$ is indecomposable with respect to semi-orthogonal decompositions, then   $[\mc T]\in \mathfrak{T}$ is a    closed points w.r. to the topology introduced in Lemma \ref{topology}. If there exists a SOD  $\mc T=\langle \mc A, \mc B \rangle$, where $\mc T \not  \cong \mc A$, then    $[\mc T]\in \mathfrak{T}$ is not a closed point in this topology. 
\end{lemma}
\bpr Let  $\mc T$ be  indecomposable and   $\mc T' \not  \cong \mc T $, then $\mc T \not \in B(\mc T')$ by the  definition \ref{def of topology}, therefore $[\mc T]$ is a closed point indeed.  

Assume that  there exists a SOD  $\mc T=\langle \mc A, \mc B \rangle$ and $\mc T \not  \cong \mc A$, it follows that $[\mc T] \in B(\mc A)/\cong$ and $[\mc T]\not = [\mc A] $, therefore all open subsets containing  $[\mc A]$, contain $[\mc T]$ as well. i. e. $[\mc A]$ is in the closure of  $[\mc T]$, and it is different from  $[\mc T]$, so $[\mc T]$ is not a closed point.  \epr

\begin{coro} The only closed point in $\frak{E}$ w.r. to the topology  introduced in lemma \ref{topology} (also with respect to the induced on $\frak{E}$ topology) is $[\mc T]=[D^b(point)]$. 
\end{coro}\bpr It is well known that $D^b(point)$ is indecomposable with respect to semi-orthogonal decomposition. Now the corollary follows from Lemma \ref{lemma for closed points}. \epr
Now we define a refinement of the topology discussed so far, in which we have many closed points, furthermore we have many discrete subsets, in particular the set of discrete derived categories (up to equivalence) will be  a discrete subset as well.  However this new topology is defined only on $\frak{PG}'$, respectively  $\frak{PG}$. 
\begin{df} \label{balls}
	For any $\varepsilon \in (0,1)$ and any $\mc T\in \frak{PG}' $ we denote \begin{gather} \norm{\mc T}^{\varepsilon}=(1-\varepsilon)\pi - \norm{\mc T}_{\varepsilon}.\end{gather}
	 For any   $\mc T\in \mathfrak{PG}'$, any $\delta>0$ we denote a subset of $\mathfrak{PG}'$ as follows:
	 \begin{gather}\label{upper epsilon}  B^{\varepsilon}_\delta(\mc T)=\{ \mc T' \in \mathfrak{PG}' : \mc T'\cong \mc T  \ \mbox{or} \ \mc T'=\langle \mc A, \mc B \rangle \ \mbox{is a SOD with} \ \mc A \cong  \mc T, \mc B \in \frak{PG}', \ \norm{\mc B}^{\varepsilon}<\delta \}.
	 \end{gather}
	 By definition we have  $B_\delta^\varepsilon(\mc T_1)=B_\delta^\varepsilon(\mc T_2)$, if $\mc T_1 \cong \mc T_2$. Furthermore from Theorem \ref{prop inequality}  it follows that  for any  $\mc T_1$, $\mc T_2 \in \frak{PG}'$ and any SOD $\mc T=\langle \mc T_1, \mc T_2 \rangle$ holds: 
	 \begin{gather} \label{squeezing norm}
	 \norm{\langle \mc T_1, \mc T_2 \rangle}^{\varepsilon}\leq \min \{\norm{ \mc T_1}^{\varepsilon}, \norm{ \mc T_2}^{\varepsilon}\}.
	 \end{gather}
\end{df}

From now on $\varepsilon$ is a real number in  $(0,1)$ and we will write just  $B_\delta(\mc T)$ instead of  $B_\delta^\varepsilon(\mc T)$.  

\begin{lemma} \label{refined topology}
		If $\mc T, \mc T' \in \frak{PG}'$, $\delta >0$ and  $\mc T'\in B_\delta(\mc T)$, then $B_{\delta'}(\mc T')\subset B_\delta(\mc T)$ for any $\delta'>0$. In particular, the family of sets $\{ B_\delta(\mc T) \}_{\mc T \in \mathfrak{PG}', delta >0}$ is a base of a topology on $\mathfrak{PG}'$, and the family of sets $\{ B_\delta(\mc T)/\cong \}_{\mc T \in \mathfrak{PG}', \delta >0}$ is a base of  a topology on $\mathfrak{PG}$.
\end{lemma}
\bpr Since $\mc T'\in B_\delta(\mc T)$, by definition   $\mc T'=\langle \mc A, \mc B \rangle$ with $\mc A \cong \mc T$, $\mc B \in \frak{PG}'$ and $\norm{\mc B}^{\varepsilon}<\delta$. 

Let $\mc T_1 \in B_{\delta'}(\mc T')$, therefore   $\mc T_1=\langle \mc C, \mc D \rangle$ with $\mc C \cong \mc T'$, $\mc D \in \frak{PG}'$.  As in the proof of Lemma \ref{topology} one dervies a    SOD  $\mc T_1= \langle  \mc A', \langle \mc B',  \mc D \rangle \rangle$ with $\mc A' \cong \mc T$, $ \mc B'\cong \mc B$. From \eqref{squeezing norm} we deduce the inequality $\norm{\langle \mc B',  \mc D \rangle}^\varepsilon\leq \norm{\mc B'}^\varepsilon < \delta $, which amounts to the required  $\mc T_1 \in  B_\delta(\mc T)$.
\epr

\begin{prop} \label{upper semi} {\rm (a)} The function below is upper semi-continuous:
	\begin{gather}
	\bd \frak{PG} & \rTo^{\norm{\cdot}^\varepsilon} & \RR & \hspace{6mm} [\mc T] & \mapsto & \norm{\mc T}^\varepsilon \ed	\end{gather}
	
	{\rm (b)} For any $x> 0$ the subset $\frak{PG}_{\geq x}=\{y \in \frak{PG} : \norm{y}^\varepsilon \geq x \}$ is a discrete subset of $\frak{PG}$ w. r. to the topology from Lemma \ref{refined topology}.
\end{prop}
\bpr (a) follows from the following application of  \eqref{squeezing norm}:  for any $\delta >0$, $\mc T \in \frak{PG}'$ holds  \begin{gather} \forall [\mc T'] \in (B_\delta  (\mc T) / \cong)\setminus \{[\mc T]\}  \qquad 0\leq \norm{\mc T'}^\varepsilon\leq \min\{\norm{\mc T}^\varepsilon, \delta\}.  \end{gather}
(b) follows from  the same formula. Indeed, from this formula one checks that for any $[\mc T] \in \frak{PG}_{\geq x}$ and any $0<\delta<x$ we have $\frak{PG}_{\geq x}\cap (B_\delta(\mc T)/\cong) = \{[\mc T]\}$. \epr	

\begin{coro}
	$\frak{DDK} \cup \{[D^b(Q)]: Q \ \mbox{is affine} \}$ is a discrete subset of $\frak{PG}$ with respect to the topology introduced in Lemma \ref{refined topology}. 
\end{coro}
\bpr From Propositions  \ref{norm of Dynkin} and \ref{norm of discrete categories} it follows that for any $[\mc T] \in \frak{DDK} \cup \{[D^b(Q)]: Q \ \mbox{is affine} \} $ has $\norm{\mc T}^{\varepsilon}=\pi (1-\varepsilon)$, hence (since the function \eqref{upper epsilon} takes values in $(0,\pi(1-\varepsilon))$) we obtain: \begin{gather} \frak{DDK} \cup \{[D^b(Q)]: Q \ \mbox{is affine} \} \subset \frak{PG}_{\geq \pi (1-\varepsilon)}.\end{gather}
On the other hand from Proposition \ref{upper semi} (b) we know that $ \frak{PG}_{\geq \pi (1-\varepsilon)}$ is a discrete subset and the corollary follows. 
\epr 
Examples of non-closed points are contained in Proposition \ref{prop maximal norms}. More precisely:
\begin{prop}
	The element $[D^b(point)]\in \frak{PG}$ is in the closure of $[\mc T]\in \frak{PG}$ for  any $\mc T \in \frak{PG}'$ which satisfies the conditions of Corollary \ref{from exc collections to maximal norm} and such that ${\rm rank}(K_0(\mc T))\geq 3$. 
\end{prop}
\bpr We will show that  $[\mc T]\in B_\delta(D^b(point))/\cong$ for any $\delta >0$. Indeed, take any $\delta>0$. From \eqref{props of Ky(l) 4}  it follows that there exists $N$ s.t.  $\pi (1-\varepsilon)-K_\varepsilon(l) < \delta $ for  $l\geq N$. Since $\mc T$ satisfies the  conditions of Corollary \ref{from exc collections to maximal norm} and ${\rm rank}(K_0(\mc T))\geq 3$, therefere there is a full exceptional collection $E_0,E_1,\dots,E_{n-1},E_{n}$ with $n\geq 2$,  s. t. $\hom^{min}(E_i,E_j)>N$ for some $i<j$. Since we can apply mutations, we can assume that  $\hom^{min}(E_{n-1},E_{n})>N$.  Now let us denote $\mc A=\langle E_0 \rangle$, $\mc B = \langle E_1,\dots, E_n \rangle$. Then we have a SOD $\mc T = \langle \mc A, \mc B \rangle$ with $\mc A \cong D^b(point)$, and $\norm{\mc B}^{\varepsilon} \leq \norm{ \langle E_{n-1}, E_n \rangle}^{\varepsilon}\leq \pi (1-\varepsilon) - K_\varepsilon(\hom^{min}(E_{n-1},E_{n}))<\delta$,  where in the latter chain of inequalities we used  \eqref{squeezing norm}, Proposition \ref{inequality for any exc pair}. Recalling the definition  of $B_\delta(D^b(point))$ (  Definition \ref{balls}) we conclude that $\mc T \in B_\delta(D^b(point))$  and the proposition follows. 
\epr
\begin{coro} \label{point in closure of rational}
	For any smooth complete rational surface surface $S$ holds $[D^b(point)] \in \cl([D^b(S)] )$.
\end{coro}

\section{Non-commutative curve-counting} \label{NCCC}

\subsection{Rescaling $\norm{}_{\frac{1}{2}}$ so that all natural numbers are values} 	Following Kontsevich-Rosenberg  \cite{KR} we denote $D^b(K(l+1))$ by  $N\PP^l$ (non-commutative projective space) for $l\geq 0$. Note that we include the case $l=0$, and   $N\PP^0$ is a non-trivial cateogry. Then if we define a function: \begin{gather} \label{dim map} \bd \left \{ \begin{array}{c} \mbox{triangulated} \\ \mbox{categories with} \\  \mbox{a phase gap}  \end{array} \right \} &  \rTo^{ \dim_{nc}} & [0,+\infty]  \ed \\  \dim_{nc}(\mc T) = \left \{  \begin{array}{l l} 
		0 & \mbox{if} \ \forall \ \mbox{full} \  \sigma \in \st(\mc T) \quad  \#(P_\sigma) < \infty \\ \frac{2}{\cos\left (\norm{\mc T}_{\frac{1}{2}}\right )} -1 & \mbox{if} \ \norm{\mc T}_{\frac{1}{2}}<\pi/2 \ \mbox{and} \  \exists \ \mbox{full} \  \sigma \in \st(\mc T) \ \mbox{s.t.} \  \#(P_\sigma) = \infty    \\
		+ \infty & \mbox{if} \  \norm{\mc T}_{\frac{1}{2}}=\pi/2  \end{array}   \right. , \nonumber   \end{gather}
using \eqref{props of Ky(l) 2}, Proposition \ref{norm of Kroneckers}, and table \eqref{table_intro}  we see that 
\begin{gather} \label{dim(NPP^l)} \dim_{nc}(N\PP^l) = l \qquad l \geq 0  \end{gather}

Thus the invariant $\dim_{nc}$  takes all natural numbers as values and due to Theorem \ref{prop inequality} and Remark \ref{stab with infinite Psigma}  we have
\begin{gather} \label{dimensions inequality} \dim_{nc}(\langle \mc A, \mc B \rangle)\geq \max \{\dim_{nc}(\mc A),\dim_{nc}(\mc B) \}, \end{gather}
whenever $\mc A$, $\mc B$ have phase gap and $\langle \mc A, \mc B \rangle$ is a semi-orthogonal decomposition.
\eqref{dimensions inequality}  ensures that whenever $\mc T$ has a finite $\dim_{nc}(\mc T)<+\infty$  and $\mc A \subset \mc T$ is a good enough embedded subcategry, then $\mc A$ has also finite $\dim_{nc}(\mc A)\leq \dim_{nc}(\mc T)<+\infty$. 
 Note that due to  
\begin{remark} \label{nc dim of Dynkin}
	Proposition  \eqref{norm of Dynkin} and table \eqref{table_intro} imply that  for an acyclic  quiver $Q$ we have $\dim_{nc}\left (D^b(Q) \right ) = 0$ iff $Q$ is Dynkin and  $\dim_{nc}\left (D^b(Q) \right ) = 1$ iff $Q$ is affine. 
\end{remark}
\begin{quest}
	Is there a category $\mc T$ with a phase gap s.t. $\dim_{nc}(\mc T) \not \in \QQ$ ?
\end{quest}

\subsection{The general definition and the question which motivates it.} \label{non-commutative curves}
Due to \eqref{dimensions inequality} we see that whenever we have $\mc T \in B_\delta( N\PP^l)$ for some real $\delta >0$ (recall that by definition \ref{balls} this means that there is   a SOD of the form $\mc T=\langle N\PP^l, \mc A \rangle $ where $\mc A$ has a phase gap ) and some integer $l\geq 0$, then $\dim_{nc}(\mc T)\geq l$. In particular if $\mc T \in B_{\delta_l}( N\PP^l)$ for arbitrary big $l$, then $\dim_{nc}(\mc T)=+\infty$, and this idea was used in Section \ref{criteria}.  Now we come to the main question of this section: 

\begin{quest} \label{main question} Let $l\geq 0$. Incidences  of the form  $\mc T \in B( N\PP^l)$ \footnote{i. e.  a SOD  $\langle N\PP^l, \mc A \rangle $ }, are a common phenomenon for this paper.   Recalling that Gromow-Witten invariants count pseudo-holomorphic curves, we view such embeddings of $N\PP^l$ into $\mc T$ as analogous to a  ``pseudo-holomorphic curve'' in the category $\mc T$, and we ask  a question: can we count such entities in a given $\mc T$, how many are they  ? 
\end{quest} 
We answer positively this  question here, and give examples. We plan to develop this idea  in  future works (starting with \cite{DK4}). 
\begin{remark} \label{motivation for genus}
 Recall that the homological dimension of  $ N\PP^l$, $l\geq 0$ is one. Also due to table \eqref{table_intro1}  we have $\st(N \PP^l)\cong \CC\times\CC$ for $l=0,1$ and  $\st(N \PP^l)\cong \CC\times\mc H$ for $l\geq 2$. Note also that, whereas the spirals in $N \PP^0$ are periodic (up to shifts there are only three exceptional objects), for $l\geq 1$ the spirals in $N \PP^l$ consist of pairwise non-isomorphic objects.	
 
 In view of these notes, we find it convenient to view $ N\PP^l$ as a non-commutative curve of genus $l$. 

\end{remark}

Let us also note:
\begin{remark} \label{full exact has lr adjoints} Let $l\geq 0$ and $\mc T$ be any  triangulated cateogry linear over $k$,  let $\bd N \PP^l & \rTo^{F} & \mc T \ed$ be any fully faithful  exact functor \footnote{recall that  an exact functor   is actually a  a pair of a functor $F'$ and  a natural isomorphism between the functors $F'\circ T_1$ and $T_2\circ F_2$, where $T_1$, $T_2$ are the translation functors of the source and the target categories, respectively}, and denote by  $\mc A$ the isomorphism closure of the  image of $F$ in $\mc T$.  Then $\mc A$ is a triangulated subcategory of $\mc T$ generated by two exceptional objects, hence due to \cite[Theorem 3.2]{B} the functor $F$ has   left and right adjoints and  there are SOD $\mc T=\langle \mc A, \mc A^{\perp} \rangle $, $\mc T=\langle  ^{\perp}\mc A, \mc A \rangle $.
	
\end{remark}

\begin{remark} \label{functor from NP tp NP} From the previous remark it follows that, if for some integer  $j\geq 0$ there exists a fully faithful functor $\bd N \PP^j & \rTo^{F} &  N \PP^l \ed$, then $l=j$ and $F$ is equivalence (from \cite{WCB2} we know that any exceptional pair in $N\PP^l$ is full).
\end{remark}
\begin{df} \label{C_l}   Let  $\mc A$, $\mc T$ be any  triangulated categories linear over $k$. And let $\Gamma \subset {\rm Aut}(\mc T)$ be a subgroup of the group of autoequivalences. We denote 
	\begin{gather} \label{C'_l,P} 
		C'_{\mc A,P}(\mc T) = \{\bd \mc A & \rTo^{F} & \mc T \ed : F \ \mbox{is fully faithful exact functor satisfying  properties} \ P \}.
	\end{gather} 
	Next we fix an equivalence relation in $C'_{\mc A ,P}(\mc T)$, and  we will be interested in the set of equivalence classes, in particular the size of this set. 
	\begin{gather}\label{equivealence}
		C_{\mc A,P}^{\Gamma}(\mc T) = C'_{\mc A,P}  (\mc T)/{\sim} \qquad F \sim F' \iff F \circ \alpha \cong \beta\circ F' \ \mbox{for some} \ \ \alpha \in {\rm Aut}(\mc A), \beta \in \Gamma  
	\end{gather}  where $ F \circ \alpha \cong \beta \circ F' $ means equivalence of exact functors between triangulated categories (this is so called graded equivalence).

\end{df}

\subsection{First non-trivial examples with  $\mc A=N\PP^l$, $l\geq 0$ and tree diferent targets: two quivers, and $D^b(\PP^2)$} \label{non-trivial examples} These examples we can prove now. We postpone  writing  down the details of the  proof of Proposition \ref{Cl(PP2)} and \ref{the numbers for two quivers} for a future work \cite{DK4}, where we plan to study also other examples. 
For the  statements till the end of this subsection  we fix the  additional properties $P$ from \eqref{C'_l,P}  as follows,  note that this  additional restriction implies that  $\mc T$  has a phase gap (Proposition \ref{prop inequality} and   Remark \ref{full exact has lr adjoints}):
\begin{gather} \label{Property P}
	\mbox{Property $P$: the left or the right orthogonal to the image of $F$ in $\mc T$ has a  phase gap.}
\end{gather}
Unless otherwise specified in this subsection we take either a trivial subgroup $\Gamma=\{\rm Id\}$,  or the      
subgroup $\Gamma = \langle S \rangle  $  generated by the Serre functor $S$ of $\mc T$ (such exists in the discussed examples), or $\Gamma = {\rm Aut}(\mc T)  $. 

\textit{Remembering that in this subsection $P$ is like  above and $\mc A = N\PP^l$ for some $l\geq 0$ we will  write  $C'_{l}(\mc T)$, $C_{l}^{\Gamma}(\mc T)$ instead of  $C'_{N\PP^l,P}(\mc T)$, $C_{N\PP^l,P}^{\Gamma}(\mc T)$. We refer to  $C_l^{\Gamma}(\mc T)$ as the set of  non-commutative curves of genus $l$ in $\mc T$ modulo the subgroup $\Gamma$. }

This paper contains many examples of $\mc T$ and $l\geq 0$, s. t. $C_{l}^{\Gamma}(\mc T) \neq \emptyset$, on the other hand  from \eqref{dimensions inequality} one sees immediately:
\begin{gather} \label{vanishings of invariants} \dim_{nc}(\mc T) \leq n  \ \  \Rightarrow \ \   C_{l}^{\Gamma}(\mc T) = \emptyset \ \ \mbox{for } \ l > n.  \end{gather}

Furthermore, in all the examples of categories $\mc T$ with $\dim_{nc}(\mc T) = \infty$ given in Section \ref{criteria} one has  $ C_{l}^{\Gamma}(\mc T) \neq \emptyset$ for infinitely many $l$.

Using  results of \cite{DK1} and \cite{DK3}  we can completely determine the invariants $C_{l}^{\{\rm Id\}}(D^b(Q))$,  $C_{l}^{\langle S \rangle}(D^b(Q))$ for  two affine quivers:

\begin{prop} \label{the numbers for two quivers}  Let $\mc T_i=D^b(Q_i)$, $i=1,2$,   where:
	\be \label{Q1} Q_1= \begin{diagram}[1em]
		&       &  \circ  &       &    \\
		& \ruTo &    & \luTo &       \\
		\circ  & \rTo  &    &       &  \circ
	\end{diagram}  \ \qquad \qquad \ Q_2= \begin{diagram}[1em]
	\circ &  \rTo  &  \circ    \\
	\uTo &        & \uTo     \\
	\circ   & \rTo  &    \circ
\end{diagram}. \ee
Then \begin{gather}
\label{D^b(Q_i) zero} \#(C_{0}^{\{\rm Id\}}(\mc T_1)=\infty, \qquad \#(C_{0}^{\{\rm Id\}}(\mc T_2)=\infty, \qquad   \\
\label{D^b(Q_i) zero serre} \#(C_{0}^{\langle S \rangle}(\mc T_1)= 3, \qquad \#(C_{0}^{\langle S \rangle}(\mc T_2)=8, \qquad   \\
\label{D^b(Q_i)} \#(C_{1}^{\{\rm Id\}}(\mc T_1)=2, \qquad \#(C_{1}^{\{\rm Id\}}(\mc T_2)=4, \qquad   \\ \label{D^b(Q_i) serre}
\#(C_{1}^{\langle S \rangle}(\mc T_1)=1, \qquad \#(C_{1}^{\langle S \rangle}(\mc T_2)=2, \qquad \\
l\geq 2, \ \Rightarrow \#(C_{l}^{\{\rm Id\}}(\mc T_1)=\#(C_{l}^{\{\rm Id\}}(\mc T_2)=0, \end{gather} 
\end{prop} 
\bpr
	The vanishings follow from (see Remark \ref{nc dim of Dynkin})  $\dim_{nc}(\mc T_1)= \dim_{nc}(\mc T_2)=1$ and \eqref{vanishings of invariants}. The rest of the proof  will be written in \cite{DK4}.
\epr
\begin{prop} \label{Cl(PP2)} Denote $\mc T = D^b(\PP^2)$. Let $\langle S \rangle \subset  {\rm Aut}(\mc T)$  be the subgroup generated by the Serre functor. Then $\forall l\geq 0$ the set $C_{l}^{\langle S \rangle}(\mc T)$ is finite and 
	\begin{gather} \label{numbers for P2} \{l\geq 0:C_{l}^{\langle S \rangle}(\mc T) \neq \emptyset \}=\{l\geq 0:C_{l}^{\rm Id}(\mc T) \neq \emptyset \}= \{3 m - 1: m \ \mbox{is a Markov number} \}.
	\end{gather}
	Furthermore for any  Markov number\footnote{Recall that a Markov number $x$ is a number $x \in \NN_{\geq 1}$ such that there exist integers $y,z$ with $x^2+y^2+z^2=3 x y z$.}  $m$ we have \begin{gather}  \label{infinite mumbers for PP2} \#\left (C_{3 m -1}^{\{\rm Id\}}(\mc T)\right )=\infty  \\ 
	\label{finite mumbers for PP2 full group} 1\leq \#\left (C_{3 m -1}^{\Aut(\mc T)}(\mc T)\right )= \#\left  \{ y:\begin{array}{l} 0\leq y<m, y \in \ZZ \ \mbox{and there exists an}\\   \mbox{ exceptional vector bundle} \  E \ \mbox{on} \ \PP^2, \\\mbox{with} \  r(E)= m, \  y=c_1(E)\end{array} \right \} \leq  m,
	\\ 
	\label{finite mumbers for PP2} 3\leq \#\left (C_{3 m -1}^{\langle S \rangle}(\mc T)\right )=3 \#\left (C_{3 m -1}^{\Aut(\mc T)}(\mc T)\right ) \leq 3 m,
	\end{gather} 
	where   $c_1(E)$, $r(E)$ are   the first Chern class (which we consider as an integer)  and the rank of $E$.
\end{prop}
\bpr
	Coming in \cite{DK4}.
\epr

\begin{coro} \label{Markov} Denote $\mc T = D^b(\PP^2)$. The first several non-trivial $C_{3m-1}^{\Aut(\mc T)}(\mc T)$ are (recall that $m$ is a Markov number and on the first row are listed the furst 9 Markov numbers): 
\begin{gather}\label{table for PP2}
\begin{array}{| c | c | c | c | c | c | c | c | c | c |}
m                                                & 1 & 2 & 5 & 13  & 29 & 34 & 89 & 169 & 194  \\ \hline
\#\left (C_{3 m -1}^{\Aut(\mc T)}(\mc T)\right ) & 1 & 1 & 2 & 2   & 2  & 2  & 2  & 2   & 2
\end{array}
\end{gather}
Furthermore, the so called Tyurin's conjecture, which is equivalent to the Markov's conjecture, (\cite[p. 100]{Rudakov1} or \cite[Section 7.2.3 ]{GorKul}) is equivalent to  the following statement:  

 For all  Markov numbers $m\neq1, m\neq 2$ we have $\#\left (C_{3 m -1}^{\Aut(\mc T)}(\mc T)\right )=2$. 
 
 Thus Markov Conjecture is true iff this statement is true.
\end{coro}
\bpr
	Coming in \cite{DK4}.
\epr

\begin{remark} \label{automatically P}
	When $\mc T = D^b(\PP^2)$ or $\mc T = D^b(Q)$ for an acyclic quiver $Q$,  then for any fully faithful  exact functor $F:N\PP^l \rightarrow \mc T$ the right and the left orthogonal to the image of $F$ is generated by an exceptional collection, therefore any such functor automatically satisfies the additional property $P$ fixed in the beginning of this subsection. This follows from the fact that every exceptional pair in $\mc T$ can be extended to a full exceptional collection in $\mc T$ (this is proved in \cite{WCB2} and \cite{GoroRuda}).
\end{remark}

 \begin{conj} Let $\mc T = D^b(S_i)$, for $i=1$  or  $i=2$ or  $i=3$, where $S_1,S_2,S_3$ are the quivers from Conjecture \ref{conj for S_1,S_2,S_3}. Remark \ref{automatically P}  holds for $\mc T$ and we take Definition \ref{C_l} with no additional restrictions $P$. We conjecture that $\dim_{nc}(\mc T)<\infty$ and $C_{l}^{\{\rm Id \}}(\mc T)$ is fnite  for all $l\geq 1$.
 \end{conj}

\subsection{Dependence on a stability condition. Semistable non-commutative curves} \label{dependence on stability condition}
\begin{df} \label{Clsigma} Let $l\in \ZZ_{\geq 1}$ and   let $\mc T$ be a  triangulated category linear over $k$ and s.t.  $\ \st(\mc T)\neq \emptyset$, let $\Gamma$ and $P$ be as in Definition \ref{C_l}. One approach to define, semi-stable w.r. to a stability condition non-commutative curves is as follows. Choose   $\sigma \in \st(\mc T)$. Now we apply the same Definition \ref{C_l} with $\mc A=N\PP^l$ and a modified additional properties, namely they are $P$ and one additional restriction depending on $\sigma$. More precisely,    let $\{s_j\}_{j\in \ZZ}$ be a Helix in $N\PP^l$ (see Section \ref{there are no Ext-nontrivial...}), then let us denote: 
	\begin{gather}  	C'_{l,P,\sigma}(\mc T) = \{F \in C'_{N\PP^l,P}(\mc T) : \#\{j: F(s_j)\in \sigma^{ss}\}=\infty  \  \},  \end{gather}	
	where $C'_{N\PP^l,P}$, is the set \eqref{C'_l,P}. The  equivalence relation in $C'_{l,P,\sigma}(\mc T)$ is the same as in \eqref{equivealence}. For completeness we write it again, to define  
	  $\sigma$-semistable non-commutative curves of genus $l$ in $\mc T$, satisfying properties $P$ and modulo $\Gamma$:
	\begin{gather}
	C_{l,P, \sigma
		}^{\Gamma}(\mc T) = C'_{l,P,\sigma}  (\mc T)/{\sim} \qquad F \sim F' \iff F \circ \alpha \cong \beta\circ F' \ \mbox{for some} \ \ \alpha \in {\rm Aut}(N \PP_l), \beta \in \Gamma  
	\end{gather}  where $ F \circ \alpha \cong \beta \circ F' $ means equivalence of exact functors between triangulated categories.
	
\end{df}
We will give two examples. In both of them the additional restrictions $P$ are empty (for both of them holds Remark \ref{automatically P}) and $\Gamma=\{\rm Id\}$, and we will write just $C_{l,\sigma}^{\{\rm Id\}}(\mc T)$ instead of  $C_{l,P,\sigma}^{\{\rm Id\}}(\mc T)$. We give proof only of the basic example:
\begin{prop} \label{kroneckers wc} Let $l\geq 1$. Then $C^{\{\rm Id \}}_{k}(N\PP^l)=\delta_{k,l}$. Let $\mc Z\subset \st(N\PP^l)$ be as in Definition \ref{def of cal Z}. Then we have $C^{\{\rm Id\}}_{l,\sigma}(N\PP^l)=1$ for $\sigma \in \cl(\mc Z)$
 and  $C^{\{\rm Id\}}_{l,\sigma}(N\PP^l)=0$ for $\sigma \not \in \cl(\mc Z)$.	\end{prop}
 \bpr The equality $C^{\{\rm Id\}}_{k}(N\PP^l)=\delta_{k,l}$ follows from Remark \ref{functor from NP tp NP}. For the proof of the rest we note fist that for any $j$ the subset in $\st(N\PP^l)$ where $s_j$ is semi-stable is closed subset.  From  Lemma \ref{sections} we know that for $\sigma \in \mc Z$ all alements in $\{s_j\}_{j
 	\in \ZZ}$ are semi-stable, therefore this holds also for  $\sigma \in \mc \cl(\mc Z)$. 
  Thus we see that  $ \sigma \in  \cl(\mc Z) \Rightarrow C_{l,\sigma}(N\PP^l)=1
 $. Recalling that $\st(N\PP^l) = \mc Z \amalg \amalg_{i\in \ZZ} \left ( \Theta_{(s_i,s_{i+1})}\setminus \mc Z \right )$ we see that if $\sigma \not \in \cl(\mc Z)$, then $\sigma \in \Theta_{(s_i,s_{i+1})}\setminus \mc Z$ for some $i\in \ZZ$. From Lemma \ref{sections} and the description of $\Theta_{(s_i,s_{i+1})}$ in Proposition \ref{lemma for f_E(Theta_E)}  we see that $s_i,s_{i+1} \in \sigma^{ss}$ and $\phi(s_{i+1})>\phi(s_{i})+1$ and then Lemma \ref{sections1} ensures that only $s_i, s_{i+1}$ are semi-stable, therefore $ \sigma \not \in  \cl(\mc Z) \Rightarrow C_{l,\sigma}(N\PP^l)=0.$
  \epr
  
  Using the description of $\st(D^b(Q_1))$ in \cite{DK3} we have proved:
  \begin{prop} \label{Q1 wc} Let $\mc T = D^b(Q_1)$ ($Q_1$ is as in Proposition \ref{the numbers for two quivers}). 
  	  As  $\sigma$ varies in $\st(\mc T)$ the number  $C_{1,\sigma}(\mc T)$ takes all possible values: $\{0,1,2 \}$ (recall that $C_{1}^{\{\rm Id\}}(\mc T)=2$). 
  \end{prop}
  \bpr Coming in a future work.
  \epr
  
  \subsection{Non-commutative Calabi-Yau curve-counting} \label{A-side curve counting}

 Now we give an example of a finite $C_{\mc A,P}^{\Gamma}(\mc T)$ (defined in Definition \ref{C_l}) coming from categories appearing naturally on the A-side.  We first introduce a CY version of the domain category $\mc A$. We pass from $D^b(K(l))$ to the new domain by  changing 
  \[
  \left\{\begin{matrix}
  \text{exceptional} \\
  \text{objects}
  \end{matrix}\right\} \longleftrightarrow 
  \left\{\begin{matrix}
  \text{spherical} \\
  \text{objects}
  \end{matrix}\right\}
  \]
  and this amounts to considering  $\mc A=CY(l)$, instead of  $\mc A=D^b(K(l))$, where $CY(l)$ is defined in \cite{Seidel123}.
  The definition is based on the quiver:
  
  
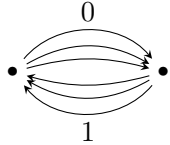
\begin{center}
  	\begin{tikzpicture}
  	\draw (-1, 0) node(a) [scale=0.8]{$\bullet$};
  	\draw (1, 0) node(b) [scale=0.8]{$\bullet$};
  	\draw [-stealth] (a) to[bend left=15] (b);
  	\draw [-stealth] (a) to[bend left=30] (b);
  	\draw [-stealth] (a) to[bend left=50] (b);
  	\draw [-stealth] (b) to[bend left=15] (a);
  	\draw [-stealth] (b) to[bend left=30] (a);
  	\draw [-stealth] (b) to[bend left=50] (a);
  	\draw (0,0.8) node{0};
  	\draw (0,-0.8) node{1};
  	
  	\end{tikzpicture}
  \end{center}
  
  In the next example we take the entire $\Gamma=\Aut(\mc T)$, more precisely:

  \begin{ex}
  	Let $\mc T=\op{Fuk}(E)$ with $E$ an elliptic curve. In this case we have a correspondence:
  	\[
  	C_{CY(n)}^{\Aut(\mc T)}(\mc T) \longleftrightarrow 
  	\left\{\begin{matrix}
  	\text{Primitive  Lagrangian} \\
  	\text{generating } \op{Fuk}(E)
  	\end{matrix}\right\}
  	\]
  	where $C_{CY(l)}^{\Aut(\mc T)}(\mc T)$ is defined in Definition \ref{C_l} and we take empty $P$. It follows that 
  	  	\[\#\left ( C_{CY(l)}^{\Aut(\mc T)}(\mc T)\right )=\#\{d | gcd(d,n) = 1, 1 \leqslant d<n \}.\]
  	
  \end{ex}
  \bpr 	Coming in a future work.
  \epr

  \begin{remark}
  For curves $S$ of higher genus we expect that  finding of the cardinality of  $C_{CY(n)}^{\Aut(\mc T)}(\mc T)$ for $\mc T = {\rm Fuk}(S)$  is related to very recent insights on  counting of geodesics - see \cite{MIZR}.
 \end{remark}

\section{A-side interpretation and holomorphic sheaves of categories} \label{A-side interpretation}

\label{A-side}

In this section we give a different point of view on the category of representations of the Kronecker 
quiver and  introduce the notion of holomorphic families of Kronecker quivers.

 We suggest  a framework in which sequences of holomorphic families of categories are viewed as sequences of extensions of non-commutative manifolds  by relating our norm to the  notion of holomorphic family of categories introduces by Kontsevich. 
Several  questions and conjectures  are posed.

First we   sketch how to interpret  $D^b(K(n))$ as a perverse sheaf of categories.
Recall that LG model of $\PP^2$ is ${\CC^*}^2$, $w=x+y+\frac{1}{xy}$ - see \cite{AKO}. 

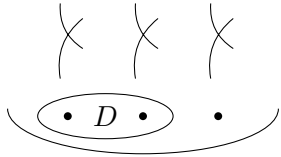
\begin{center}
\begin{tikzpicture}[scale=1]

\draw[xshift=0cm] (-0.1,0.2) to[bend left=30] (0.2,1)(-0.1,1.2) to[bend right=30] (0.2,0.6);
\draw[xshift=-1cm] (-0.1,0.2) to[bend left=30] (0.2,1)(-0.1,1.2) to[bend right=30] (0.2,0.6);
\draw[xshift=1cm] (-0.1,0.2) to[bend left=30] (0.2,1)(-0.1,1.2) to[bend right=30] (0.2,0.6);

\begin{scope}[shift={(0,-0.3)}]
\fill (0,-0) circle (1.5pt);
\fill (1,-0) circle (1.5pt);
\fill (-1,-0) circle (1.5pt);
\draw (-0.5,-0) node{$D$} ellipse (0.9 and 0.3);
\draw (-1.8,0.1) arc (180:360:1.8 and 0.6);
\end{scope}

\end{tikzpicture}
\end{center}

The category $D^b(K(3))$ can be obtained  by taking  the  part of the Landau Ginzburg model over a disc $D$ which contains 2 singular fibers.

 A surgery on the manifold:
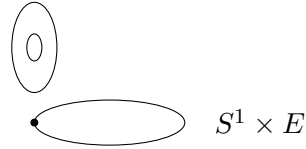
\begin{center}
\begin{tikzpicture}[scale=1]
\draw (0,1) ellipse (0.3 and 0.6) ellipse (0.1 and 0.18);
\fill (0,-0) circle (1.5pt); 
\draw (0,0) arc (180:-180:1 and 0.3) (3,0) node{$S^1 \times E$};
\end{tikzpicture}
\end{center}

results in  changing  the Floer homology  $\op{HF}(L_1, L_2)=3$  to $\op{HF}(L_1, L_2)=4$.
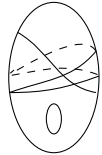
\begin{center}
\begin{tikzpicture}[scale=1]

\draw ellipse (0.6 and 1);
\draw (0,-0.55)  ellipse (0.1 and 0.2);

\draw (-0.48,0.6) to[out=-30,in=170] (0.56,-0.2);
\draw (-0.59,-0.2) to[out=20,in=230] (0.56,0.35);
\draw[dashed] (0.56,0.35) to[out=120,in=220] (-0.53,0.1);
\draw[dashed] (-0.56,0.03) to[out=10,in=160] (0.58,0.04);
\draw (0.6,0.02) to[out=200,in=0] (-0.59,-0.2);

\end{tikzpicture}
\end{center}
As a result we get $D^b(K(4))$. By similar surgeries we can get all quivers from 
$K(0) = A_1 + A_1 $ to 
$K(n)$. To  interpret $D^b(K(n))$ as a perverse sheaf of categories one considers a  locally constant sheaf of categories over a graph $\Gamma$ shown on the picture below, the picture  encodes also the data about the sheaf, in particular  $p_1$, $p_2$ denote spherical functors (see \cite{AHH}, \cite{AKO}):
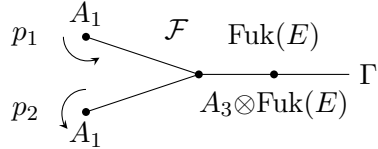
\begin{center}
\begin{tikzpicture}[scale=1]

\draw (0,0) -- (2,0) node[right]{$\Gamma$};
\draw (0,0) -- (-1.5,0.5)
(0,0) -- (-1.5,-0.5)
;

\draw (1,0.5) node{Fuk$(E)$}
(1,-0.4) node{$A_3 \otimes$Fuk$(E)$}
(-1.5,0.8) node{$A_1$}
(-1.5,-0.8) node{$A_1$}
(-2.3,-0.5) node{$p_2$}
(-2.3,0.5) node{$p_1$}
;
\draw [shift={(-1.5,-0.5)}, -stealth]
(0,0.3) arc (90:220:0.32);
;
\draw [shift={(-1.5,0.5)}, -stealth,]
(-0.3,0) arc (180:300:0.32);
;
\fill (0,0) circle (1.5pt)
(1,0) circle (1.5pt)
(-1.5,0.5) circle (1.5pt)
(-1.5,-0.5) circle (1.5pt)
;

\draw (-0.3,0.6) node{$\mcf$};

\end{tikzpicture}
\end{center}
 The category of global sections $H^0(\Gamma, \mcf)$  of the  sheaf $\mcf$  is the same as   $D^b(K(n)).$
The surgeries are recorded by the changes of the spherical functors $p_1$, $p_2$.

The category    $D^b(K(4))$ can be interpreted also as part of the LG model of $\PP^3$,   ${\CC^*}^3$, $w=x+y+z+\frac{1}{xyz}$:
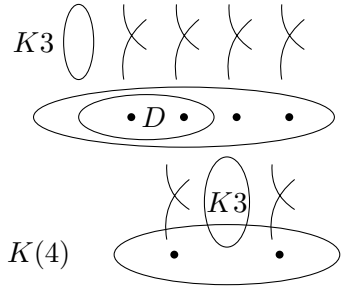
\begin{center}
\begin{tikzpicture}[scale=1]

\draw[xshift=-0.7cm] (-0.1,0.2) to[bend left=30] (0.2,1)(-0.1,1.2) to[bend right=30] (0.2,0.6);
\draw[xshift=0cm] (-0.1,0.2) to[bend left=30] (0.2,1)(-0.1,1.2) to[bend right=30] (0.2,0.6);
\draw[xshift=0.7cm] (-0.1,0.2) to[bend left=30] (0.2,1)(-0.1,1.2) to[bend right=30] (0.2,0.6);
\draw[xshift=1.4cm] (-0.1,0.2) to[bend left=30] (0.2,1)(-0.1,1.2) to[bend right=30] (0.2,0.6);

\begin{scope}[shift={(0,-0.3)}]
\fill (1.4,-0) circle (1.5pt);
\fill (0.7,-0) circle (1.5pt);
\fill (0,-0) circle (1.5pt);
\fill (-0.7,-0) circle (1.5pt);
\draw (-0.4,-0) node{$D$} (-0.5,0) ellipse (0.9 and 0.3);
\draw (-1.4,1) ellipse (0.2 and 0.5) (-2,1) node{$K3$};
\draw ellipse (2 and 0.4);
\end{scope}

\end{tikzpicture}
\end{center}

\begin{center}
\begin{tikzpicture}[scale=1]

\draw[xshift=-0.7cm] (-0.1,0.2) to[bend left=30] (0.2,1)(-0.1,1.2) to[bend right=30] (0.2,0.6);
\draw[xshift=0.7cm] (-0.1,0.2) to[bend left=30] (0.2,1)(-0.1,1.2) to[bend right=30] (0.2,0.6);
\draw (0,0.7) ellipse (0.3 and 0.6)  node{$K3$};

\draw ellipse (1.5 and 0.4);
\fill (0.7,-0) circle (1.5pt);
\fill (-0.7,-0) circle (1.5pt);

\draw (-2.5,0) node{$K(4)$};

\end{tikzpicture}
\end{center}
We make a surgery on the fiber -  a $K3$ surface: 
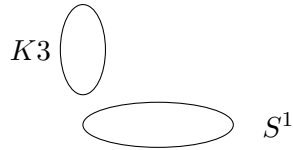
\begin{center}
\begin{tikzpicture}[scale=1]

\draw (0,1) ellipse (0.3 and 0.6)  (-0.7,1) node{$K3$};

\draw (0,0) arc (180:-180:1 and 0.3) (2.6,0) node{$S^1$};

\end{tikzpicture}
\end{center}
This surgery amounts to change from
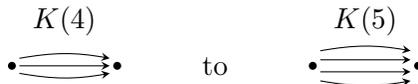
\begin{center}
\begin{tikzpicture}[scale=1]

\draw (0,0) node{to};

\begin{scope}[shift={(-2,0)}]
\fill (0.7,-0) circle (1.5pt);
\fill (-0.7,-0) circle (1.5pt);
\draw (0,0.6) node{$K(4)$};
\draw[-stealth] (-0.6,0) to (0.6,0);
\draw[-stealth] (-0.6,0.1) to[bend left=10] (0.6,0.1);
\draw[-stealth,yscale=-1] (-0.6,0.1) to[bend left=10] (0.6,0.1);
\end{scope}

\begin{scope}[shift={(2,0)}]
\fill (0.7,-0) circle (1.5pt);
\fill (-0.7,-0) circle (1.5pt);
\draw (0,0.6) node{$K(5)$};
\draw[-stealth] (-0.6,0.08) to (0.6,0.08);
\draw[-stealth, yscale=-1] (-0.6,0.08) to (0.6,0.08);
\draw[-stealth] (-0.6,0.2) to[bend left=10] (0.6,0.2);
\draw[-stealth,yscale=-1] (-0.6,0.2) to[bend left=10] (0.6,0.2);
\end{scope}

\end{tikzpicture}
\end{center}
The  Landau Ginzburg models with $K3$ surfaces in the fibers  can be interpreted as  perverse sheafs of categories, encoded in the following picture - see \cite{AHH}:
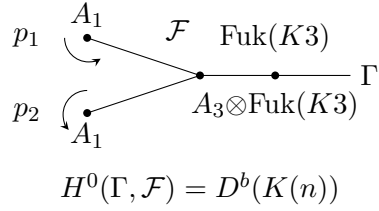
\begin{center}
\begin{tikzpicture}[scale=1]
\draw (0,0) -- (2,0) node[right]{$\Gamma$};
\draw (0,0) -- (-1.5,0.5)
(0,0) -- (-1.5,-0.5)
;

\draw (1,0.5) node{Fuk$(K3)$}
(1,-0.4) node{$A_3 \otimes$Fuk$(K3)$}
(-1.5,0.8) node{$A_1$}
(-1.5,-0.8) node{$A_1$}
(-2.3,-0.5) node{$p_2$}
(-2.3,0.5) node{$p_1$}
;
\draw [shift={(-1.5,-0.5)}, -stealth]
(0,0.3) arc (90:220:0.32);
;
\draw [shift={(-1.5,0.5)}, -stealth,]
(-0.3,0) arc (180:300:0.32);
;
\fill (0,0) circle (1.5pt)
(1,0) circle (1.5pt)
(-1.5,0.5) circle (1.5pt)
(-1.5,-0.5) circle (1.5pt)
;

\draw (-0.3,0.6) node{$\mcf$};

\draw (0,-1.5) node{$H^0(\Gamma,\mcf)=D^b(K(n))$};
\end{tikzpicture}.
\end{center}

\begin{remark} 
The property of having a phase gap, which we require  in this paper to define the norm, can also be interpret as existence of a CY form with certain properties.

Namely let $Y$ be a LG model, $\Omega$  is a CY form on $Y$. Let $L$ be a Lagrangian s.t.
$\theta_ 1 \leq arg \Omega|_ {L} \leq \theta_2 $. Assume that there exits  a form $\beta$ on $Y$ s.t.

(1) $\beta = d \alpha$, ($\alpha$  is an $n-1$ form ),
(2) $Re \beta | _{L}  > 0.$
(3) $\alpha \to  0$ when $\omega \to  0$.

Then there are no stable lagrangians $L$ with
$\theta_ 1 \leq arg \Omega|_ {L} \leq \theta_2 $. 
In other words existence of such forms $\Omega$ and $\alpha$ lead to gaps in phases.
\end{remark}

One more direction for future research is holomorphic families of categories, in particular holomorphic families of Kronecker quivers. 

Holomorphic families of categories over $X$ with fiber $K(n)$ should be  defined by homomorphisms $\varphi_i: \mco(U_i) \to \op{HH}^0(D^b(K(n)))$ in the zero-th Hochschild cohomology of $D^b(K(n))$  where  $\{U_i\}$ is a covering of $X$ by open sets. We use the following picture for such a holomorphic  family of categories:

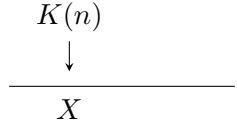
\begin{center}
\begin{tikzpicture}[scale=1]

\draw (-1.5,0) -- (1.5,0);
\draw (-0.7,-0.3) node{$X$};
\draw[-stealth] (-0.7,0.6) node[above]{$K(n)$}--(-0.7,0.2);

\end{tikzpicture}
\end{center}

The holomorphic sheaves of categories are enhanced by perverse sheaves of stability conditions - see \cite{KK} for defining 
morphisms and the gluing between the categories on intersecting opens that defines the sheaf.

The case of holomorphic family of $K(2)$ is the classical case of conic bundles:
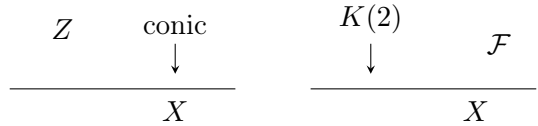
\begin{center}
\begin{tikzpicture}[scale=1]
\begin{scope}[shift={(-2,0)}]
\draw (-1.5,0) -- (1.5,0);
\draw (0.7,-0.3) node{$X$};
\draw (-0.8,0.8) node{$Z$};
\draw[-stealth] (0.7,0.6) node[above]{conic}--(0.7,0.2);
\end{scope}
\begin{scope}[shift={(2,0)}]
\draw (-1.5,0) -- (1.5,0);
\draw (0.7,-0.3) node{$X$};
\draw (1,0.6) node{$\mcf$};
\draw[-stealth] (-0.7,0.6) node[above]{$K(2)$}--(-0.7,0.2);
\end{scope}
\end{tikzpicture}
\end{center}

The global sections $H^0(X, \mcf)$ are $\db{Z}$. Similarly $H^0(X, \mcf)$ with $K(n)$ for $n\geq 3$ produces a new non-commutative variety. 

 Iterating the  procedure described above results in a family of categories over a family of categories.  Some questions addressing relations between   the norms of the  fibers and of the gobal sections follow: 
\begin{quest}
Under what condition $||\mcc||_{\varepsilon} \geqslant ||H^0(X,\mcf)||_{\varepsilon} ?$ (here $\mc C$ is the category in the fiber)
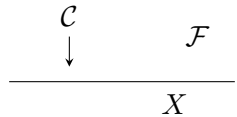
\begin{center}
\begin{tikzpicture}[scale=1]
\draw (-1.5,0) -- (1.5,0);
\draw (0.7,-0.3) node{$X$};
\draw (1,0.6) node{$\mcf$};
\draw[-stealth] (-0.7,0.6) node[above]{$\mcc$}--(-0.7,0.2);
\end{tikzpicture}
\end{center}
\end{quest}

\begin{quest} Let us consider  a tower of families of categories and each of the fiber categories $\mcc_i$ has non maximal (Recall the relation of $\norm{\cdot}^{\varepsilon}$ and $\norm{\cdot}_{\varepsilon}$ in Definition \ref{balls}) $\norm{\cdot}^{\varepsilon}$. Is it  true that if the category in the  combined fiber    has  a Rouquier dimension  \cite{BFK} equal to one then the  norm of this category  is    non-maximal $\norm{\cdot}^{\varepsilon}$ ?

\end{quest}

We summarise the proposed analogy in the table bellow.

\FloatBarrier
\begin{table}[ht]
\setlength{\tabcolsep}{1mm} 
\def\arraystretch{1.5}
\centering
\begin{tabular}{|>{\centering\arraybackslash}m{1.5in}|>{\centering\arraybackslash}m{2.8in}|}
  \hline
Galois theory & Norms
\\ \hline
\begin{tikzpicture}[scale=1]
\draw(0,2.4);
\draw (0,2) node(a){$X_2$};
\draw (0,1) node(b){$X_1$};
\draw (0,0) node(c){$X$};
\path [-stealth](a) edge node[right=0.3]{finite} (b) (b) edge node[right=0.3]{finite}(c);
\draw (0,-1) node{The sequence of }
(0,-1.5) node{finite coverings}
(0,-2) node{is finite};;

\end{tikzpicture}  & \begin{tikzpicture}[scale=1]
\draw(0,2.4);
\begin{scope}[shift={(-2,0)}]
\draw (0,2) node(a){$X_2$};
\draw (0,1) node(b){$X_1$};
\draw (0,0) node(c){$X$};
\path [-stealth](a) edge node[right=0.3]{$\mcc_2$ \quad $||\mcc_2||^\epsilon <$ max}  node[right=2]{} (b) (b) edge node[right=0.3]{$\mcc_1$ \quad $||\mcc_1||^\epsilon <$ max} node[right=2]{}    (c);
\end{scope}

\draw (2.6,1) node{$\left( \begin{matrix}\mcc_2 \\ \downarrow \\ \mcc_1  \end{matrix} \right)=\mcc$}; 
\draw (3.7,1) node{$ \begin{matrix} X_2 \\ \downarrow \\ X  \end{matrix} $};  \qquad 

\draw[shift={(0.4,0)}] (0,-1) node{Rouquier dim $(\mcc)=1$}
(0,-1.5) node{$\Downarrow$}
(0,-2) node{$||\mcc||^\epsilon<$ max};;
\end{tikzpicture}
\\ \hline
  \end{tabular}
\begin{center}
\begin{tikzpicture}[scale=1]

;

\end{tikzpicture}
\end{center}
\begin{tikzpicture}[overlay]

\end{tikzpicture}
\end{table}
\FloatBarrier

\begin{quest} Do we have a similar theory as the classical theory of conic bundles for sheaves of categories with fibers categories of   representations of Kronecker 
quivers or any other quiver category with a Rouquier dimension \cite{BFK} equal to one?

\end{quest}

In a certain way our norm can be seen as analogue of height function defined in \cite{BP1}.
We expect that some higher analogues of this norm  for higher Rouquier dimensions can be defined.
In fact in this paper we  only scratch the surface  proposing  a   possible  approach  to  ``noncommutative Galois theory''  - representing ``noncommutative manifolds'' (categories) as a 
sequence of perverse sheaves of categories and holomorphic families of categories.

It will be interesting to  study categories which can be represented  as  a tower of holomorphic families of categories with nonmaximal norms. 
One example  of such  category  is $D^b(\PP^1 \times ... \times\PP^1)$. 

\begin{quest} Characterise projective varieties $X$
whose derived categories $D^b(X)$  can be 
represented as tower of holomorphic families of categories with non-maximal norms.

1) Under what conditions are these projective varieties $X$ rational?
(It is rather clear that a nontrivial condition is needed since every hyperelliptic curve can be seen as such a tower. )
 
 2) Can the existence of   tower of holomorphic families of categories
 be represented in terms of modular forms?

 3) We conjecture that a  unirational variety $U$ can be represented as a tower

\bd  Z &\rTo^{F_1} & Z_1 & \rTo^{F_1} Z_2 & \cdots & \rTo & U \ed
where $\norm{D^b(F_i)}^{\varepsilon}<\pi (1-
\varepsilon)$ (see \cite{BT}), and $Z$ is rational.  And another question is, does there exists a cohomology theory, which determines, if $U$ is rational.

\end{quest}

  In the end  we put a  question coming  from the interplay between towers of sheaves of categories and stability conditions. 
 Let $K(N_1,N_2)$ be a category obtained from a family  where the base is a category of representations of the 
 Kronecker quiver $K(N_1)$   and the fibers are the   category   of representations of the 
 Kronecker quiver  $K(N_2)$. Similarly we have $K(N_1, N_2,..., N_l)$ 
 denoting extensions of extensions.

\begin{quest} Is the moduli space of stability condition  of $K(N_1, N_2,.., N_l)$ a bundle over 
 Hilbert modular variety? 
   \end{quest}
 It would  be interesting to investigate  the connection of the geometry of these Hilbert modular surfaces with the norms we have defined as well as new modular identities coming from 
 wall-crossing formulas.

\let\oldaddcontentsline\addcontentsline
\renewcommand{\addcontentsline}[3]{}

\let\addcontentsline\oldaddcontentsline

\end{document}